\def\no{\if01}
\def\iftwelvept{\no}

\def\ifusepdf{\no}
\def\ifpsfont{\no}

\iftwelvept
\documentclass[leqno,12pt]{amsart}
\else
\documentclass[leqno,11pt]{amsart}
\fi
\usepackage{amssymb}
\usepackage{amscd}
\usepackage{latexsym}
\usepackage{verbatim}
\usepackage[all]{xy}

\ifusepdf
\usepackage{hyperref}
\else\fi
\ifpsfont
\usepackage[T1]{fontenc}
\usepackage{times}
\else\fi

\iftwelvept
\else\fi

\theoremstyle{plain}
\newtheorem{Theorem}{Theorem}[section]

\newtheorem{Proposition}[Theorem]{Proposition}
\newtheorem{Lemma}[Theorem]{Lemma}
\newtheorem{Corollary}[Theorem]{Corollary}

\newtheorem{Claim}{Claim}[Theorem]

\theoremstyle{definition}

\newtheorem{Definition}[Theorem]{Definition}
\newtheorem{Remark}[Theorem]{Remark}
\newtheorem{Example}[Theorem]{Example}

\newcommand{\AI}{A_\infty}
\newcommand{\bA}{\bar{A}}
\newcommand{\Ph}{\partial_{\textup{Hoch}}}
\newcommand{\PH}{\partial^{\textup{Hoch}}}
\newcommand{\Art}{\operatorname{Art}_k}
\newcommand{\Set}{\operatorname{Sets}}
\newcommand{\MC}{\operatorname{MC}}
\newcommand{\Def}{\mathsf{Spf}}
\newcommand{\DEF}{\mathsf{Spf}}

\newcommand{\ZZ}{{\mathbb{Z}}}

\newcommand{\CC}{{\mathbb{C}}}

\newcommand{\NN}{{\mathbb{N}}}

\newcommand{\EE}{\mathbb{E}}

\newcommand{\Gr}{\mathbb{F}}

\newcommand{\OO}{{\mathcal{O}}}

\newcommand{\Hom}{\operatorname{Hom}}
\newcommand{\Ext}{\operatorname{Ext}}

\newcommand{\Ker}{\operatorname{Ker}}
\newcommand{\Coker}{\operatorname{Coker}}

\newcommand{\Mod}{\operatorname{Mod}}

\newcommand{\Map}{\operatorname{Map}}

\newcommand{\Fun}{\operatorname{Fun}}
\newcommand{\Alg}{\operatorname{Alg}}

\newcommand{\End}{\operatorname{End}}

\newcommand{\CAlg}{\operatorname{CAlg}}
\newcommand{\dGr}{\mathsf{Gr}}

\newcommand{\Coder}{\operatorname{Coder}}

\newcommand{\Aut}{\operatorname{Aut}}

\newcommand{\Sym}{\operatorname{Sym}}
\newcommand{\Proof}{{\sl Proof.}\quad}
\newcommand{\QED}{{\unskip\nobreak\hfil\penalty50\quad\null\nobreak\hfil
{$\Box$}\parfillskip0pt\finalhyphendemerits0\par\medskip}}

\begin{document}

\title{Period mappings for noncommutative algebras}

\author{Isamu Iwanari}

\address{Mathemathtical Institute, Tohoku University, Japan}

\email{iwanari@math.tohoku.ac.jp}

%\date{version Apr/2016}

\maketitle

\setcounter{tocdepth}{1}
\tableofcontents

\section{Introduction}

Differential graded (dg) algebras naturally appear in
the center of noncommutative geometry.
One of main approaches to noncommutative geometry is
to regard a pretriangulated dg category (or some kind of
stable $(\infty,1)$-categories) as a space
and think that two spaces coincide if two categories are equivalent.
This approach unifies 
many branches of mathematics such as algebraic geometry, symplectic geometry, representation theory and mathematical physics.
In many interesting cases, such a category of interest
admits a single compact generator,
and the dg category is quasi-equivalent to the dg category of dg modules
over a (not necessarily commutative) dg algebra.
Thus, dg algebras can be viewed as incarnations of noncommutative spaces.
For example, by \cite{BV} the dg category of quasi-coherent complexes
on a quasi-compact and separated scheme admits a single compact generator.

The purpose of this paper is to construct a period mapping for
deformations of a dg algebra.
We generalize Griffiths' period mapping for deformations
of algebraic varieties to noncommutative setting.
Before proceeding to describe formulations and results of this paper,
we would like to briefly review the Griffiths' construction of the period
mapping.
Let $X_0$ be a complex smooth projective variety (more generally, a complex compact
K\"ahler manifold).
Let $f:X\to S$ be a deformation of $X_0$ such that $X_0$
is the fiber over a point $0\in S$, and $S$ is a base complex manifold.
For each fiber $X_s$ over $s\in S$, consider the Hodge filtration
$F^rH^i(X_s,\CC)\subset H^i(X_s,\CC)$ on the singular cohomology.
The derived pushforward $R^if_*\CC_X$ of the
constant sheaf $\CC_X$ is a local system.
We assume that $S$ is small so that $S$ is contractible.
There is the natural identification $H^i(X_s,\CC)\simeq H^i(X_0,\CC)$,
and one has the pair
\[
(F^rH^i(X_s,\CC), H^i(X_s,\CC)\simeq  H^i(X_0,\CC)).
\]
Therefore, to any point $s$ one can associate the subspace $F^rH^i(X_s,\CC)$
in $V:=H^i(X_0,\CC)$ via $H^i(X_s,\CC)\simeq  H^i(X_0,\CC)$.
It gives rise to the classifying map called the (infinitesimal)
period mapping
\[
\mathbf{P}:S \to \textup{Grass}(V,n)
\]
where $n=\dim F^rH^i(X_0,\CC)$, and
$\textup{Grass}(V,n)$ is the Grassmannian which parametrizes
the $n$-dimensional subspaces of $V$.
%If we consider all subspaces $\{F^rH^i(X_s,\CC)\}_{r\ge 0}$, then we have the period mapping to the flag manifold.

In order to generalize $\mathbf{P}$, it is necessary to consider a
noncommutative analogue of Hodge filtrations.
Remember that by 
Hochschild-Kostant-Rosenberg theorem the periodic cyclic homology $HP_*(X)$ of a smooth scheme $X$
over a field of characteristic zero
can be described in terms of the algebraic de Rham cohomology of $X$.
Namely, $HP_i(X)=\prod_{n:odd}H_{dR}^n(X)$ for odd $i$, and $HP_i(X)=\prod_{n:even}H^n_{dR}(X)$ for even $i$.
The image of the 
natural map $NH_*(X)\to HP_*(X)$ from the negative cyclic homology $HN_*(X)$
determines
the subspace which can be experessed as a product of Hodge filtrations (see
for instance \cite{Wei}, Remark~\ref{classHodge}).
Note that cyclic homology theories can be defined for
more general objects such as dg algebras, dg categories, etc.
Let $A$ be a dg algebra over a field of characteristic zero.
Let $CC^-_\bullet(A)$ and $CC^{per}_\bullet(A)$
denote its negative cyclic complex and its periodic cyclic complex
respectively. It is for this reason that
it is natural to think of the natural map of complexes
$CC^-_\bullet(A)\to CC^{per}_\bullet(A)$
as a noncommutative analogue of Hodge complex.
Hence, intuitively speaking,
a period mapping should carry a deformation $\tilde{A}$ of $A$
to the deformation $CC^-_\bullet(\tilde{A})\to CC^{per}_\bullet(\tilde{A})$
of $CC^-_\bullet(A)\to CC^{per}_\bullet(A)$ endowed with a trivialization
of $CC^{per}_\bullet(\tilde{A})$.

We employ the deformation theory by means of dg Lie algebras.
The basic idea of so-called derived deformation theory is
that 
any reasonable deformation problem in characteristic zero
is controled by a dg Lie algebra.
To a dg Lie algebra $L$ over a field $k$ of characteristic zero,
one can associate
a functor $\Def_{L}:\Art\to \textup{Sets}$, from the category of artin local $k$-algebras with residue field $k$
to the category of sets, which assigns to any artin local algebra $R$ with the maximal ideal $m_R$
to the set of equivalence classes of Maurer-Cartan elements of $L\otimes m_R$, see Section~\ref{premoduli}
 (see \cite{Hin0}, \cite{DAGX} for a homotopical foundation of this approach).
Deformation theory via dg Lie algebras fits in with derived moduli theory in derived algebraic
geometry (see e.g., the survey \cite{To}). In fact, the Chevalley-Eilenberg cochain complex of a dg Lie
algebra plays the role of the ring of functions on a formal neighborhood of
a point of a derived moduli space.
Another remarkable
 advantage is that it allows one to study deformation problems
using homological methods of dg Lie algebras.
That is to say, once one knows that a deformation problem
$\Art\to \textup{Sets}$ which assigns $R$ to the isomorphism classes
of deformations to $R$ is of the form $\Def_L$,
one can use homological algebra concerning $L$.
It has been fruitful.
For example, the famous applications can be found in Goldman-Millson's
local study on certain moduli spaces \cite{GM}
and Kontsevich's deformation quantization of a Poisson
manifold \cite{Kon}.

One of main dg Lie algebras in this paper is
the dg Lie algebra $C^\bullet(A)[1]$, that is
the shifted Hochschild cochain complex
 of a dg algebra $A$ over the field $k$.
It controls curved
$A_\infty$-deformations of $A$. Namely, the
functor $\Def_{C^\bullet(A)[1]}:\Art\to \textup{Sets}$
associated to $C^\bullet(A)[1]$
can be identified with the functor which assigns to $R$ the set
of isomorphism classes of curved $A_\infty$-deformations of $A$ to $R$,
that is, curved $A_\infty$-algebras $A\otimes_kR$ over $R$ whose
reduction $A\otimes_kR/m_R$ are identified with $A$
(see Section~\ref{premoduli} for details).
Let $k[t]$ be the dg algebra (with zero differential)
generated by $t$ of cohomological degree $2$.
The complex $CC^-_\bullet(A)$ naturally
comes equipped with the action of $k[t]$,
and $CC^{per}_\bullet(A)$ comes equipped with the action of $k[t,t^{-1}]$.
We then conisder a pair
$(W, \phi:CC_\bullet^{per}(A)\otimes_kR\simeq W\otimes_{R[t]}R[t,t^{-1}])$
such that $W$ is a deformation of dg $k[t]$-module $CC^-_\bullet(A)$ to $R$,
and $\phi$ is an isomorphism of dg $R[t,t^{-1}]$-modules.
There is a dg Lie algebra $\mathbb{F}$ such that $\Def_{\mathbb{F}}(R)$
can be naturally identified with the set of
isomorphism classes of 
such pairs.
We can think of $\Def_{\mathbb{F}}$ as a formal
Sato Grassmannian generalized to the complex level (see Section~\ref{peri}).

From a naive point of view, our main construction of a period mapping
may be described as follows (see Section~\ref{construction}, Theorem~\ref{moduli} for details):

\begin{Theorem}
\label{main1}
We construct an $L_\infty$-morphism
$\mathcal{P}:C^\bullet(A)[1]\to \mathbb{F}$
such that the induced morphism of functors
$\Def_{\mathcal{P}}:\Def_{C^\bullet(A)[1]}\to \Def_{\mathbb{F}}$
has the following modular interpretation:
If one identifies $\alpha\in \Def_{C^\bullet(A)[1]}(R)$
with a curved $A_\infty$-deformation $\tilde{A}_\alpha$,
then $\Def_{\mathcal{P}}$ sends $\tilde{A}_\alpha$
to the pair
$(CC_\bullet^-(\tilde{A}_\alpha),CC_\bullet^{per}(A)\otimes_kR\simeq CC_\bullet^-(\tilde{A}_\alpha)\otimes_{R[t]}R[t,t^{-1}])$ in $\Def_{\mathbb{F}}(R)$.
(See Section~\ref{mappingsp} for $L_\infty$-morphisms)
\end{Theorem}
Note that the pair
\[
(CC_\bullet^-(\tilde{A}_\alpha),CC_\bullet^{per}(A)\otimes_kR\simeq CC_\bullet^-(\tilde{A}_\alpha)\otimes_{R[t]}R[t,t^{-1}])
\]
associated to $\tilde{A}_\alpha$ in Theorem~\ref{main1} is a
counterpart of
$(F^rH^i(X_s,\CC), H^i(X_s,\CC)\simeq  H^i(X_0,\CC))$
appeared in the construction of the period mapping for deformations of $X_0$.
For this reason, we shall refer to 
$\Def_{\mathcal{P}}:\Def_{C^\bullet(A)[1]}\to \Def_{\mathbb{F}}$
as the period mapping for $A$ (one may think of the $L_\infty$-morphism
$\mathcal{P}$ as
a Lie algebra theoretic realization of the period mapping).
The construction of the $L_\infty$-morphism
$\mathcal{P}:C^\bullet(A)[1]\to \mathbb{F}$ (and the period mapping)
is carried out in several steps.
One of the key inputs is
a modular interpretation of the Lie algebra
action of the shifted Hochschild cochain complex
$C^\bullet(A)[1]$ on Hochschild chain complex $C_\bullet(A)$.
Note that if a dg algebra $A$ is deformed, then it induces
a deformation of Hochschild chain complex $C_\bullet(A)$.
Since deformations of the complex $C_\bullet(A)$ is controled
by the endomorphism dg Lie algebra $\End(C_\bullet(A))$,
there should be a corresponding morphism $C^\bullet(A)[1]\to \End(C_\bullet(A))$
of dg Lie algebras. We find that it is given by
the action $L:C^\bullet(A)[1]\to \End(C_\bullet(A))$ of $C^\bullet(A)[1]$ on $C_\bullet(A)$ which was studied
by Tamarkin-Tsygan \cite{TT} generalizing the calculus structure on
the pair $(HH^*(A),HH_*(A))$ given by Daletski-Gelfand--Tsygan,
to the level of complexes (see Section~\ref{LM}).
(It is crucial for our construction to work with complexes.)
Also, it allows us to make use of homological algebra of homotopy calculus
operad (see Section~\ref{construction}).
Another input is a result on a nullity (Propsoition~\ref{triviality}).
It says that the morphism $L((t)):C^\bullet(A)[1]\to \End_{k[t,t^{-1}]}(CC_\bullet^{per}(A))$ of dg Lie algebras induced by $L$ 
(see Section~\ref{cylicify}) is
null-homotopic. This fact plays a role analogous to Ehresmann fibration
theorem.

\vspace{2mm}

Although we have discussed about general dg algebras,
we would like to focus on the important case of interest where
the dg algebra $A$ is smooth and proper
(see Section~\ref{smoothproper}).
In that case, the degeneration of an analogue of Hodge-to-de Rham
spectral sequence on cyclic homology theories was conjectured
by Kontsevich-Soibelman and was proved by Kaledin (see Section~\ref{Kaledin}).
As in the commutative world, the degeneration is pivotal.
We obtain (see Section~\ref{concludeSP}):

\begin{Theorem}
\label{main2}
Suppose that $A$ is smooth and proper.
Then the dg Lie algebra $\mathbb{F}$ is equivalent to
the dg Lie algebra
$\End_{k[t,t^{-1}]}(HH_*(A)((t)))/\End_{k[t]}(HH_*(A)[[t]])[-1]$
equipped with zero differential and zero bracket, and the $L_\infty$-morphism
in Theorem~\ref{main1} is
\[\mathcal{P}:C^\bullet(A)[1]\to \End_{k[t,t^{-1}]}(HH_*(A)((t)))/\End_{k[t]}(HH_*(A)[[t]])[-1].
\]
Here $HH_*(A)$ is the graded vector space of Hochschild homology,
and $HH_*(A)[[t]]$ is the graded vector space $HH_*(A)\otimes_kk[t]$
(the cohomological degree of $t$ is $2$).
The graded vector space
$HH_*((t))$ is $HH_*(A)\otimes_{k}k[t,t^{-1}]$.
\end{Theorem}

By using the properties on our period mapping and its codomain
(so-called the period domain), one can deduce the property of the domain of the period mapping.
Indeed, applying the period mapping as an $L_\infty$-morphism,
we prove (see Section~\ref{BTTsection}):

\begin{Theorem}
\label{main3}
Suppose that the dg algebra $A$ is smooth, proper and Calabi-Yau.
Then the dg Lie algebra $C^\bullet(A)[1]$ is quasi-abelian.
Namely, it is quasi-isomorphic to an abelian dg Lie algebra.
In particular, a curved $A_\infty$-deformation of $A$ is unobstructed.
\end{Theorem}

Theorem~\ref{main3} is a Bogomolov-Tian-Todolov
theorem for $C^\bullet(A)[1]$.
It is also applicable to a dg category which is quasi-equivalent to the dg category of dg modules over $A$. This unobstructedness was formulated
in Katzarkov-Kontsevich-Pantev \cite[4.4.1]{KKP} with the outline of
an approach which uses
an action of the framed little disks operad in the Calabi-Yau case.
The proof of Theorem~\ref{main3}
in this paper is different from the approach in {\it loc. cit}.,
for example, we do not employ framed little disks operad.
Nevertheless, it would be interesting to compare them.

We obtain the following infinitesimal Torelli theroem as a direct but interesting consequence from our moduli-theoretic construction and Calabi-Yau property:

\begin{Theorem}
\label{main4}
Suppose that $A$ is smooth, proper and Calabi-Yau.
The period mapping
$\Def_{\mathcal{P}}:\Def_{C^\bullet(A)[1]}\to \Def_{\mathbb{F}}$,
which carries a curved $A_\infty$-deformation
$\tilde{A}$ to the associated pair
\[
(CC_\bullet^-(\tilde{A}),CC_\bullet^{per}(A)\otimes_kR\simeq CC_\bullet^-(\tilde{A})\otimes_{R[t]}R[t,t^{-1}]),
\]
is a monomorphism.
Here as in Theorem~\ref{main1}
we implicitly identifies an element in $\Def_{C^\bullet(A)[1]}(R)$
with a curved $A_\infty$-deformation.
\end{Theorem}

It might be worth considering the case when the dg algebra $A$ comes from
a (quasi-compact and separated) smooth scheme $X$ over the field $k$
(cf. Remark~\ref{comparisonscheme}).
The tangent space $\Def_{C^\bullet(A)[1]}(k[\epsilon]/(\epsilon^2))$
of $\Def_{C^\bullet(A)[1]}$ is isomorphic to the second Hochschild
cohomology $HH^2(A)$ as vector spaces.
According to Hochschild-Kostant-Rosenberg theorem, there is an
isomorphism
\[
\Def_{C^\bullet(A)[1]}(k[\epsilon]/(\epsilon^2))\simeq HH^2(A)\simeq H^0(X,\wedge^2T_X)\oplus H^1(X,T_X)\oplus H^2(X,\OO_X)
\]
where $T_X$ is the tangent bundle, and $\OO_X$ is the structure sheaf.
Informally speaking, the space $H^1(X,T_X)$ parametrizes
deformations of the scheme $X$ (in the commutative world),
and the subspace $H^0(X,T_X^{\wedge 2})\oplus H^2(X,\OO_X)$ is
involved with
(twisted) deformation quantizations of $X$, which we think of as the noncommutative part.
Thus, in that case, Theorem~\ref{main4} means that
the associated pair is deformed faithfully along deformations of the algebra not only to commutative directions but also to noncommutative directions.

The author would like to thank H. Minamoto and Sho Saito for interesting discussions related to the topics of this paper. He thanks A. Takahashi
for pointing out to him the paper \cite{KKP}.
The author is partly supported by Grant-in-Aid for Young Scientists JSPS.

\vspace{3mm}

\begin{comment}
This paper is organized as follows: In Section 2,
we review curved $A_\infty$-algebras, Hochschild (co)chain complexes,
cyclic complexes and fix our convention for them.
In Section 3, we study a modular
interpretation of the action $L$ of $C^\bullet(A)[1]$ on $C_\bullet(A)$.
For this,
the brief review on deformation theory from the viewpoint of dg Lie algebras
and the discussion on examples are also included.
In Section 4, we construct a period mapping (cf. Theorem~\ref{main1}).
In Section 5, we consider the smooth and proper case.
We apply the degeneration of Hodge-to-de Rham spectral sequence
to obtain Theorem~\ref{main2} (see Section~\ref{concludeSP}). We give applications in Section 6 and 7.
Using the period mapping we prove a noncommutative Bogomolov-Tian-Todolov theorem (cf. Theorem~\ref{generalBTT}, Corollary~\ref{BTT}). In Section 7, we prove the infinitesimal Torelli theorem for our period mapping.
\end{comment}

\section{$A_\infty$-algebra and Hochschild complex}

Our main interest lies in differential graded algebras and their deformations.
But, if we consider (explicit) deformations of differential
graded algebras,
curved $A_\infty$-algebras naturally appear as deformed objects.
Thus, in this Section
we review the basic definitions and facts about curved
$A_\infty$-algebras and Hochschild complexes which we need later.
We also recall cyclic complexes of $A_\infty$-algebras.

\vspace{2mm}

{\it Convention.}
When we consider $\ZZ$-graded modules,
differential graded modules, i.e., complexes, etc.,
we will use the cohomological grading which are denoted by upper indices.
But, the homological grading is familiar
to homology theories such as Hochschild homology, cyclic homology.
Thus, when we treat Hochschild, cyclic homology theories etc.,
we use the homological grading which are denoted by lower indices.

\subsection{}
Let $k$ be a unital commutative ring. In this paper we usually treat
the case when $k$ is either a field of characteristic zero or
an artin local ring over a field of characteristic zero.
Let $A$ be a $\ZZ$-graded $k$-module.
Unless stated otherwise,
we assume that
each $k$-module $A^n$ is a free $k$-module.
If $a\in A$ is a homogeneous element, we shall denote by $|a|$ the (cohomological) degree
of $a$ (any element $a$ appeared in $|a|$ will be implicitly
assumed to be homogeneous).
We put $TA:=\oplus_{n\ge 0}A^{\otimes n}$
where $A^{\otimes n}$ denotes the $n$-fold tensor product.
We shall use the standard symmetric monoidal structure of the category of $\ZZ$-graded
$k$-modules. In particular, the symmetric structure induces
$a\otimes b\mapsto (-1)^{|a||b|}b\otimes a$.
We often write $\otimes$ for the tensor product $\otimes_k$ over the
base $k$.
The graded $k$-module $TA$ has a (graded) coalgebra structure
given by a comultiplication $\Delta:TA\to TA\otimes TA$ where
\[
\Delta(a_1,\ldots,a_n)=\sum_{i=0}^{n}(a_1,\ldots,a_i)\otimes(a_{i+1},\ldots,a_n),
\]
and a counit $TA\to A^{\otimes 0}=k$ determined by the projection.
Here we write $(a_1,\ldots,a_n)$ for $a_1\otimes a_2\otimes\ldots\otimes a_n$.

A $k$-linear map $f:TA\to TA$ is said to be a coderivation
if two maps $(1\otimes f+f\otimes 1)\circ\Delta, \Delta \circ f:TA\rightrightarrows TA\otimes TA$ coincide.
Let us denote by $\Coder(TA)$ the graded $k$-module of coderivations from
$TA$ to $TA$, where the grading is defined as the grading on maps of graded modules. The composition with the projection $TA\to A$ induces an isomorphism
\[
\Coder(TA)\to \Hom_k(TA,A)
\]
where the right-hand side is the hom space of $k$-linear maps,
and the inverse is given by $\Sigma_{i=0}^{n}\Sigma_{j=0}^{n-j}1^{\otimes j}\otimes f_i\otimes 1^{\otimes n-i-j}$ for $\{f_i\}_{i\ge 0}\in \prod_{i\ge 0}\Hom(A^{\otimes i},A)\simeq \Hom(TA,A)$ (cf. \cite[Proposition 1.2]{GJ}).

\subsection{}
\label{Hochscochaindef}
For a graded $k$-module $A$, the suspension $s A=A[1]$ is the same $k$-module
endowed with the shift grading $(A[1])^i=A^{i+1}$.
Put $BA=T(A[1])$.
A curved $A_\infty$-structure on a graded $k$-module $A$ is a coderivation $b:BA\to BA$ of degree $1$
such that $b^2=b\circ b=0$.
Let $b_i$ denote the composite
$(A[1])^{\otimes i}\hookrightarrow BA\stackrel{b}{\to} BA\to A[1]$
with the projection.
A curved $A_\infty$-structure $b$ on $A$ is called an $A_\infty$-structure
if $b_0=0$.
A curved $A_\infty$-structure $b$ on $A$ is called a differential graded (dg for  short)
algebra structure if $b_0=0$ and $b_i=0$ for $i>2$.
A curved $A_\infty$-algebra (resp. $A_\infty$-algebra, dg algebra) is a pair $(A,b)$ with $b$ a curved $A_\infty$-structure (resp. an $A_\infty$-structure,
dg algebra) and we often abuse notation
by writing $A$ for $(A,b)$ if no confusion is likely to arise.

Let $s:A\to A[1]$ be the ``identity map'' of degree $-1$
which identifies $A^i$ with $(A[1])^{i-1}$.
The map $s$ induces $A^{\otimes n}\simeq (A[1])^{\otimes n}$.
Therefore we have a $k$-linear map $m_n:A^{\otimes n}\to A$
such that the following diagram
\[
\xymatrix{
A^{\otimes n}  \ar[r]^{m_n} \ar[d]_{s^{\otimes n}} & A  \\
(A[1])^{\otimes n} \ar[r]^{b_n} & A[1] \ar[u]^{s^{-1}}.
}
\]
commutes.
We adopt the sign convention \cite{GJ} which especially
implies
\[
(f\otimes g)(a\otimes a')=(-1)^{|a||g|}f(a)\otimes g(a')
\]
where $a, a'\in A$ and $f$ and $g$ are $k$-linear maps of degree $|f|$
and $|g|$ respectively.
To be precise, we write $s^{\otimes n}$ for $(s\otimes 1^{\otimes n-1})(1\otimes s \otimes 1^{n-2})\cdots (1^{\otimes n-1}\otimes s)$.
(In the literature the author adopts $s^{\otimes n}=(1^{\otimes n-1}\otimes s)(1^{\otimes n-2} \otimes s\otimes  1)\cdots (s\otimes 1^{\otimes n-1})$.)
Put $[a_1|\ldots |a_n]=sa_1\otimes\ldots \otimes sa_n$.
Applying this sign rule to $s$ we have
$b_n[a_1|\ldots|a_n]=(-1)^{\Sigma_{i=1}^{n-1}(n-i)|a_i|}m_n(a_1,\ldots,a_n)$.
A curved $A_\infty$-algebra $(A,b)$
is said to be unital if there is an element $1_A$ of degree zero (called a unit)
such that (i) $m_1(1_A)=0$, (ii) $m_2[a|1_A]=m_2[1_A|a]=a$ for any $a\in A$,
and (iii) for $n\ge 3$, $m_{n}(a_1,\ldots,a_n)=0$ for one of $a_i$ equals $1_A$.

\subsection{}
An $A_\infty$-morphism $f:(A,b)\to (A',b')$ between curved $A_\infty$-algebras
is a differential graded (=dg) coalgebra map $f:(BA,b)\to (BA',b')$, for which
we often write an $A_\infty$-morphism $f:A\to A'$.
In the literature authors often impose the condition $f(k)=k\subset BA'$.
But we emphasize that when we treat curved $A_\infty$-algebras and deformation theory
it is natural to drop this condition (see also \cite{Lo}).
Consider the natural projection $BA'\to A'[1]$. 
Then
the composition with the projection $BA'\to A'[1]$ gives rise to an isomorphism
\[
\Hom_{\textup{coalg}}(BA,BA')\to \Hom_{k}(BA,A'[1])
\]
where the left-hand side indicates the hom set of graded
coalgebra maps and right-hand side indicates hom set of $k$-linear maps.
The inverse is given by
$BA \stackrel{\Delta_n}{\to}BA^{\otimes n}\stackrel{f^{\otimes n}}{\to} (A'[1])^{\otimes n}$
for $f:BA\to A'[1]$ where $\Delta_n$ is the $(n-1)$-fold iteration of comultiplication for each $n\ge 0$.
Therefore, a dg coalgebra map $f:BA\to BA'$
is determined by the family of graded $k$-linear maps $\{f_n:(A[1])^{\otimes n} \hookrightarrow BA\to BA'\to A'[1]\}_{n\ge 0}$ that satisfies a certain relation between $b_i$'s and $f_j$'s
corresponding to the compatibility with respect to differentials.
We shall refer to $f_n$ as the $n$-th component of $f$.
If $f_1$ induces an isomorphism $A[1]\to A'[1]$ of graded $k$-modules,
we say that $f$ is an $A_\infty$-isomorphism.
Suppose that $(A,b)$ is an $A_\infty$-algebra.
Then since the condition $b_0=0$ implies $b_1^2=m_1^2=0$, we can define the cohomology $H^*(A,b_1)$ of
$(A,b_1)$ which we regard as a graded $k$-module.
An $A_{\infty}$-morphism $f:(A,b)\to (A',b')$ of two $A_\infty$-algebras
is said to be quasi-isomorphism if $f$ induces
an isomorphism $H^*(A,b_1)\to H^*(A',b'_1)$ of graded $k$-modules.
For two unital curved $A_\infty$-algebras $A$ and $A'$, and an $A_\infty$-morphism $f:A\to A'$,
we say that $f$ is unital if $f(1_A)=1_{A'}$ and $f_n[a_1|\ldots|a_n]=0$
if one of $a_i$ equals $1_A$.

\subsection{}
We define a Hochschild cochain complex of an $A_\infty$-algebra.
Let $(A,b)$ be a unital $A_\infty$-algebra.
The graded $k$-module $\Coder(BA)$ has a differential given by
$\PH(f):=[b,f]=b\circ f-(-1)^{|f|}f\circ b$ such that $\PH\circ \PH=0$ .
We shall refer to 
\[
C^\bullet(A)[1]:=\Coder(BA)\simeq \Hom_k(BA,A[1])\simeq \prod_{n\ge 0}\Hom_k((A[1])^{\otimes n},A[1])
\]
as the shifted Hochschild cochain complex of $A$. Its cohomology $H^*(C^\bullet(A)[1])$
is called the shifted Hochschild cohomology of $A$.
We easily see that the bracket $[-,-]_G$ determined by the graded commutator
and the differential $\PH$ exhibit $C^\bullet(A)[1]$ as a dg Lie algebra.
For the definition of dg Lie algebra, see e.g. \cite[I. 3.5]{I}.
The bracket $[-,-]_G$ is call the Gerstenhaber bracket,
and by an abuse of notation we usually write $[-,-]$ for $[-,-]_G$.
Let us regard $f\in C^\bullet(A)[1]$ as a family of $k$-linear maps
$\{f_n\}_{n\ge0}\in \prod_{n\ge 0}\Hom_k((A[1])^{\otimes n},A[1])$
where $f_n:A[1])^{\otimes n}\to A[1]$. We say that $f$ is normalized
if $f_n[a_1|\ldots|a_n]=0$ if one of $a_i$ equals $1_A$ for $n\ge 1$.
The normalized elements constitutes a dg Lie subalgebra $\overline{C}^\bullet(A)[1]$ of $C^\bullet(A)[1]$ (we can easily check this by the definition of $\partial$
and the bracket).
Moreover, by \cite[Theorem 4.4]{Laz} there is a deformation retract $\overline{C}^\bullet(A)[1] \hookrightarrow C^\bullet(A)[1]\twoheadrightarrow \overline{C}^\bullet(A)[1]$, which induces quasi-isomorphisms.
If $A'$ is a dg algebra that is quasi-isomorphic to the $A_\infty$-algebra
$A$,
the complex $C^\bullet(A)[1]$ is quasi-isomorphic (via zig-zag) to $C^\bullet(A')[1]$ which is the usual (shifted)
Hochschild cochain complex, that is quasi-isomorphic to
the derived Hom complex $\textup{RHom}_{A'-\textup{bimodule}}(A',A')[1]$ (see \cite[Section 4]{Laz}).

\subsection{}
\label{Hochschildchain}
Following \cite[Section 3]{GJ}
we recall the Hochschild chain complex of a unital curved $A_\infty$-algebra
that generalizes the the Hochschild chain complex of a unital dg-algebra
(see {\it loc}. {\it cit}. for details).
Let $(A,b)$ be a unital curved $\AI$-algebra and $k\to A$ a unital
$A_\infty$-morphism which sends $1\in k$ to $1_A$. Put $\bar{A}=\Coker(k\to A)$
and let $\bar{b}:B\bA\to B\bA$ be the curved $\AI$-structure on $\bA$
associated to
$(A,b)$.
Consider the graded $k$-module 
\[
C_\bullet(A):=A\otimes B\bar{A}.
\]
We will define a differential on $C_\bullet(A)$.
Note first that graded module $B\bar{A}\otimes A \otimes B\bar{A}$ has a natural $B\bA$-bi-comodule structure given by $1_{B\bar{A}\otimes A} \otimes \Delta_{B\bar{A}}$ and $\Delta_{B\bar{A}}\otimes 1_{A \otimes B\bar{A}}$.
We will define a differential
\[
d:B\bar{A}\otimes A \otimes B\bar{A}\to B\bar{A}\otimes A \otimes B\bar{A}
\]
which commutes with the coderivation of $B\bar{A}$, that is, the equality
$(\bar{b}\otimes 1+1\otimes d)\circ (\Delta_{B\bar{A}}\otimes 1_{A \otimes B\bar{A}})=(\Delta_{B\bar{A}}\otimes 1_{A \otimes B\bar{A}})  \circ d$ holds,
and a similar condition holds for the right coaction (we think of it as
a $\bA$-bimodule structure on $A$).
As in the case of coderivation on $BA$, this data is uniquely determined by
the family of $k$-linear maps
\[
\{d_{m,n}:(\bA[1])^{\otimes m}\otimes A\otimes (\bA[1])^{\otimes n} \to A\}_{m,n\ge 0}.
\]
We define $d$ by giving the family $\{d_{m,n}\}_{m,n\ge0}$ such that the diagram
\[
\xymatrix{
(\bA[1])^{\otimes m}\otimes A\otimes (\bA[1])^{\otimes n} \ar[r]^(0.8){d_{m,n}} \ar[d]_{1^{\otimes m}\otimes s \otimes 1^{\otimes n}} & A \ar[d]^s \\
(\bA[1])^{\otimes m}\otimes A[1] \otimes (\bA[1])^{\otimes n} \ar[r]_(0.8){\bar{b}_{m+1+n}} & A[1].
}
\]
commutes.
The map $d_{m,n}$
sends $[a_1|\ldots|a_m]\otimes a \otimes [a_1'|\ldots|a_n']$
to
\[
(-1)^{\Sigma_{l=1}^m (|a_l|-1)}s^{-1}\bar{b}_{m+1+n}[a_1|\ldots|a_m|a|a_1'|\ldots|a_n'].
\]
Let $\Delta^R$ and $\Delta^L$ denote the canonical right and left coactions
of $B\bA$ on $B\bar{A}\otimes A \otimes B\bar{A}$.
Let $S_{2341}$ be the symmetric permutation
\[
B\bA \otimes (B\bar{A}\otimes A \otimes B\bar{A})\to (B\bA \otimes A \otimes  B\bA) \otimes B\bar{A},
\]
and put $\Phi(A):=\Ker(\Delta^R-S_{2341}\circ \Delta^L)\subset B\bar{A}\otimes A \otimes B\bar{A}$.
Since $\Delta^R$ and $S_{2341}\circ \Delta^L$
commute with differentials on the domain and the target,
$\Phi(A)$ inherits a differential from $B\bA \otimes A \otimes  B\bA$.
Moreover, $A\otimes B\bA\stackrel{1\otimes \Delta_{B\bA}}{\longrightarrow} A\otimes B\bA\otimes B\bA\stackrel{S_{312}}{\longrightarrow}B\bA\otimes A\otimes B\bA$
is injective, and its image is $\Phi(A)$.
Thus $C_\bullet(A)\simeq \Phi(A)$ inherits a differential $\partial_{\textup{Hoch}}$, and
we define a Hochschild chain complex of $A$ to be
a differential graded (dg) $k$-module $C_\bullet(A)$ endowed with $\partial_{\textup{Hoch}}$.
The Hochschild chain complex has a differential $B$ of degree $-1$ called Connes' operator given by
\[
B(a\otimes [a_1|\ldots |a_n])=\sum_{i=0}^n(-1)^{(\epsilon_{i-1}+|a|-1)(\epsilon_n-\epsilon_{i-1})}1_A\otimes [a_i|\ldots|a_n|a|a_{1}|\ldots|a_{i-1}],
\]
that satisfies $B\Ph+\Ph B=0$ and $B^2=0$,
where $\epsilon_i=\Sigma_{r=1}^{i}(|a_r|-1)$.
The last two identities mean that $(C_\bullet(A),\Ph,B)$ is a mixed complex
in the sense of \cite{Kas}. Put another way, if we let $\Lambda$ be
the dg algebra $k[\epsilon]/(\epsilon^2)$ with zero differential
such that the (cohomological) degree of the generator $\epsilon$ is $-1$,
then $(C_\bullet(A),\Ph)$ is a dg-$\Lambda$-module where
the action of $\epsilon$ is induced by $B$.

\subsection{}
\label{cyclic}
To the mixed complex $(C_\bullet(A),\Ph,B)$ we associate
its negative cyclic chain complex and its periodic cyclic chain complex.
Let 
\[
CC^-_\bullet(A):=(C_\bullet(A)[[t]], \Ph+tB),\ \ 
CC^{per}_\bullet(A):=(C_\bullet(A)((t)), \Ph+tB),
\]
where $t$ is a formal variable of (cohomological) degree two.
We consider the graded module $C_\bullet(A)[[t]]$ to be $\prod_{i\ge 0}C_\bullet(A)\cdot t^i$.
If $C_l(A)$ denotes the part of (homological) degree $l$
of $C_\bullet(A)$,
then the part of (cohomoligcal) degree $r$ of $C_\bullet(A)[[t]]$
is $\prod_{i\ge 0, 2i-l=r}C_l(A)\cdot t^i$.
The graded module $C_\bullet(A)((t))$ is regarded
as $\bigcup_{l\in \ZZ}\prod_{i\ge l}C_\bullet(A)\cdot t^i$.
The identities $\Ph^2=B\Ph+\Ph B=B^2=0$ implies $(\Ph+tB)^2=0$.
We call $CC^-_\bullet(A)$ (resp. $CC^{per}_\bullet(A)$)
the negative cyclic complex of $A$
(resp. the periodic cyclic complex of $A$).
The cohomology $HN_n(A):=H^{(-n)}(CC^-_\bullet(A))$
and $HP_n(A):=H^{(-n)}(CC^{per}_\bullet(A))$ is called
the negative cyclic homology and periodic cyclic homology respectively.
The peridic cyclic homology is periodic in the sense
that $HP_n(A)=HP_{n+2i}(A)$ for any integer $i$.
The module $HN_*(A)$
can be identified with the Ext-group $\Ext_{\Lambda}(k,C_\bullet(A))$.
To be precise, we let $k$ be the mixed complex of $k$ placed in degree zero
equipped with the trivial action of $\epsilon\in \Lambda$.
The category of dg $\Lambda$-module admits the projective model structure
where weak equivalences are quasi-isomorphisms, and fibrations are
degreewise surjective maps. We choose a cofibrant resolution $K$ of $k$:
\[
\cdots \stackrel{\cdot \epsilon}{\to} k\epsilon \stackrel{0}{\to} k \stackrel{\cdot \epsilon}{\to} k\epsilon \stackrel{0}{\to} k.
\]
Then the graded module of Hom complex $\Hom_{\Lambda}(K,C_\bullet (A))$
can naturally be identified with $C_\bullet(A)[[t]]$ (if
we denote by $1_{2i}$ the unit $1\in k$ placed in (cohomological) degree
$-2i$, $f:K\to C_\bullet(A)$ corresponds
to $\Sigma_{i=0}^{\infty}f(1_{2i})t^{i}$ in $C_\bullet(A)[[t]]$).
Unwinding the definition we see that 
the differential of the Hom complex corresponds to
the differential $\Ph+tB$. Similarly, the Hom complex $\Hom_{\Lambda}(K,k)$
is identified with $k[t]$.
Here we regard $k[t]$ as the dg algebra with zero differential
such that the (cohomological) degree of $t$ is $2$.
This dg algebra is the Koszul dual of $\Lambda$.
Let $\iota:k\to K$ be a canonical section of the resolution $K\to k$.
Then the natural $k[t]$-module structure on $CC^-_\bullet(A)[[t]]$
corresponds to
\[
\Hom_{\Lambda}(K,k)\otimes \Hom_{\Lambda}(K,C_\bullet(A))\to \Hom_{\Lambda}(K,C_\bullet(A)),\ \ \phi\otimes f\mapsto f \circ \iota\circ \phi.
\]

\section{Lie algebra action and modular interpretation}
\label{LM}

\subsection{}
\label{Lieaction}
Let $k$ be a base field of characteristic zero.
Let $A$ be a unital $\AI$-algebra or a unital dg algebra.
In what follows we write $C^\bullet(A)[1]$ 
for the dg Lie algebra $\overline{C}^\bullet(A)[1]$ of the normalized
Hochschild cochain complex,
and $C^\bullet(A):=(C^\bullet(A)[1])[-1]$.
The dg Lie algebra $C^\bullet(A)[1]$ acts on the Hochschild chain complex
$C_\bullet(A)$. Our reference for this action is
Tamarkin-Tsygan \cite{TT}.
It is a part of data of a homotopy calculus algebra
on $(C^\bullet(A),C_\bullet(A))$ (see Section~\ref{calcdef} for algebra over calculi). We describe the action of the dg Lie algebra
$C^\bullet(A)[1]$ on $C_\bullet(A)$ in detail.

Let $P\in C^\bullet(A)[1]$ be a homogeneous element,
and suppose $P$ lies in $\Hom((A[1])^{\otimes l},A[1])\subset \oplus_{n\ge0}\Hom((A[1])^{\otimes n},A[1])$.
Let $|P|$ be the degree of $P$ in the cochain complex $C^\bullet(A)$
(thus the degree of $P$ in $C^\bullet(A)[1]$ is $|P|-1$).
Let $\mu_i=\Sigma_{r=0}^{i}(|a_r|-1)=|a_0|-1+\epsilon_i$.
A linear map $L_P:C_\bullet(A)\to C_\bullet(A)$ of degree $|P|-1$,
which we regard as an element of the endomorphism complex $\End_k(C_\bullet(A))$
of degree $|P|-1$, is defined by the formula
\begin{eqnarray*}
L_P(a_0\otimes[a_1|\ldots|a_n])
=\sum_{0\le j,\ j+l\le n}(-1)^{(|P|-1)\mu_j}a_0\otimes [a_1|\ldots|a_j|P[a_{j+1}|\ldots|a_{j+l}]|\ldots |a_n] \\
+\sum_{0\le i\le n \atop P\ \textup{includes}\ a_0}(-1)^{\mu_i(\mu_n-\mu_i)}
s^{-1} P([a_{i+1}|\ldots|a_n|a_0|\ldots|a_j ])\otimes [ a_{j+1}|\ldots|a_i]
\end{eqnarray*}

It gives rise to a graded $k$-linear map
\[
L:C^\bullet (A)[1]\to \End_k(C_\bullet(A))
\]
which carries $P$ to $L_P$.

\begin{Proposition}{$($cf. \cite[3.3.2]{TT}.$)$}
\label{Laction}
Let $b:BA\to BA$ denote the $\AI$-structure of $A$.
Let $P,Q$ be elements in $C^\bullet(A)[1]$.
The followings hold:

\begin{enumerate}
\renewcommand{\labelenumi}{(\theenumi)}

\item $L_{[P,Q]}=[L_P,L_Q]_{\textup{End}}:=L_P\circ L_Q-(-1)^{|L_P||L_Q|}L_Q\circ L_P$,

\item $\partial^{\textup{End}}L_P-L_{\PH P}=0$ where $\partial^{\textup{End}}L_P=\Ph\circ L_P-(-1)^{|L_P|}L_P\circ \Ph$,

\item $[B,L_P]_{\textup{End}}=0$.

\end{enumerate}
\end{Proposition}

\begin{Remark}
By (1) and (2) of this Proposition $L:C^\bullet (A)[1]\to \End_k(C_\bullet(A))$ is a map of dg Lie algebras, where $\End_k(C_\bullet(A))$ is endowed with
the bracket given by $[f,g]=f\circ g-(-1)^{|f||g|}g\circ f$.
The condition (3) means that $L$ factors through the dg Lie algebra
$\End_{\Lambda}(C_{\bullet}(A))$ of the endomorphism dg Lie algebra
of the $\Lambda$-module $C_{\bullet}(A)$.
\end{Remark}

\proof
(1) and (3) are nothing but \cite[3.3.2]{TT},
(but unfortunately its proof is omitted).
Thus, for the reader's convenience we here give the proof
since we will use these formulas.

We will prove (1). We first calculate $L_P\circ L_Q$.
We may and will assume that
$P$ and $Q$ belong to $\Hom((A[1])^{\otimes u},A[1])$ and
$\Hom((A[1])^{\otimes v},A[1])$ respectively.
For ease of notation we write $(a_0,a_1,\ldots,a_n)$ for $a_0\otimes [a_1|\ldots| a_n]$
in $A\otimes (\bar{A}[1])^{\otimes n}$.
\begin{eqnarray*}
L_P\circ L_Q=\sum_{i,j} A(P,Q,i,j)+\sum_{i,j} B(P,Q,i,j)+\sum_{i,j} C(P,Q,i,j)+\sum_{i,j} D(P,Q,i,j) \\
+\sum_{i,j} E(P,Q,i,j)+ \sum_{i,j} F(P,Q,i,j) +\sum_{i,j} G(P,Q,i,j) +\sum_{i,j} H(P,Q,i,j),
\end{eqnarray*}
where
\begin{align*}
A(P,Q,i,j)&=(-1)^{(|Q|-1)\mu_i+(|P|-1)\mu_j}(a_0,\ldots,P[a_{j+1}|\ldots],\ldots ,Q[a_{i+1}|\ldots],\ldots), \\
B(P,Q,i,j)&=(-1)^{(|Q|-1)\mu_i+(|P|-1)\mu_j}(a_0,\ldots,P[a_{j+1}|\ldots|Q[a_{i+1}|\ldots]|\ldots],\ldots), \\
C(P,Q,i,j)&=(-1)^{(|Q|-1)\mu_i+(|P|-1)(\mu_j+|Q|-1)}(a_0,\ldots,Q[a_{i+1}|\ldots],\ldots ,P[a_{j+1}|\ldots],\ldots), \\
D(P,Q,i,j)&=(-1)^{(|Q|-1)\mu_i+(\mu_j+|Q|-1)(\mu_n-\mu_j)}(s^{-1}P[a_{j+1}|\ldots|a_0|\ldots|Q[a_{i+1}|\ldots]|\ldots],\ldots,a_j), \\
E(P,Q,i,j)&=(-1)^{(|Q|-1)\mu_i+(\mu_j+|Q|-1)(\mu_n-\mu_j)}(s^{-1}P[a_{j+1}|\ldots|a_0|\ldots],\ldots,Q[a_{i+1}|\ldots],\ldots a_j), \\
F(P,Q,i,j)&=(-1)^{(|Q|-1)\mu_i+\mu_j(\mu_n-\mu_j+|Q|-1)}(s^{-1}P[a_{j+1}|\ldots|Q[a_{i+1}|\ldots]|\ldots|a_0|\ldots],\ldots,a_j), \\
G(P,Q,i,j)&=(-1)^{\mu_i(\mu_n-\mu_i)+(|P|-1)(|Q|-1+\mu_n-\mu_i+\mu_j)}(s^{-1}Q[a_{i+1}|\ldots|a_0|\ldots],\ldots,P[a_{j+1}|\ldots],\ldots,a_i), \\
H(P,Q,i,j)&=(-1)^{\mu_i(\mu_n-\mu_i)+(\mu_i-\mu_j)(|Q|-1+\mu_n-\mu_i+\mu_j)}(s^{-1}P[a_{j+1}|\ldots |Q[a_{i+1}|\ldots|a_0|\ldots]|\ldots],\ldots,a_j).
\end{align*}
Here $i$ and $j$ run over an adequate range
which depend on the type of terms in each summation.
Interchanging $(P,i)$ and $(Q,j)$ we put
\[
A(Q,P,j,i)=
(-1)^{(|Q|-1)\mu_i+(|P|-1)\mu_j}(a_0,\ldots,Q[a_{i+1}|\ldots],\ldots ,P[a_{j+1}|\ldots],\ldots).
\]
Other terms $B(Q,P,j,i), C(Q,P,j,i),\ldots$ are defined in a similar way.
For example,
\[
C(Q,P,j,i)=(-1)^{(|P|-1)\mu_j+(|Q|-1)(\mu_i+|P|-1)}(a_0,\ldots,P[a_{j+1}|\ldots],\ldots ,Q[a_{i+1}|\ldots],\ldots).
\]
Then we have
\begin{eqnarray*}
L_Q\circ L_P=\sum_{i,j} A(Q,P,j,i)+\sum_{i,j} B(Q,P,j,i)+\sum_{i,j} C(Q,P,j,i)+\sum_{i,j} D(Q,P,j,i) \\
+\sum_{i,j} E(Q,P,j,i)+ \sum_{i,j} F(Q,P,j,i) +\sum_{i,j} G(Q,P,j,i) +\sum_{i,j} H(Q,P,j,i)
\end{eqnarray*}
Notice that
\begin{eqnarray*}
A(P,Q,i,j)-(-1)^{(|P|-1)(|Q|-1)}C(Q,P,j,i)&=&0, \\
C(P,Q,i,j)-(-1)^{(|P|-1)(|Q|-1)}A(Q,P,j,i)&=&0, \\
E(P,Q,i,j)-(-1)^{(|P|-1)(|Q|-1)}G(Q,P,j,i)&=&0, \\
G(P,Q,i,j)-(-1)^{(|P|-1)(|Q|-1)}E(Q,P,j,i)&=&0.
\end{eqnarray*}
Therefore, the terms $A,C,E$ and $G$ do not appear in $L_{P}\circ L_Q-(-1)^{|L_P||L_Q|} L_P\circ L_Q$.

Next we will calculate $L_{[P,Q]}$.
Note first that $L_{[P,Q]}=L_{P\circ Q}-(-1)^{(|P|-1)(|Q|-1)}L_{Q\circ P}$.
The composition $P\circ Q\in \Hom((A[1])^{\otimes u+v-1},A[1])$ is given by
\[
P\circ Q[a_1|\ldots |a_n]=\sum_{i}(-1)^{(|Q|-1)(\Sigma_{r=1}^{i}(|a_r|-1)}P[a_1|\ldots|Q[a_{i+1}|\ldots|a_{i+v}]|\ldots|a_{u+v-1}].
\]
Using this we see that $L_{P\circ Q}(a_0,a_1,\ldots,a_n)$ is equal to
\begin{eqnarray*}
& &\sum_{i,j}(-1)^{(|P|-1+|Q|-1)\mu_j+(|Q|-1)(\mu_i-\mu_j)}(a_0,\ldots ,P[a_{j+1}|\ldots|Q[a_{i+1}|\ldots]|\ldots],\ldots) \\
&+&\sum_{i,j}(-1)^{\mu_j(\mu_n-\mu_j)+(|Q|-1)(\mu_n-\mu_j+\mu_i)}(s^{-1}P[a_{j+1}|\ldots|a_0|\ldots|Q[a_{i+1}|\ldots]|\ldots],\ldots,a_j) \\
&+&\sum_{i,j}(-1)^{\mu_j(\mu_n-\mu_j)+(|Q|-1)(\mu_i-\mu_j)}(s^{-1}P[a_{j+1}|\ldots|Q[a_{i+1}|\ldots]|\ldots|a_0|\ldots],\ldots,a_j) \\
&+&\sum_{i,j} (-1)^{\mu_j(\mu_n-\mu_j)+(|Q|-1)(\mu_i-\mu_j)}(s^{-1}P[a_{j+1}|\ldots |Q[a_{i+1}|\ldots|a_0|\ldots]|\ldots],\ldots,a_j)
\end{eqnarray*}
where $i$ and $j$ run over adequate range in each summation.
Thus by an easy computation of signs, we see that it is equal to
\[
\sum_{i,j}B(P,Q,i,j)+\sum_{i,j}D(P,Q,i,j)+\sum_{i,j}F(P,Q,i,j)+\sum_{i,j}H(P,Q,i,j).
\]
Similarly, we see that $L_{Q\circ P}$ is equal to
\[
\sum_{i,j}B(Q,P,j,i)+\sum_{i,j}D(Q,P,j,i)+\sum_{i,j}F(Q,P,j,i)+\sum_{i,j}H(Q,P,j,i).
\]
Therefore we deduce that $L_{[P,Q]}=[L_P,L_Q]$.

Next we prove (2).
By a calculation of the differential $\Ph$ below
(Lemma~\ref{explicitHoch}),
$L_b=\Ph$.
Combined with $\PH P=[b,P]$, we see that (2) follows from (1).

Finally, we prove  (3). We compare $B\circ L_P$ with $L_P\circ B$.
\begin{align*}
&B\circ L_P(a_0,a_1,\ldots,a_n) \\
&  =\sum_{i,j}(-1)^{(|P|-1)\mu_i+(\mu_j+|P|-1)(\mu_n-\mu_j)}(1_A,a_{j+1},\ldots,a_n,a_0,\ldots,P[a_{i+1}|\ldots],\ldots,a_j) \\
&  \ \ \ \ +\sum_{i,j}(-1)^{(|P|-1)\mu_i+\mu_j(\mu_n-\mu_j+|P|-1)}(1_A,a_{j+1},\ldots,P[a_{i+1}|\ldots],\ldots,a_0,\ldots,a_j) \\
&   \ \ \ \ +\sum_{i,j}(-1)^{\mu_i(\mu_n-\mu_i)+(\mu_i-\mu_j)(\mu_n-\mu_i+\mu_j+|P|-1)}(1_A,a_{j+1},\ldots,P[a_{i+1}|\ldots|a_0|\ldots],\dots,a_j).
\end{align*}
On the other hand, Since $P$ is normalized we have
\begin{align*}
&L_P\circ B(a_0,a_1,\ldots,a_n) \\
&  =\sum_{i,j}(-1)^{\mu_j(\mu_n-\mu_j)+(|P|-1)(1+\mu_n-\mu_j+\mu_i)}(1_A,a_{j+1},\ldots,a_n,a_0,\ldots,P[a_{i+1}|\ldots],\ldots,a_j) \\
&   \ \ \ \ +\sum_{i,j}(-1)^{\mu_j(\mu_n-\mu_j)+(|P|-1)(1+\mu_i-\mu_j)}(1_A,a_{j+1},\ldots,P[a_{i+1}|\ldots],\ldots,a_j).
\end{align*}
where $P[\ldots]$ is allowed to include $a_0$ in the second summation.
Comparing $B\circ L_P$ with $L_P\circ B$ we see that $B\circ L_P-(-1)^{(|P|-1)}L_P\circ B=0$.
\QED

\subsection{}
\label{cylicify}
We continue to assume that $(A,b)$ be a unital $\AI$-algebra (or a unital
dg algebra).
Consider the pair $(C^\bullet(A)[1], C_\bullet(A))$
of the shifted (normalized) Hochschild cochain complex and the Hochschild chain complex
of $A$. According to Proposition~\ref{Laction}, the dg Lie algebra
$C^\bullet(A)[1]$ acts on the mixed complex $C_\bullet(A)$.
That is, there is an (explicit) morphism $L:C^\bullet(A)[1]\to \End_\Lambda(C_\bullet(A))$ of dg Lie algebras.
Let us consider the morphism of dg Lie algebras
$\End_{\Lambda}(C_\bullet(A))\to \End_{k}(\Hom_{\Lambda}(K,C_\bullet(A)))$
which carries $f\in C_\bullet(A)\to C_\bullet(A)$ to $\{\phi \mapsto f\circ \phi\}\in \End_{k}(\Hom_{\Lambda}(K,C_\bullet(A)))$,
where $K$ is the cofibrant resolution of $k$  (see Section~\ref{cyclic}).
Remember that
there exists the natural isomorphism $\Hom_{\Lambda}(K,C_\bullet(A))\simeq (C_\bullet(A)[[t]],\Ph+tB)$. Moreover, $\{\phi \mapsto f\circ \phi\}$ commutes
with the action of $\Hom_{\Lambda}(K,k)\simeq k[t]$. Consequently,
we have the composition of maps of dg Lie algebras
\[
C^\bullet(A)[1]\stackrel{L}{\to} \End_{\Lambda}(C_\bullet(A))\to \End_{k[t]}(C_\bullet(A)[[t]]):=\End_{k[t]}((C_\bullet(A)[[t]],\Ph+tB)).
\]
We easily observe that the image of $P\in C^\bullet(A)[1]$ in $\End_{k[t]}(C_\bullet(A)[[t]])$
is $L_P[[t]]$. We shall denote by $L[[t]]$ this composite.
We write $k[t^\pm]$ for $k[t,t^{-1}]$.
By tensoring with $\otimes_{k[t]}k[t^\pm]$
we also have
\[
L((t)): C^\bullet(A)[1] \to \End_{k[t^\pm]}(C_\bullet(A)((t))):=\End_{k[t^\pm]}((C_\bullet(A)((t)),\Ph+tB)).
\]

\subsection{}
\label{premoduli}
To a dg Lie algebra one can associate a deformation functor.
There are several formalisms of deformation theories.
We employ the simplest form that dates back to
Schlessinger's formalism \cite{Schl}.
(There is an enhanced formalism, see e.g. Hinich \cite{Hin0}.)
We review the setting of deformation theories
we shall employ
in this paper.
Let $E$ be a nilpotent dg Lie algebra.
An element $x$ of degree one in
$E$ equipped with a differential $d$ and a bracket $[-,-]$ is said to be a Maurer-Cartan element if the Maurer-Cartan equation $dx+\frac{1}{2}[x,x]=0$ holds.
We denote by $\MC(E)$ the set of Maurer-Cartan elements.
The space $E^0$ of degree zero is a (usual) Lie algebra.
Let $\exp(E^0)$ be the exponential group associated to the nilpotent
Lie group $E^0$ whose product is given by the Baker-Campbell-Hausdorff
product on $E^0$.
The Lie algebra $E^0$ acts on $E^1$ by $E^0\to \End_k(E^1),\ \alpha\mapsto [\alpha,-]-d\alpha$.
It gives rise to an action of $\exp(E^0)$ on the space $E^1$ given by
\[
\mu \mapsto e^\alpha\bullet \mu :=e^{ad(\alpha)}\mu-\int_0^1(e^{ad(s\alpha)}d\alpha) ds.
\]
for $\mu$ in $E^1$.
This action preserves Maurer-Cartan elements.
Let $x,y\in MC(E)$. We say that $x$ is gauge equivalent to $y$
if both elements coincide in $\MC(E)/\exp(E^0)$.
The element $x$ is gauge equivalent to $y$ if there is
$\alpha\in E^0=\exp(E^0)$ such that $e^\alpha\bullet x=y$.
There is another more theoretical definition:
If $\Omega_1=k[u,du]$ denotes the dg algebra
of 1-dimensional polynomial differential forms (see Section~\ref{mappingsp}),
the coequalizer of two degeneracies $d_0,d_1:\MC(\Omega_1\otimes E)\rightrightarrows \MC(E)$ determined by $u=du=0$
and $u=1,\ du=0$ is isomorphic to
$\MC(E)/\exp(E^0)$. This definition can also be applicable to
$L_\infty$-algebras.
Let $\phi:E\to E'$ be a map of nilpotent dg Lie algebras (or more generally, an $L_\infty$-morphism, cf. Section~\ref{mappingsp}). Then $\phi$ induces $\MC(E)\to \MC(E')$ which commutes with gauge actions.
In particular, it gives rise to $\MC(E)/\exp(E^0) \to \MC(E')/\exp((E')^0)$.

Let $\Art$ be the category of artin local $k$-algebras with
residue field $k$.
For $R$ in $\Art$, we write $m_R$ for the maximal ideal of $R$.
Let $E$ be a dg Lie algebra.
The bracket on $E\otimes m_R$ is defined by $[e\otimes m, e'\otimes m']=[e,e']\otimes mm'$.
Note that the dg Lie algebra $E\otimes m_R$ is nilpotent.
For a morphism $R\to R'$ in $\Art$, the induced map $E\otimes_km_R\to E\otimes_km_{R'}$ determines a map
$\MC(E\otimes m_R)/\exp(E^0\otimes m_R)\to \MC(E\otimes m_{R'})/\exp(E^0\otimes m_{R'})$.
We then define a functor
\[
\Def_E:\Art\to \Set
\]
which carries $R$ to $\MC(E\otimes m_R)/\exp(E^0\otimes m_R)$,
where $\Set$
is the category of sets. 
By this notation
we think of $\Def_E$ as an analogue of formal schemes in scheme theory.

If $\phi:E\to E'$ is a quasi-isomorphism of dg Lie algebras,
the induced morphism $\Def_{\phi}:\Def_{E}\to \Def_{E'}$
is an equivalence (this is a theorem stated and proved
by Deligne, Goldman and  Millson, Hinich, Kontsevich, Fukaya, Getzler and others).

\vspace{2mm}

Let us consider two examples of deformation problems.

\begin{Definition}
For a complex $(E,d)$ (i.e. a dg $k$-module) and $R$ in $\Art$,
a deformation of the complex $E$ to $R$ is a dg $R$-module $(E\otimes_k R,\tilde{d})$ such that
the reduction of the differential $\tilde{d}$ of $E\otimes_kR$ to the complex $E\otimes_kR/m_R$ exhibits $E\otimes_kk$ as the initial complex $(E,d)$.
Let $(E\otimes_kR,\tilde{d}_1)$ and $(E\otimes_kR,\tilde{d}_2)$
be two deformations of $E$ to $R$.
An isomorphism of these deformations is
an isomorphism $E\otimes_kR\to E\otimes_kR$ of dg $R$-modules
which induces the identity $E\otimes_kR/m_R\to E\otimes_kR/m_R$.
Note that by nilpotent Nakayama lemma, every
homomorphism $E\otimes_kR\to E\otimes_kR$ of dg $R$-modules
which induces the identity $E\otimes_kR/m_R\to E\otimes_kR/m_R$,
is an isomorphism.

\end{Definition}

Let us consider the functor $\mathsf{D}_E:\Art\to \Set$ which carries
$R$ to the set of isomorphism classes of deformations of $E$ to
$R$ (the functoriality is defined in the natural way, that is,
a morphism $R\to R'$ in $\Art$ induces the
base change $E\otimes_kR\otimes_RR'$).

\begin{Proposition}
\label{compdeform}
There is a natural equivalence of functors $\Def_{\End_{k}(E)}\to \mathsf{D}_E$.
For each $R$ in $\Art$ it carries a Maurer-Cartan element $x$ in $\End_k(E)\otimes m_R$
to a deformation $(E\otimes_kR,d\otimes_k R+x)$, where the differential is given by
$(d\otimes_kR+x)(e\otimes r)=d(e)\otimes r+x(e)r$.
\end{Proposition}

In Introduction, the following
deformations are referred to as curved $A_\infty$-deformations of $A$.

\begin{Definition}
\label{algdeformdef}
Let $(A,b)$ be an $\AI$-algebra with unit $1_A$
(or a dg algebra) and $R$ an artin local $k$-algebra in $\Art$.
A deformation of $(A,b)$ to $R$ is
a graded $R$-module $A\otimes_kR$ endowed with a curved $\AI$-structure
$\tilde{b}:BA\otimes_kR\to BA\otimes_kR$ (over $R$) such that
(i) the reduction $BA\otimes_kR/m_R\to BA\otimes R/m_R$ is $b: BA\to BA$ via the canonical
isomorphism $BA\otimes_kR/m_R\simeq BA$, and
(ii) the element $1_A\otimes 1_R\in A\otimes_kR$ is a unit.
Here $BA\otimes_kR$ is a graded coalgebra over $R$.
Let $(A\otimes_kR,\tilde{b}_1)$ and $(A\otimes_kR,\tilde{b}_2)$
be two deformations of $(A,b)$ to $R$.
An isomorphism $(A\otimes_kR,\tilde{b}_1)\to (A\otimes_kR,\tilde{b}_2)$
of deformations is an $A_\infty$-isomorphism $BA\otimes_kR\to BA\otimes_kR$
(of dg coalgebras) over $R$, such that
the reduction $BA\otimes_kR/m_R\to BA\otimes_kR/m_R$ is the identity.
(As in the case of complexes, a homomorphism of dg coalgebras over $R$ whose reduction is the identity, is an isomorphism.)
Note that a deformed $\AI$-structure is allowed to be curved.
\end{Definition}

Let $\mathsf{DAlg}_{A}$ denote the functor $\Art\to \Set$
which carries $R$ to the set of isomorphism classes of
deformations of $(A,b)$ to $R$.

\begin{Proposition}
\label{algdeform}
There is a natural equivalence $\Def_{C^\bullet(A)[1]}\to \mathsf{DAlg}_{A}$.
It carries a Maurer-Cartan element $x\in \MC(C^\bullet(A)[1]\otimes_k m_R)$
to a curved $\AI$-structure $b\otimes_k R+x:BA\otimes_kR\to BA\otimes_kR$.
$($Here $C^\bullet(A)[1]$ is the normalized Hochschild cochain complex,
$x:BA\to BA\otimes_km_R$,
and $(b\otimes_k R+x)(a\otimes r)=b(a)\otimes r+x(a)r$.$)$
\end{Proposition}

\begin{Remark}
\label{moritaremark}
The dg Lie algebra $C^\bullet(A)[1]$ is derived Morita invariant (cf. \cite{Kel}).
Let $\textup{Perf}_A$ be the dg category (or some model
of a stable $(\infty,1)$-category)
of
perfect dg (left) $A$-modules.
A dg $A$-module $M$ is said to be perfect
if at the level of the homotopy category,
$M$ belongs to the smallest triangulated
subcategory of the triangulated category of dg $A$-modules, which contains
$A$ and is closed under direct summands.
If there is a quasi-equivalence $\textup{Perf}_A\simeq \textup{Perf}_B$,
there is an equivalence $C^\bullet(A)[1]\simeq C^\bullet(B)[1]$
(more precisely, an $L_\infty$-quasi-isomorphism).
In particular, $\Def_{C^\bullet(A)[1]}\simeq \Def_{C^\bullet(B)[1]}$.
Thus, it seems that it is natural to describe
$\Def_{C^\bullet(A)[1]}$ in terms of category theory.
Namely, it is worthwhile to attempt to describe $\Def_{C^\bullet(A)[1]}$
as the functor of deformations
of the dg category (or stable $(\infty,1)$-category)
$\textup{Perf}_A$ (see Keller-Lowen's work \cite{KL}
for the progress of this problem).
At the time of writing this paper,
the author does not know a category-theoretic 
formulation that fits into nicely with nilpotent deformations to {\it curved $A_\infty$-algebras}.
Therefore, we choose the approach using curved deformations of the algebra
$A$.
\end{Remark}

\begin{Remark}
\label{comparisonscheme}
To get the feeling of deformations of a $A_\infty$-algebra $A$
(or a dg algbera), let us consider the situation when
$A$ comes from quasi-compact separated scheme $X$ over $k$.
According to \cite{BV} the derived category $D(X)$ of (unbounded) quasi-coherent complexes on
$X$ admits a compact generator, that is, a compact object $C$
such that for any $D\in D(X)$
the condition $\Ext_{D(X)}^*(C,D)=0$
implies that $D$ is quasi-isomorphic to $0$.
It follows that $D(X)$ is quasi-equivalent to the dg category 
of dg (left) modules over some dg algebra $A$ as dg-categories,
namely, $A$ comes from $X$. 
Indeed, one can choose $A$ to
be the endomorphism dg algebra of $C$ (after a suitable resolution).
(Here we regard $D(X)$ as a dg category.)
For simplicity, suppose that $X$ is smooth and proper over $k$.
On one hand, Hochschild-Kostant-Rosenberg theorem implies
that $HH^n(A)=\oplus_{i+j=n}H^j(X,T_X^{\wedge i})$.
On the other hand, by the definition and Propsoition~\ref{algdeform}
\[
\mathsf{DAlg}_A(k[\epsilon]/(\epsilon^2))\simeq\Def_{C^\bullet(A)[1]}(k[\epsilon]/(\epsilon^2))\simeq HH^2(A),
\]
where $\Def_{C^\bullet(A)[1]}(k[\epsilon]/(\epsilon^2))$
is the tangent space of $\Def_{C^\bullet(A)[1]}$.
Consequently, we have
\[
\mathsf{DAlg}_A(k[\epsilon]/(\epsilon^2))\simeq\Def_{C^\bullet(A)[1]}(k[\epsilon]/(\epsilon^2))\simeq H^0(X,\wedge^2T_X)\oplus H^1(X,T_X)\oplus H^2(X,\OO_X).
\]
While there are several interpretations,
we review one modular interpretation of the space of
the right hand side, though details remain to be elucidated.
The subspace $H^1(X,T_X)$ parametrizes (first-order)
deformations of the scheme $X$ that induce deformations
of $D(X)$ (see Remark~\ref{moritaremark}).
Morally, the subspace $H^0(X,\wedge^2T_X)\oplus  H^2(X,\OO_X)$
should be considered to be the part of noncommutative deformations.
For example, a Poisson structure on $X$ belongs
to $H^0(X,\wedge^2T_X)$,
and it gives rise to a deformation quantization of $X$.
The space $H^2(X,\OO_X)$ is identified with the space of liftings of the zero $H^2_{\textup{\'et}}(X,\mathbb{G}_m)$ to $H^2_{\textup{\'et}}(X\times_kk[\epsilon]/(\epsilon^2),\mathbb{G}_m)$.
By a main theorem of \cite{ToA}, such a lifting can be thought of as
a (first-order) deformation of the structure sheaf $\OO_X$ as a derived Azumaya algebra. In other words, this sort of deformations may be described
as deformations of twisted sheaves (see \cite{Yu}).
\end{Remark}

\subsection{}
Proposition~\ref{compdeform} and~\ref{algdeform} are well-known to experts.
In particular, deformations of (curved) deformations of algebraic structure
is one of main subjects of Hochschild cohomology.
We present the proof
for reader's
convenience because we are unable to find the literature
that fits in with our curved $A_\infty$-setting.
Proposition~\ref{algdeform} is a consequence of the following
two Claims.

\begin{Claim}
\label{deformclaim1}
Let $\operatorname{Defor}_{\textup{Alg}}(A,R)$ be the set of deformations $\tilde{A}=(A\otimes_kR,\tilde{b})$ of $(A,b)$ to $R$.
There is a natural bijective map
\[
\MC(C^\bullet(A)[1]\otimes_km_R)\stackrel{\sim}{\longrightarrow} \operatorname{Defor}_{\textup{Alg}}(A,R)
\]
which carries a coderivation $f:BA\to BA\otimes_k m_R$ of degree one
satisfying the Maurer-Cartan equation $\PH f+\frac{1}{2}[f,f]=0$
to a coderivation $\PH\otimes_kR+f\otimes R:BA\otimes_kR\to BA\otimes_kR$.
\end{Claim}

\proof
Let $f:BA\to BA\otimes_k m_R$ be a coderivation of degree one,
that is, an element of $(C^\bullet(A)[1]\otimes m_R)^1$.
We abuse notation by writing $f$ for $f\otimes R:BA\otimes_kR\to BA\otimes_kR$
given by $f(x\otimes r)=f(x)r$.
Then $(b\otimes_kR+f)^2=(b\otimes_kR)^2+b\circ f+f\circ b+f\circ f=\PH f+\frac{1}{2}[f,f]$.
Thus $b\otimes_kR+f$ is a square-zero coderivation if and only if
$f$ is a Maurer-Cartan element. The map $f\mapsto b\otimes_kR+f$
is injective. In addition, note that $1_A\otimes_k1_R$ is a unit
in $(A\otimes_kR, b\otimes_kR+f)$ exactly when
$f$ is normalized. It remains to be proved that the map is surjective.
To this end, it suffices to observe only that
if a square-zero coderivation $\tilde{b}:BA\otimes_kR\to BA\otimes_kR$ over $R$
whose reduction
%\[
%BA\simeq BA\otimes_kk\subset BA\otimes_kR\stackrel{\tilde{b}}{\to}BA\otimes_kR\to BA\otimes_kR/m_R\simeq BA
%\]
is $b$, then $(\tilde{b} - b\otimes_kR)\circ \iota$ belongs to $C^\bullet(A)[1]\otimes_km_R$ where $\iota:BA\otimes_kk\to BA\otimes_kR$
is the canonical inclusion.
\QED

\begin{Claim}
\label{deformclaim2}
Let $\operatorname{Aut}(A,R)$ be the group of automorphisms $BA\otimes_kR\to BA\otimes_kR$ of the graded coalgebra $BA\otimes_kR$ (not equipped with any coderivation)
whose reduction $BA\stackrel{\sim}{\to}BA$ is the identity.
There is a natural isomorphism
\[
\exp(C^\bullet(A)[1]\otimes_km_R)^0\stackrel{\sim}{\longrightarrow} \Aut(A,R)
\]
which carries a coderivation $d:BA\to BA\otimes_k m_R$ of degree zero
to an automorphism $e^d=\Sigma_{i=0}^\infty\frac{1}{i!}d^i:BA\otimes_kR\to BA\otimes_kR$.

Consider the action of $\Aut(A,R)$ on $\operatorname{Defor}_{\textup{Alg}}(A,R)$ given by $\tilde{b}\mapsto c\circ \tilde{b}\circ c^{-1}$ for $c\in \Aut(A,R)$.
Then the gauge action of $\exp(C^\bullet(A)[1]\otimes_km_R)^0$
on $\MC(C^\bullet(A)[1]\otimes_km_R)$ commutes with the action of $\Aut(A,R)$
on $\operatorname{Defor}_{\textup{Alg}}(A,R)$.
\end{Claim}

\Proof
Since $(\Delta_{BA}\otimes_km_R)\circ d=(d\otimes 1+1\otimes d)\circ \Delta_{BA}:BA\to BA\otimes_kBA\otimes_km_R$, $e^d=\Sigma_{i=0}^\infty\frac{1}{i!}d^i$
and $e^{d\otimes1+1\otimes d}=e^d\otimes e^d$ satisfy the commutativity
\[
\xymatrix{
BA\otimes_k R \ar[r]^(0.3){\Delta_{BA\otimes R}} \ar[d]_{e^d} & (BA\otimes_kR)\otimes_R (BA\otimes_kR) \ar[d]^{e^{d\otimes 1+1\otimes d}} \\
BA\otimes_k R\ar[r]_(0.3){\Delta_{BA\otimes R}} & (BA\otimes_kR)\otimes_R (BA\otimes_kR)
}
\]
which means that $e^d$ is an automorphism of the graded coalgebra $BA\otimes_kR$. By the construction, the reduction $BA\otimes_kR/m_R\to BA\otimes_kR/m_R$
is the identity.

Next we prove our claim by induction on the length of $R$.
The case $R=k$ is obvious.
Let
\[
0 \to I \to R\to R'\to 0
\]
be an exact sequence where $R\to R'$ is a sujective map of artin local $k$-algebras
and $I$ is the kernel such that $I\cdot m_R=0$.
We assume that $(C^\bullet(A)[1]\otimes_km_{R'})^0\to  \Aut(A,R')$,
which carries $d$ to $e^d$, is an isomorphism.
We will prove that $(C^\bullet(A)[1]\otimes_km_{R})^0\to  \Aut(A,R)$
is an isomorphism.
Let $d\in (C^\bullet(A)[1]\otimes_km_{R})^0$ and $h\in (C^\bullet(A)[1]\otimes_kI)^0$.
Then using $I\cdot m_R=0$, we see that
$e^{d+h}=1+(d+h)+\frac{1}{2}(d+h)^2+\cdots=e^d+h$.
Combined this equality with the assumption on induction we deduce that
$(C^\bullet(A)[1]\otimes_km_{R})^0\to  \Aut(A,R)$ is injective.
To prove that the map is surjective, let $f\in \Aut(A,R)$
and take $d\in C^\bullet(A)[1]\otimes_km_R$ such that
$f=d\ \textup{mod}\ I$ in $(C^\bullet(A)[1]\otimes_km_{R'})^0\simeq \Aut(A,R')$.
If $h=f-e^d:C\to C\otimes_kI$ is coderivation,
then $e^{d+h}=e^d+h=f$ implies that the map is surjective.
Note that $(f\otimes f-e^d\otimes e^d) \circ \Delta_{BA}=(\Delta_{BA}\otimes_kI)\circ (f-e^d):BA\to BA\otimes BA\otimes_k I$.
Put $\bar{f}=f-1$ and $\bar{e}^d=e^d-1$.
Then 
\begin{eqnarray*}
f\otimes f-e^d\otimes e^d &=& (1+\bar{f})\circ (1+\bar{f})- (1+\bar{e}^d)\circ (1+\bar{e}^d) \\
&=& (f-e^d)\otimes 1+1\otimes(f-e^d) +\bar{f}\otimes\bar{f}-\bar{e}^d\otimes \bar{e}^d
\end{eqnarray*}
The term $\bar{f}\otimes \bar{f}-\bar{e}^d\otimes \bar{e}^d=(\bar{f}-\bar{e}^d)\otimes \bar{f}+\bar{e}^d\otimes (\bar{f}-\bar{e}^d)$
is zero since $\bar{f}-\bar{e}^d\in C^\bullet(A)[1]\otimes_kI$,
$\bar{f}, \bar{e}^d\in C^\bullet(A)[1]\otimes_k m_R$ and $I\cdot m_R=0$.
Therefore, $f-e^d$ is a coderivation.

Finally, we prove that
the gauge action of $\exp(C^\bullet(A)[1]\otimes_km_R)^0$
on $\MC(C^\bullet(A)[1]\otimes_km_R)$ commutes with the action of $\Aut(A,R)$
on $\operatorname{Defor}_{\textup{Alg}}(A,R)$.
Let $d\in (C^\bullet(A)[1]\otimes_km_{R})^0$ and $f \in \MC(C^\bullet(A)[1]\otimes_km_{R})$.
The automorphism $e^d:BA\otimes_kR\to BA\otimes_kR$ acts on $(b\otimes_kR+f)$
by $(b\otimes_kR+f)\mapsto e^d\circ (b\otimes_kR+f)\circ e^{-d}$.
Using $e^x\circ e^y\circ e^{-x}=e^{ad(x)}(y)$ and $(b\otimes_kR+f)^2=0$
we have
\begin{eqnarray*}
e^d\circ (b\otimes_kR+f)\circ e^{-d} &=& e^{[d,-]}(b\otimes_kR+f) \\
&=& \sum_{i\ge 0}\frac{[d,-]^i}{i!}(b\otimes_kR+f) \\
&=& b\otimes_kR+f+\sum_{i\ge 0}\frac{[d,-]^i}{(i+1)!}([d,b\otimes_kR]+[d,f]) \\
&=& b\otimes_kR+f+\sum_{i\ge 0}\frac{[d,-]^i}{(i+1)!}([d,f]-\PH\otimes_kR(d)) \\
&=& b\otimes_kR+e^d\bullet f
\end{eqnarray*}
where $e^d\bullet (-)$ indicates the (gauge) action of $d\in \exp(C^\bullet(A)[1]\otimes_km_{R})$ on $\MC(C^\bullet(A)[1]\otimes_km_{R})$,
see Section~\ref{premoduli}.
Hence the desired compatibility follows.
\QED

Proposition~\ref{compdeform} follows from the following Claims.
The proofs of them are analogous to those of Claim~\ref{deformclaim1} and~\ref{deformclaim2}, and are easier.

\begin{Claim}
\label{deformclaim3}
Let $E$ be a complex (i.e., a dg $k$-module with a differential $d$).
Let $\operatorname{Defor}_{\textup{com}}(E,R)$ be the set of deformations $\tilde{E}=(E\otimes_kR,\tilde{d})$ of $E$ to $R$.
There is a natural bijective map
\[
\MC(\End_{k}(E)\otimes_km_R)\stackrel{\sim}{\longrightarrow} \operatorname{Defor}_{\textup{com}}(E,R)
\]
which carries $f:E\to E\otimes_k m_R$ of degree one
satisfying the Maurer-Cartan equation
to a differential $d\otimes_kR+f\otimes R:E\otimes_kR\to E\otimes_kR$.
\end{Claim}

\Proof
As in Claim~\ref{deformclaim1},
our claim follows from the observation that
if $\iota:E\to E\otimes R$ is the natural inclusion,
a $R$-linear map $\tilde{d}:E\otimes R\to E\otimes R$ of degree $1$,
whose reduction is $d$, defines a deformation of $E$
if and only if
$(\tilde{d}-d\otimes R)\circ \iota:E\hookrightarrow E\otimes R\to E\otimes m_R$
is a Maurer-Cartan elememt in $\End_{k}(E)\otimes_km_R$.
\QED

\begin{Claim}
\label{deformclaim4}
Let $\operatorname{Aut}_{\textup{gr}}(E,R)$ be the group of automorphisms $E\otimes_kR\to E\otimes_kR$ of the graded $R$-module $E\otimes_kR$
whose reduction $E\stackrel{\sim}{\to}E$ is the identity.
There is a natural isomorphism of groups
\[
\exp(\End_k(E)\otimes_km_R)^0\stackrel{\sim}{\longrightarrow} \Aut_{\textup{gr}}(E,R)
\]
which carries $g:E\to E\otimes_k m_R$ of degree zero
to an automorphism $e^g:E\otimes_kR\to E\otimes_kR$.
Moreover, the gauge action of $\exp(\End_k(E)\otimes_km_R)^0$
on $\MC(\End_k(E)\otimes_km_R)$ commutes with the action of $\Aut_{\textup{gr}}(E,R)$
on $\operatorname{Defor}_{\textup{com}}(E,R)$ that is defined in
a similar way to that of $\Aut(A,R)$ on $\operatorname{Defor}_{\textup{Alg}}(A,R)$.
\end{Claim}

\Proof
Given $\phi:E\otimes R\to E\otimes R$ in
$\operatorname{Aut}_{\textup{gr}}(E,R)$,
we let $f:=\phi-\textup{id}:E\otimes R\to E\otimes m_R$.
Consider the $R$-linear map
$g:=\log(1+f)=\Sigma_{n=1}^\infty(-1)^{n-1}\frac{f^n}{n}:E\otimes R\to E\otimes m_R$ (it is a finite sum since $m_R$ is nilpotent).
Then we have $e^g=\phi$.
Conversely, if $h$ belongs to $\End_k(E)\otimes_km_R$,
then the exponential $e^{h\otimes R}:E\otimes R\to E\otimes R$ of $h\otimes R:E\otimes R\to E\otimes m_R$ defines an element in
$\operatorname{Aut}_{\textup{gr}}(E,R)$.
This correspondence yields the bijection between $\operatorname{Aut}_{\textup{gr}}(E,R)$ and $\End_k(E)\otimes_km_R$.
The second claim can be proved in a similar way to Claim~\ref{deformclaim2}.
\QED

\subsection{}
\label{premoduli2}
Let $\tilde{A}:=(A\otimes R,\tilde{b})$ be a deformation of $(A,b)$ to $R$.
Then it gives rise to
\begin{itemize}
\item the
Hochschild chain complex $C_\bullet(\tilde{A})$,
\item the negative cyclic complex $C_\bullet(\tilde{A})[[t]]$,
\item the periodic cyclic complex $C_\bullet(\tilde{A})((t))$.
\end{itemize}
Note that $C_\bullet(\tilde{A})$ is a dg $R$-module
(which is $(A\otimes_kB\bA)\otimes_kR$ as a graded $R$-module).
It follows from the definition of Hochschild chain complexes in Section~\ref{Hochschildchain} that the the complex $C_\bullet(\tilde{A})$ over $R$ is a deformation of the complex $C_\bullet(A)$.
Likewise, 
 $C_\bullet(\tilde{A})[[t]]$ is a dg $R[t]$-module, and the dg $R[t]$-module $C_\bullet(\tilde{A})[[t]]$ is a deformation of the dg $k[t]$-module
 $C_\bullet(A)[[t]]$ to $R$.
To be precise, by a deformation $M$
of the dg $k[t]$-module $C_\bullet(A)[[t]]$
(resp. the dg $k[t^\pm]$-module $C_\bullet(A)((t))$)
to $R$ we mean
a graded $R[t]$-module $M=(C_\bullet(A)\otimes_kR)[[t]]\simeq C_\bullet(A)[[t]]\otimes_kR$
(resp. a graded $R[t^\pm]$-module $M=(C_\bullet(A)\otimes_kR)((t))\simeq C_\bullet(A)((t))\otimes_kR$)
equipped with a differential $\tilde{\partial}$
 whose reduction $(C_\bullet(A)\otimes_kR/m_R)[[t]]\simeq C_\bullet(A)[[t]] \to (C_\bullet(A)\otimes_kR/m_R)[[t]]\simeq C_\bullet(A)[[t]]$ (resp. $(C_\bullet(A)\otimes_kR/m_R)((t))\to (C_\bullet(A)\otimes_kR/m_R)((t))$) is $\Ph+tB$.
An isomorphism $((C_\bullet(A)\otimes_kR)[[t]],\tilde{\partial}_1)\to ((C_\bullet(A)\otimes_kR)[[t]],\tilde{\partial}_2)$  of deformations
of the dg $k[t]$-module $C_\bullet(A)[[t]]$ to $R$
is an isomorphism of dg $R[t]$-modules whose reduction
is the identity. An isomorphism of deformations of the dg $k[t^\pm]$-module $C_\bullet(A)((t))$ is defined in a similar way.

Let $\mathsf{D}_{C_{\bullet}(A)[[t]]}$ denote a functor
$\Art\to \Set$ which carries $R$ to the set of isomorphism classes of
deformations of dg $k[t]$-module $C_{\bullet}(A)[[t]]$ to $R$.
Let $\mathsf{D}_{C_{\bullet}(A)((t))}$ denote a functor
$\Art\to \Set$ which carries $R$ to the set of isomorphism classes of
deformations of dg $k[t^\pm]$-module $C_{\bullet}(A)((t))$ to $R$.
The following is a version of Proposition~\ref{compdeform},
which follows from formal extensions of Claim~\ref{deformclaim3} and Claim~\ref{deformclaim4}.

\begin{Proposition}
\label{compdeform2}
There is a natural equivalence of functors $\Def_{\End_{k[t]}(C_\bullet(A)[[t]])}\to \mathsf{D}_{C_\bullet(A)[[t]]}$.
It carries a Maurer-Cartan element $x$ in $\End_{k[t]}(C_\bullet(A)[[t]])\otimes m_R$
to a deformation $((C_\bullet(A)\otimes_kR)[[t]],(\Ph+tB)\otimes_k R+x\otimes R)$.
% where the differential is given by$((\Ph+tB)\otimes_kR+x)(e\otimes r)=(\Ph+tB)(e)\otimes r+x(e)r$.
%If one formally replaces $[[t]]$ by $((t))$, the same assertion holds for $\mathsf{D}_{C_\bullet(A)((t))}$.
\end{Proposition}

An isomorphism between two deformations $\tilde{A}_1=(A\otimes_kR,\tilde{b}_1)\stackrel{\sim}{\to} \tilde{A}_2=(A\otimes_kR,\tilde{b}_2)$ of $(A,b)$
induces an isomorphism of deformations $C_\bullet(\tilde{A}_1)\stackrel{\sim}{\to} C_\bullet(\tilde{A}_2)$ of the complex $C_\bullet(A)$,
and an isomorphism of deformations $C_\bullet(\tilde{A}_1)[[t]]\stackrel{\sim}{\to} C_\bullet(\tilde{A}_2)[[t]]$ of the dg $k[t]$-module $C_\bullet(A)[[t]]$
to $R$. Consequently, we obtain a natural transformation of functors
\[
\mathcal{P}:\mathsf{DAlg}_{A}\to \mathsf{D}_{C_{\bullet}(A)}
\]
which sends the isomorphism class of a deformation of $\tilde{A}$ of $(A,b)$ to
the isomorphism class of
the deformation $C_\bullet(\tilde{A})$ of $C_\bullet(A)$
for each $R$.
Similarly, we define functors
\[
\mathcal{Q}:\mathsf{DAlg}_{A}\to \mathsf{D}_{C_{\bullet}(A)[[t]]}\ \ \ \textup{and}\ \ \ \mathcal{R}:\mathsf{DAlg}_{A}\to \mathsf{D}_{C_{\bullet}(A)((t))}
\]
which carry a deformation of $\tilde{A}$ of $(A,b)$ to
the deformation $C_\bullet(\tilde{A})[[t]]$ of $C_\bullet(A)[[t]]$
and the deformation $C_\bullet(\tilde{A})((t))$ of $C_\bullet(A)((t))$
respectively for each $R$.

\begin{Proposition}
\label{perioddeform1}
Through the equivalences $\mathsf{DAlg}_{A}\simeq \Def_{C^\bullet(A)[1]}$, $\mathsf{D}_{C_\bullet(A)} \simeq \Def_{\End_k(C_\bullet(A))}$,
$\mathsf{D}_{C_\bullet(A)[[t]]}\simeq \Def_{\End_{k[t]}(C_\bullet(A)[[t]])}$,
and $\mathsf{D}_{C_\bullet(A)((t))}\simeq \Def_{\End_{k[t^\pm]}(C_\bullet(A)((t)))}$,
the functors $\mathcal{P}$, $\mathcal{Q}$ and $\mathcal{R}$ can be identified
with the functors
\begin{eqnarray*}
&&\Def_L:\Def_{C^\bullet(A)[1]}\to \Def_{\End_k(C_\bullet(A))}, \\
&&\Def_{L[[t]]}:\Def_{C^\bullet(A)[1]}\to \Def_{\End_{k[t]}(C_\bullet(A)[[t]])},
 \\
&&\Def_{L((t))}:\Def_{C^\bullet(A)[1]}\to \Def_{\End_{k[t^\pm]}(C_\bullet(A)((t)))},
\end{eqnarray*}
associated to $L$, $L[[t]]$ and $L((t))$ respectively.
\end{Proposition}

We prove this Proposition.
First, we find an explicit formula of the differential on
Hochschild chain complex $C_\bullet(A)$ (see Section~\ref{Hochschildchain}).
Let $(A,b)$ be a unital curved $\AI$-algebra.
We denote by $b_l:(A[1])^{\otimes l}\to A[1]$ the $l$-th component of $b$.

\begin{Lemma}
\label{explicitHoch}
The differential on the Hochschild chain complex $C_\bullet(A)$ of $A$
is determined by the formula
\begin{align*}
&\partial_{\textup{Hoch}}(a_0\otimes[a_1|\ldots|a_n]) \\
&\ \ =\sum_{0\le i\le n \atop b_l\textup{includes}\ a}(-1)^{(\epsilon_i+|a_0|-1)(\epsilon_n-\epsilon_i)}s^{-1}b_l[a_{i+1}|\ldots|a_n|a_0|a_{1}|\ldots|a_{l-n+i-1}]\otimes [a_{l-n+i}|\ldots|a_i] \\
&\ \ \ \ \ +\sum_{j+l\le n} (-1)^{|a_0|-1+\epsilon_j} a_0\otimes [a_1|\ldots|b_{l}[a_{j+1}|\ldots|a_{j+l}]|a_{l+j+1}|\ldots|a_n].
\end{align*}
\end{Lemma}

\Proof
We have an isomorphism $A\otimes B\bA\simeq \Phi(A)\subset B\bA\otimes A\otimes B\bA$ given by $S_{312}\circ (1\otimes \Delta_{B\bA})$ (see Section~\ref{Hochschildchain}). In explicit terms,
it carries $a\otimes [a_1|\ldots|a_n]$ to
\[
\sum_{i=0}^{n} (-1)^{(\epsilon_i+|a|)(\epsilon_n-\epsilon_i)}[a_{i+1}|\ldots|a_n]\otimes a\otimes [a_1|\ldots |a_i]
\]
where $\epsilon_i=\Sigma_{l=1}^{i}(|a_l|-1)$.
The counit map $\eta:B\bA\to k$ induces the
inverse isomorphism $\Phi(A)\subset B\bA\otimes A\otimes B\bA\stackrel{\eta\otimes 1\otimes 1}{\longrightarrow} A\otimes B\bA$.
Let $d$ be the differential on $B\bA\otimes A\otimes B\bA$
given by $\{d_{m,n}\}_{m,n\ge0}$ in Section~\ref{Hochschildchain}: 
\begin{align*}
 &d([a_1|\ldots|a_i]\otimes a\otimes [a'_1|\ldots|a'_j]) \\
 &\ \ =  \sum_{s+l\le i}(-1)^{\epsilon_s}[a_1|\ldots|a_s|b_l[a_{s+1}|\ldots|a_{s+l}]|\ldots |a_i] \otimes a \otimes [a_1'|\ldots|a_j'] \\
 &\ \ \  +
\sum_{s<l,i<s+l}(-1)^{\epsilon_s}(-1)^{\epsilon_i-\epsilon_s}[a_1|\ldots|a_s]\otimes s^{-1}b_l[a_{s+1}|\ldots |a_i|a|a_1'|\ldots |a'_{l-i+s-1}]\otimes [a'_{l-i+s}|\ldots |a'_j] \\
&\ \ \  +\sum_{t+l\le j}(-1)^{\epsilon_i+|a|}(-1)^{\epsilon'_t}[a_1|\ldots|a_i]\otimes a\otimes [a_1'|\ldots|a_t'|b_l[a_{t+1}'|\ldots|a_{t+l}']|\ldots |a'_j]
\end{align*}
where $b_l$ is the $l$-th component of curved $A_\infty$-structure of $\bar{A}$
(we abuse notation),
and $\epsilon'_t=\Sigma_{l=1}^t(|a_l'|-1)$.
Taking account of the above two formulas
we obtain the following formula of the differential on $C_\bullet(A)=A\otimes B\bA$
\begin{align*}
&\partial_{\textup{Hoch}}(a_0\otimes[a_1|\ldots|a_n])=(\eta\otimes 1\otimes 1)\circ d \circ (S_{312}\circ 1\otimes \Delta_{B\bA}) (a_0\otimes[a_1|\ldots|a_n]) \\
&\ =\sum_{0\le i\le n \atop b_l\textup{includes}\ a}(-1)^{(\epsilon_i+|a_0|)(\epsilon_n-\epsilon_i)+(\epsilon_n-\epsilon_i)}s^{-1}b_l[a_{i+1}|\ldots|a_n|a_0|a_{1}|\ldots|a_{l-n+i-1}]\otimes [a_{l-n+i}|\ldots|a_i] \\
&\ \ \ \ +\sum_{j+l\le n} (-1)^{1+|a_0|+\epsilon_i} a_0\otimes [a_1|\ldots|b_{l}[a_{j+1}|\ldots|a_{j+l}]|a_{l+j+1}|\ldots|a_n].
\end{align*}
\QED

{\it Proof of Proposition~\ref{perioddeform1}}.
In this proof, we will prove a more refined statement.
Let $R$ be an artin local $k$-algebra in $\Art$ and $m_R$ its maximal ideal.
Let $\MC_L:\MC(C^\bullet(A)[1]\otimes_km_R)\to \MC(\End_k(C_\bullet(A))\otimes_km_R)$
be the map induced by $L$.
Claim~\ref{deformclaim1} (resp. Claim~\ref{deformclaim3})
describes a canonical bijective correspondence between
$\MC(C^\bullet(A)[1]\otimes_km_R)$ and the set of deformations of $(A,b)$ to $R$ (resp. between $\MC(\End_k(C_\bullet(A))\otimes_km_R)$ and the set of deformations of $C_\bullet(A)$ to $R$).
Using these correspondences we can identify $\MC_L$ with the map
from the set of deformations of the $\AI$-algebra $(A,b)$ to that of
the complex $C_\bullet(A)$, which carries $\tilde{A}=(A\otimes_kR,\tilde{b})$
to $C_\bullet(\tilde{A})$ as follows.
Indeed, if we put $\tilde{b}=b\otimes_kR+f\otimes R$ with $f\in \MC(C^\bullet(A)[1]\otimes_km_R)$, then the differential on $C_\bullet(\tilde{A})=(A\otimes_kR)\otimes_R(B\bar{A}\otimes_kR)$ is given by the formula in Lemma~\ref{explicitHoch}
where each $b_l$ is replaced by $\tilde{b}_l=b_l\otimes_kR+f_l\otimes R$
where $f_l:(A[1])^{\otimes l}\to A[1]\otimes m_R$ is the $l$-th term of $f$.
Observe that it coincides with $\Ph\otimes_kR+L_f\otimes R$ where $L_f:C_\bullet(A)\to C_\bullet(A)\otimes_km_R$ is defined in Section~\ref{Lieaction}!

Next through correspondences,
 the computation similar to Lemma~\ref{explicitHoch}
shows that $(C^\bullet(A)[1]\otimes_km_R)^0\to (\End_k(C_\bullet(A))\otimes_km_R)^0$ induced by $L$ is identified with
$\Aut(A,R)\to \Aut_{\textup{gr}}(C_\bullet(A),R)$ which carries $BA\otimes_kR\stackrel{\sim}{\to} BA\otimes_kR$ (an automorphism of the graded coalgebras) to $C_\bullet(A)\otimes_kR\stackrel{\sim}{\to} C_\bullet(A)\otimes_kR$ (the induced automorphism of the graded $R$-module).

Finally,
\begin{eqnarray*}
\mathcal{P}:\mathsf{DAlg}_{A}(R)&\simeq& \MC(C^\bullet(A)[1]\otimes_km_R)/\exp(C^\bullet(A)[1]\otimes_km_R)^0 \\ &\to& \MC(\End_k(C_\bullet(A))\otimes_km_R)/\exp(\End_k(C_\bullet(A))\otimes_km_R)^0 \simeq \mathsf{D}_{C_\bullet(A)}(R)
\end{eqnarray*}
is naturally isomorphic to $\Def_L$. The cases of $\mathcal{Q}$ and
$\mathcal{R}$ are similar.
\QED

\section{Period mapping and homotopy calculus}
\label{construction}

\subsection{}
Let $A$ be a unital dg algebra over a field $k$ of characteristic zero.
Henceforth dg algebras are assumed to be unital in this paper.
In Section~\ref{LM} we considered a morphism $\mathcal{P}:\mathsf{DAlg}_{A}\to \mathsf{D}_{C_{\bullet}(A)}$ of functors, and its cyclic versions
$\mathcal{Q}:\mathsf{DAlg}_{A}\to \mathsf{D}_{C_{\bullet}(A)[[t]]}$,
$\mathcal{R}:\mathsf{DAlg}_{A}\to \mathsf{D}_{C_{\bullet}(A)((t))}$.
Proposition~\ref{perioddeform1}
allows us to consider them as the functors associated to homomorphism
$L,L[[t]],L((t))$ of dg Lie algebras. In particular, $\mathcal{R}$
can be identified with $\DEF_{L((t))}:\DEF_{C^\bullet(A)[1]}\to \DEF_{\End_{k[t^\pm]}(C_\bullet(A))}$. 
The purpose of this section is to construct a period mapping
for deformations of a dg algebra which can be described
in a moduli-theoretic fashion.
We first prove the following:

\begin{Proposition}
\label{triviality}
The morphism $L((t)):C^\bullet(A)[1]\to \End_{k[t^\pm]}(C_\bullet(A)((t)))$ of 
dg Lie algebras is null homotopic. See Section~\ref{mappingsp} for the mapping
space between dg Lie algebras.
\end{Proposition}

\begin{Remark}
Proposition~\ref{triviality} especially means that any deformation of $A$ induces no non-trivial deformation
of the periodic cyclic complex $C_\bullet(A)((t))$ at any level (including
higher homotopy and derived structure).
If we think of $C_\bullet(A)((t))$ as a ``topological invariant'' of $A$
(keep in mind that Hochschild-Kostant-Rosenberg theorem),
this is an analogue of Ehresmann's fibration theorem.
%If one replaces $((t))$ by $[[t]]$ in Theorem~\ref{triviality}, then an analogues result does not hold in general.
\end{Remark}

The proof of Proposition~\ref{triviality}
is completed in the end of Section~\ref{trivialityend}.

\subsection{}
\label{twomore}
To prove Proposition~\ref{triviality}, we need two more algebraic
structures on $(C^\bullet(A),C_\bullet(A))$
(see \cite{TT}).
We briefly review them.
We adopt notation in Section~\ref{Hochscochaindef}.
Recall that $C^\bullet(A)=C^\bullet(A)[1][-1]=\prod_{l\ge0}\Hom_k((\bar{A}[1])^{\otimes l},A)$
% \simeq \prod_{l\ge0}\Hom_k(\bar{A}^{\otimes l},A)$
where $\bar{A}=A/k$.
% and the last map is induced by$\Pi s^{\otimes l}$.
Let $P\in \Hom_k((\bar{A}[1])^{\otimes p},A)\simeq \Hom_k((\bar{A}[1])^{\otimes p},A[1])$
and $Q\in \Hom_k((\bar{A}[1])^{\otimes q},A)\simeq \Hom_k((\bar{A}[1])^{\otimes p},A[1])$ be two elements in $C^\bullet(A)$ (the isomorphism is given by $s:A\to A[1]$). 
The tensor product of maps, and compositions with the 2-nd component $b_2:A[1]\otimes A[1]\to A[1]$
(of the $\AI$-structure) and $s^{-1}:A[1]\to A$ induce
\[
\Hom_k((\bar{A}[1])^{\otimes p},A[1])\otimes \Hom_k((\bar{A}[1])^{\otimes q},A[1])\to \Hom_k((\bar{A}[1])^{\otimes p+q},A[1])\to \Hom_k((\bar{A}[1])^{\otimes p+q},A).
\]
We define the cup product $P\cup Q$ to be the image of this composite.
The cup product exhibits $C^\bullet(A)$ as a (non-commutative) dg algebra
with unit $1\in \Hom(k,A)$.

For $P$ in $\Hom_k((\bar{A}[1])^{\otimes p},A)\subset C^\bullet(A)$ we define a contraction map $I_P:C_\bullet(A)\to C_\bullet(A)$ of degree $|P|$ as follows.
\begin{eqnarray*}
\Hom_k((\bar{A}[1])^{\otimes p},A)\otimes (A\otimes (\bar{A}[1])^{\otimes n}) &\stackrel{s\otimes 1^{\otimes n+1}}{\longrightarrow}& \Hom_k((\bar{A}[1])^{\otimes p},A[1])\otimes (A\otimes (\bar{A}[1])^{\otimes n}) \\
&\to&  A\otimes A[1]\otimes (\bar{A}[1])^{\otimes n-p} \\
&\stackrel{s\otimes 1}{\to}& A[1]\otimes A[1]\otimes (\bar{A}[1])^{\otimes n-p} \\
&\stackrel{b_2\otimes 1}{\to}& A[1]\otimes (\bar{A}[1])^{\otimes n-p} \\
 &\stackrel{s^{-1}\otimes 1}{\to}& A \otimes (\bar{A}[1])^{\otimes n-p}
\end{eqnarray*}
where the second arrow is given by the pairing of
$\Hom((\bar{A}[1])^{\otimes p},A[1])$ and the first $p$ factors
in $(\bar{A}[1])^{\otimes n}$ for $p\le n$
 (if otherwise it is defined to be zero).
We denote the induced map by
$I:C^\bullet(A)\to \End_k(C_\bullet(A)),\ P\mapsto I_P$.
Explicitly, $I_P(a_0\otimes [a_1|\ldots|a_n])$
has the form $\pm a_0P[a_1|\ldots|a_p]\otimes[a_{p+1}|\ldots |a_n]$
for $p\le n$.

The five operations $\cup$, the bracket $[-,-]_G$ on $C^\bullet(A)[1]$,
the dg Lie algebra map $L$, the contraction $I$, and Connes' operator $B$
constitute a calculus structure on $(HH^\bullet(A),HH_\bullet(A))$
(cf. \cite{TT}).
Here is the definition of calculus.

\begin{Definition}
\label{calcdef}
\begin{enumerate}
\renewcommand{\labelenumi}{(\theenumi)}

\item A graded $k$-module $V$ is a Gerstenharber algebra if it
is endowed with a (graded) commutative and associative product $\wedge:V\otimes V\to V$
of degree zero and a Lie bracket $[-,-]:V\otimes V\to V$
of degree $-1$ which exhibits $V[1]$ as a graded Lie algebra.
These satisfy the Leibniz rule
\[
[a,b\wedge c]=[a,b]\wedge c+(-1)^{(|a|+1)|b|}b\wedge[a,c]
\]
where $a,b,c$ are homogeneous elements of $V$.

\item A precalculus structure is a pair of a Gerstenharbar algebra
$(V,\wedge,[-,-])$ and a graded $k$-module $W$ together with
a module structure 
\[
i:V\otimes W\to W
\]
of the commutative algebra $(V,\wedge)$, and a module structure
\[
l:V[1]\otimes W\to W
\]
of (graded) Lie algebra $V[1]$ such that
\[
i_al_b-(-1)^{|a|(|b|-1)}l_bi_a=i_{[a,b]},\ \ \ \textup{and}\ \ \ l_{a\wedge b}=l_ai_b+(-1)^{|a|}i_al_b
\]
where $i_a:W\to W$ is determined by $i(a\otimes(-))$ for $a\in V$,
and $l_b:W\to W$ is determined by $l(b\otimes(-))$ for $b\in V$.

\item A calculus is a precalculus $(V,W,\wedge,[-,-],i,l)$ 
endowed with a linear map $\delta:W\to W$ of degree $-1$ such that
\[
\delta^2=0\ \ \ \textup{and}\ \ \delta i_a-(-1)^{|a|}i_a\delta=l_a.
\]
We call the 7-tuple $(V,W,\wedge,[-,-],i,l,\delta)$ a calculus algebra
(or we also refer to it as a structure of calculus on $(V,W)$).

\item A 5-tuple
$(V,W,[-,-],l,\delta)$ is said to be a $\textup{Lie}^\dagger$-algebra
when
$[-,-]$ determines a Lie bracket on $V[1]$ as in (1), 
$l:C[1]\otimes W\to W$ is an action of $C[1]$ on $W$,
$\delta:W\to W$ is a linear map of degree $-1$ such that
\[
\delta^2=0\ \ \ \textup{and}\ \ \delta l_a-(-1)^{|a|-1}l_a\delta=0.
\]

\end{enumerate}

\end{Definition}

\subsection{}
According to the work of
Daletski, Gelfand, Tamarkin and Tsygan (see \cite{TT},
\cite{DTT1}, \cite{DTT2}),
$(\cup, [-,-]_G,I,L,B)$ determines a structure of calculus
$(V=HH^*(A),W=HH_*(A))$, that is,
\[
(HH^*(A),HH_*(A),\cup, [-,-]_G,I,L,B)
\]
is an important example of a calculus algebra.
The operations
$(\cup, [-,-]_G,I,L,B)$ fail to
determine a structure of calculus on $(C^\bullet(A),C_\bullet(A))$.
For example, $\cup$ and $[-,-]_G$ do not satisfy the Leibniz rule.
To endow $(C^\bullet(A),C_\bullet(A))$ with algebraic structures
coming from $(\cup, [-,-]_G,I,L,B)$ in a suitable way,
one needs the machinery of operads.

Kontsevich and Soibelman constructed a 2-colored dg operad $\mathcal{KS}$
(consisting of two colors) and a natural action of $\mathcal{KS}$
on the pair $(C^\bullet(A),C_\bullet(A))$ for an $\AI$-algebra $A$,
see \cite[11.1, 11.2, 11.3]{KS}.
This action of the operad $\mathcal{KS}$ generalizes
the solution of Deligne's conjecture on
an action of the little disks operad on $C^\bullet(A)$.
We also refer to
Dolgushev-Tamarkin-Tsygan
\cite[Section 4]{DTT1} for the detailed study on $\mathcal{KS}$
in the case of (dg) algebras.
The recent work of Horel
\cite{Hor} treats a generalization
to the case of ring spectra, which is based on the factorization homology
and Swess-cheese operad conjecture.
The operad $\mathcal{KS}$ is closely related to calculi.
Let $\mathbf{calc}$ be the 2-colored graded operad defined by
generators $\wedge, [-,-], i,l,\delta$ and their relations
in Definition~\ref{calcdef}
(we refer the reader to \cite{LV} for dg operads,
\cite{Hin}, \cite[Section 3.5, 3.6]{TT} for colored operads,
modules over an algebra over an
operad, and \cite{DTT1} for $\mathbf{calc}$).
The operad $\mathbf{calc}$ has the suboperad called Gerstenhaber operad,
that is generated by $\wedge$ and $[-,-]$ and relations in
Definition~\ref{calcdef} (1), see \cite[13.3.12]{LV} for the Gerstenhaber operad.
Let $\mathbf{Lie}^\dagger$ be the 2-colored graded suboperad
of $\mathbf{calc}$ generated by
$[-,-],l,\delta$ and their relations
in Definition~\ref{calcdef} (4).
A $\mathbf{Lie}^\dagger$-algebra is a $\textup{Lie}^\dagger$-algebra.
According to Proposition~\ref{Laction}, $(C^\bullet(A),C_\bullet(A),[-,-],L,B)$ is a $\mathbf{Lie}^\dagger$-algebra in the category of complexes.
Thanks to \cite[Thereom 2, Threorem 3, Proposition 8]{DTT1}, the homology
operad $H_*(\mathcal{KS})$ is quasi-isomorphic to $\mathbf{calc}$.
Moreover, $\mathcal{KS}$ is formal, namely, it is quasi-isomorphic to $H_*(\mathcal{KS})\simeq \mathbf{calc}$.
Let $\mathcal{C}alc$ denote a natural cofibrant replacement by the cobar-bar
resolution
$\mathcal{C}alc:=\textup{Cobar}(\textup{Bar}(\mathbf{calc})) \to \mathbf{calc}$.
There is a quasi-isomorphism $\mathcal{C}alc\to \mathcal{KS}$.
Thus, we obtain a $\mathcal{C}alc$-algebra $(C^\bullet(A),C_\bullet(A))$.
We summarize the result as

\begin{Proposition}
\label{goodcalculus}
There exists a $\mathcal{C}alc$-algebra $(C^\bullet(A),C_\bullet(A))$,
that is, an action of $\mathcal{C}alc$ on $(C^\bullet(A),C_\bullet(A))$
whose underlying $H_*(\mathcal{C}alc)\simeq \mathbf{calc}$-algebra
is the calculus $(HH^*(A),HH_*(A))$.
Moreover, the operad $\mathcal{C}alc$ has a suboperad
$\mathcal{L}ie^\dagger$ endowed with a quasi-isomorphism
$\mathcal{L}ie^\dagger\to \mathbf{Lie}^\dagger$ which commutes with
the counit quasi-isomorphism  $\mathcal{C}alc\to \mathbf{calc}$,
such that the pullback of the $\mathcal{C}alc$-algebra $(C^\bullet(A),C_\bullet(A))$
to $\mathcal{L}ie^\dagger$ coincides with
the pullback of the $\mathbf{Lie}^\dagger$-algebra
 $(C^\bullet(A),C_\bullet(A),[-,-],L,B)$ to $\mathcal{L}ie^\dagger$.
\end{Proposition}

\Proof
The statement of the first half is \cite[Theorem 2, 3, Proposition 8]{DTT1}.
Namely, we can obtain an action of the operad $\mathcal{C}alc$
on the pair of complexes
$(C^\bullet(A),C_\bullet(A))$ from that of $\mathcal{KS}$
whose underlying action of $H_*(\mathcal{KS})\simeq H_*(\mathcal{C}alc)\simeq \mathbf{calc}$ on $(HH^*(A),HH_*(A))$
can be identified with the above calculus $(HH^*(A),HH_*(A))$.
According to \cite[Theorem 4]{DTT1} and its proof,
the $\mathcal{C}alc$-algebra $(C^\bullet(A),C_\bullet(A))$
is quasi-isomorphic to another action of $\mathcal{C}alc$ on $(C^\bullet(A),C_\bullet(A))$ whose restriction to $\mathcal{L}ie^\dagger$ comes from
the $\mathbf{Lie}^\dagger$-algebra
given by operations $L$, $[-,-]$, $B$ on $(C^\bullet(A),C_\bullet(A))$.
\QED

\subsection{}
\label{extendlie}
Let $(C,M,[-,-], l,\delta)$ be a $\mathbf{Lie}^\dagger$-algebra in the
symmetric monoidal category of complexes of $k$-modules
(in particular, $C,C', M$ and $M'$ are complexes).
In explicit terms,
a $\mathbf{Lie}^\dagger$-algebra amounts to data consisting of
(i) a bracket $C[1]\otimes C[1]\to C[1]$ of degree 0 which exhibits $C[1]$
as a dg Lie algebra, (ii) $\delta:M\to M$ of degree $-1$ with $\delta^2=0$
and $d_M \delta+\delta d_M=0$,
(iii) an action of $l:C[1]\otimes M\to M$ of the dg Lie algebra $C[1]$
on the complex $M$
which commutes with $\delta$ (i.e., $\delta l_c-(-1)^{|c|-1}l_c \delta=0$
for $c\in C$). Here $d_M$ is the differential of $M$.
The pair $(M,\delta)$ is a mixed complex.
Thus it gives rise to a dg $k[t^\pm]$-module $(M((t)), d_M+t\delta)$.
Roughly speaking, we may say that a $\mathbf{Lie}^\dagger$-algebra
is a dg Lie algebra $C[1]$ together with its action on a mixed complex $M$.
To a $\mathbf{Lie}^\dagger$-algebra $(C,M,[-,-], l,\delta)$
we associate a dg Lie algebra map
\[
l((t)):C[1]\to \End_{k[t^\pm]}(M((t)))
\]
which carries $c$ to $l_c((t)):M((t))\to M((t)),\ mt^n\mapsto l_c(m)t^n$.

\subsection{}
Henceforce we will use some model category structures.
Appropriate references are \cite{Hov}, \cite[Appendix]{HTT}.
We emply
the model category structure of the category of
algebras over a dg operad \cite[2.3.1, 2.6.1]{Hin}.
The category of $\mathbf{Lie}^\dagger$-algebras denoted by $\Alg_{\mathbf{Lie}^\dagger}$ admits a combinatorial model category structure
where a morphism $(C,M,[-,-], l,\delta)\to (C',M',[-,-]', l',\delta')$
is a weak equivalence (resp. a fibration) if both $C\to C'$ and $M\to M'$ are quasi-isomorphisms (resp. degreewise surjective maps) of complexes.
Similarly, the category of dg Lie algebras (resp. dg $k[t^\pm]$-modules)
admits
a combinatorial model category structure
where a morphism
is a weak equivalence (resp. a fibration) if the underlying map of complexes
is a quasi-isomorphism (resp. a degreewise surjective map).
Note that every object is fibrant.

\subsection{}
\label{mappingsp}
Let us recall the mapping space between two dg Lie algebras,
cf. \cite{Hin0}.
Let $V$ and $W$ be two dg Lie algebras.
Let $B_{com}V$ be the bar construction associated to $V$,
which is a unital
cocommutative dg coalgebra over $k$.
Here coalgebras are assumed to be conilpotent.
Explicitly,
$B_{com}V$ is $\oplus_{i\ge 0}\operatorname{Sym}^n(V[1])$ as a graded cocommutative coalgebra. Here $\Sym^n$ indicates the $n$-fold symmetric product, and
$B_{com}(-)$ is different from $B(-)$ in the previous sections.
 The differential is the sum of two differentials; the first comes from the
 differential of $V$, the second is determined by $[-,-]:(V\wedge V)[2]\simeq V[1]\wedge V[1] \to V[1][1]$. It gives rise to a functor
$B_{com}:\textup{dgLie}\to \textup{dgcoAlg}$ from the category of dg Lie algebras
to the category of unital cocommutaitve dg coalgebras.
There is a left adjoint $Cob_{com}:\textup{dgcoAlg} \to \textup{dgLie}$
of $B_{com}$ given by the cobar construction (see e.g. \cite[2.2.1]{Hin0}).
For any dg Lie algebra $V$, the counit map $Cob_{com}B_{com}V\to V$
gives a (canonical) cofibrant replacement of $V$.
Let $\Omega_n$ denote the commutative
dg algbera of polynomial differential forms on the standard $n$-simplex.
Namely, it is
\[
\Omega_n:=k[u_0,\ldots u_n,du_0,\ldots du_n]/(\Sigma_{i=0}^{n}u_i-1, \Sigma_{i=0}^n du_i)
\]
where $k[u_0,\ldots u_n,du_0,\ldots du_n]$
is the free commutative graded algebra generated by $u_0,\ldots, u_n$
and  $du_0,\ldots, du_n$ with $|u_i|=0$, $|du_i|=1$ for each $i$,
and the differential carries $u_i$ to $du_i$ (see e.g. \cite{BG}).
If we consider
the family $\Omega_\bullet=\{\Omega_n \}_{n\ge0}$ of
commutative dg algebras, they form
a simplicial commutative dg algebras in the natural way.
Using the simplicial commutative dg algebra $\Omega_\bullet$ and
the simplicial model category structure of $\textup{dgLie}$
\cite[2.4]{Hin0}
we obtain a Kan complex $\Hom_{\textup{dgLie}}(Cob_{com}B_{com}V, \Omega_\bullet\otimes W)$,
that is a model of the mapping space.
In \cite{Hin0}, the definition of simplicial model categories is
slightly weaker than the standard one. But by \cite[1.4.2]{Hin2}
this Kan complex is naturally homotopy equivalent to the Hom simplicial set
in the simplicial category associated to the underlying model category via Dwyer-Kan hammock localization.
We prefer to work with another presentation of this model.
Let $C$ be a
dg coalgebra with a comultiplication $\Delta:C\to C\otimes C$
and $V$ a dg Lie algebra.
Let $\Hom(C,V)$ be the Hom complex between the underlying complexes of $C$
and $V$. It is endowed with the following convolution Lie bracket:
\[
[f,g]:C\stackrel{\Delta}{\to} C\otimes C\stackrel{f\otimes g}{\longrightarrow} V\otimes V\stackrel{[-,-]}{\longrightarrow} V
\]
for $f,g\in \Hom(C,V)$.
For the dg Lie algebra $\Hom(C,V)$ we define
the set of Maurer-Cartan elements $\MC(C,V):=\textup{MC}(\Hom(C,V))$.
According to \cite[2.2.5]{Hin0} there is a natural isomorphism
\[
\Hom_{\textup{dgLie}}(Cob_{com}B_{com}V, \Omega_n \otimes W)\simeq 
\MC(\overline{B_{com}V}, \Omega_n \otimes W)
\]
where $\overline{B_{com}V}$ is the kernel of the counit $B_{com}V\to k$.
(The isomorphism is induced by the composition with the inclusion
$\overline{B_{com}V}[-1]\hookrightarrow Cob_{com}B_{com}V$.)
In particular,
\[
\Map(V,W):=\Hom_{\textup{dgLie}}(Cob_{com}B_{com}V, \Omega_\bullet \otimes W)\simeq  \MC(\overline{B_{com}V},\Omega_\bullet \otimes W).
\]
In explicit terms, an element of
the 0-th term $\MC(\overline{B_{com}V},W)$ of $\MC(\overline{B_{com}V},\Omega_\bullet\otimes W)$
corresponds to a family of linear maps $\operatorname{Sym}^n(V[1]) \simeq (\wedge^nV)[n] \to W[1]$
of degree zero ($n\ge1$) that satisfies a certain relation coming from
the Maure-Cartan equation. It is sometimes called an $L_\infty$-morphism
in the literature. Thus we refer to an element of $\MC(\overline{B_{com}V},W)$ as an $L_\infty$-morphism. When $V[1]=\operatorname{Sym}^1(V[1]) \to W[1]$ is a quasi-isomorphism, we shall call it
an $L_\infty$-quasi-isomorphism.
Any equivalence class of a morphism from $V$ to $W$
of two dg Lie algebras can be represented by an $L_\infty$-morphism.
Since $\Map(V,W)$
is a Kan complex, for two $L_\infty$-morphisms $f$ and $g$
corresponding
to two vertices of this Kan complex, the space of homotopies/morphisms
from $f$ to $g$ makes sense (cf. \cite[1.2.2]{HTT}).

Put $\Omega_1=k[u,du]$. We write $W[u,du]$ for $k[u,du]\otimes W$ (do not confuse it with $W\otimes k[u,du]$, that gives rise to the different sign rule).
The degree of $du$ is $1$.
An element of $\MC(\overline{B_{com}V},W[u,du])$
represents a morphism (or a homotopy) between two $L_\infty$-morphisms.
Face maps $d_0,d_1:\MC(\overline{B_{com}V},W[u,du])\rightrightarrows \MC(\overline{B_{com}V},W)$ are induced by the composition with maps $p_0$
and $p_1$ of dg Lie algebras:
\[
p_0:W[u,du]\to W,\ \ \textup{and}\  p_1:W[u,du]\to W
\]
given by $p_0(u)=p_0(du)=0$, and $p_1(u)=1,\ p_1(du)=0$.
Thus, a homotopy/morphism from an $L_\infty$-morphism $f$ to an $L_\infty$-morphism $g$
can be represented by
$\phi \in \MC(\overline{B_{com}V},W[u,du])$
such that $d_0(\phi)=f$ and $d_1(\phi)=g$.
We say that $f$ is equivalent to $g$ if a morphism exists between them.

\subsection{}

\begin{Lemma}
\label{mapfunct}
Let $(C,M,[-,-], l,\delta)$ and $(C',M',[-,-]', l',\delta')$ be two $\mathbf{Lie}^\dagger$-algebras in the
category of complexes of $k$-modules.
Let $f:(C,M,[-,-], l,\delta)\to(C',M',[-,-]', l',\delta')$
be a weak equivalence of $\mathbf{Lie}^\dagger$-algebras.
If $f$ is a trivial fibration (i.e., both underlying maps $C\to C'$ and $M\to M'$
are surjective quasi-isomorphisms), then there exist a dg Lie algebra $F$ 
and weak equivalences
\[
\End_{k[t^\pm]}(M((t))) \leftarrow F \rightarrow \End_{k[t^\pm]}(M'((t)))
\]
such that the left arrow is injective, and the right arrow is surjective.

Moreover, there exists a homotopy equivalence between the mapping space from $L((t)):C[1]\to \End_{k[t^\pm]}(M((t)))$ to the zero morphism (see Section~\ref{mappingsp})
and the mapping space from $L'((t)):C'[1]\to \End_{k[t^\pm]}(M'((t)))$ to the zero morphism
via the zig-zag (see the proof for the construction).
\end{Lemma}

\Proof
We suppose that $f$ is a trivial fibration.
It follows that the induced morphism $M((t))\to M'((t))$
is a trivial fibration of dg $k[t^\pm]$-modules, i.e., a surjective quasi-isomorphism. Note also that every dg $k[t^\pm]$-module is cofibrant
with respect to the projective model structure.
Consider the pullback diagram of Hom complexes
\[
\xymatrix{
F \ar[r] \ar[d] & \End_{k[t^\pm]}(M'((t))) \ar[d] \\
\End_{k[t^\pm]}(M((t))) \ar[r] & \Hom_{k[t^\pm]}(M((t)),M'((t)))
}
\]
where the lower horizontal arrow and the right vertical arrow
are induced by the composition with $M((t))\to M'((t))$ respectively.
Since the model category of dg modules admits a complicial model
structure, the lower horizontal arrow is a trivial fibration, i.e.,
a surjective quasi-isomorphism. Thus,
the upper horizontal arrow is also a trivial fibration.
The right vertical arrow is
also a weak equivalence, i.e., a (injective) quasi-isomorphism
since $M((t))\to M'((t))$ is so.
By the 2-out-of-3 property, the left vertical arrow is
a (injective) weak equivalence.
Note that $F$ is a dg Lie subalgebra of $\End_{k[t^\pm]}(M((t)))$
consisting of those linear maps $M((t))\to M((t))$ such that
the composite $M((t))\to M((t))\to M'((t))$ factors through
the quotient $M((t))\to M'((t))$. Moreover, the upper horizontal
arrow is a morphism of dg Lie algebras.
Thus, we see the first assertion.
To prove the second claim, note first that
$f$ is a trivial fibration between $\mathbf{Lie}^\dagger$-algebras,
so that $L((t)):C[1] \to \End_{k[t^\pm]}(M((t)))$ factors through
$F\subset \End_{k[t^\pm]}(M((t)))$.
We then have the commutative diagram
\[
\xymatrix{
C[1] \ar[r]^(0.5){L((t))} \ar[d] & F \ar[d]  \\
C'[1] \ar[r]_(0.3){L'((t))} & \End_{k[t^\pm]}(M'((t))) 
}
\]
where the left vertical arrow is the morphism of dg Lie algebras
induced by $f$.
Thus, we have a zig-zag of homotopy equivalences of mapping spaces of dg Lie algebras
\[
\Map(C[1],E)\leftarrow \Map(C[1],F) \to \Map(C[1],E') \leftarrow \Map(C'[1],E')
\]
where $E=\End_{k[t^\pm]}(M((t)))$ and $E'=\End_{k[t^\pm]}(M'((t)))$.
Consequently, we obtain a homotopy equivalence
between the mapping space from $L((t))$ to $0$ and that of $L'((t))$ to $0$.
\QED

Here is another $\infty$-categorical proof of the above Lemma.

\begin{Lemma}
\label{alternative}
Let $(C,M,[-,-], l,\delta)$ and $(C',M',[-,-]', l',\delta')$ be two $\mathbf{Lie}^\dagger$-algebras in the
category of complexes of $k$-modules.
Let $f:(C,M,[-,-], l,\delta)\to(C',M',[-,-]', l',\delta')$
be a weak equivalence of $\mathbf{Lie}^\dagger$-algebras.
Then there exists a natural homotopy equivalence between the mapping space from $L((t)):C[1]\to \End_{k[t^\pm]}(M((t)))$ to the zero morphism (see Section~\ref{mappingsp})
and the mapping space from $L'((t)):C'[1]\to \End_{k[t^\pm]}(M'((t)))$ to the zero morphism
via the zig-zag (see the proof for the construction).
\end{Lemma}

\Proof
Let $U:\textup{dgLie}\rightleftarrows \textup{dgAlg}$
be the adjoint pair where the right adjoint functor
from the category of dg algebras
$\textup{dgAlg}\to \textup{dgLie}$ carries 
a dg algebra $A$ to the dg Lie algebra $A$ with the graded commutator.
The left adjoint is given by the free functor of the
universal envelopping algebras.
If $\textup{dgAlg}$ and $\textup{dgLie}$ are endowed with the projective
model structures respectively, the adjoint pair is a Quillen adjunction
since the right adjoint preserves trivial fibrations and fibrations.
The morphism of dg Lie algebras
$L((t)):C[1]\to \End_{k[t^\pm]}(M((t))$ and
$L'((t)):C'[1]\to \End_{k[t^\pm]}(M'((t))$
induce morphisms of dg algebras
$U(L((t))):U(C[1])\to \End_{k[t^\pm]}(M((t))$ and
$U(L'((t))):U(C'[1])\to \End_{k[t^\pm]}(M'((t))$
respectively.
If necessary,
we replace $C[1]$ and $C'[1]$ by the cofibrant models.
By the Quillen adjunction,
we may work with dg algebras instead of dg Lie algebras.
Namely, it is enough to prove that
there exists a homotopy equivalence between the mapping space from $U(L((t))):U(C[1])\to \End_{k[t^\pm]}(M((t)))$ to $U(0):U(C[1])\to \End_{k[t^\pm]}(M((t)))$
and the mapping space from $U(L'((t)))\circ U(f):U(C[1])\stackrel{U(f)}{\to} U(C'[1])\to \End_{k[t^\pm]}(M'((t)))$ to $U(0):U(C[1])\to \End_{k[t^\pm]}(M'((t)))$.
For simplify the notation, we put $F=U(L((t)))$, $F'=U(L'((t)))\circ U(f)$, and $A=U(C[1])$.
We introduce some categories.
Let $\Mod_A$ be the category of dg $k[t^\pm]\otimes A$-modules
and let $\Mod$ be the category of dg $k[t^\pm]$-modules.
We consider them to be (combinatorial)
model categories which are endowed with projective model structures.
(A morphism is a weak equivalence if the underlying map is a quasi-isomorphism,
and a morphism is a fibration if the underlying map is a termwise surjective
map. Every dg $k[t^\pm]$-module is cofibrant.).
Taking the hammock localization, fibrant replacement of simplicial
categories \cite[1.1.4]{HTT}, and the simplicial nerve functor
\cite[1.1.5]{HTT}
we associate to the model categories $\Mod_A$ and $\Mod$ $\infty$-categories
$\Map_A^\infty$ and $\Map^\infty$ respectively. We also obtain
the functor $\Mod_A^\infty\to \Mod^\infty$ from
the forgetful functor $\Mod_A\to \Mod$.
This functor can be decomposed as
$\Map_A^\infty\stackrel{u}{\to} \tilde{\Map}_A^\infty \stackrel{v}{\to} \Map^\infty$
where $u$ is a trivial cofibration, and $v$ is a fibration
with respect to Joyal model structure (see \cite[2.2.5]{HTT}).
Let $\Mod_{A}^{\infty}\times^h_{\Mod^\infty}\{M((t))\}$
denote the homotopy fiber over $M((t))\in \Mod^\infty$.
We may take $\tilde{\Mod}_{A}^{\infty}\times_{\Mod^\infty}\{M\}$ as
one of its model.
Let $\Map(A,\End_{k[t^\pm]}(M((t))))$
denote the mapping space from $A$ to $\End_{k[t^\pm]}(M((t)))$
which we regard an $\infty$-category that is an $\infty$-groupoid.
By the categorical characterization of the endomorphism algebras
\cite[6.1.2.41]{HA},
there is a natural equivalence
$\Map(A,\End_{k[t^\pm]}(M((t))))\simeq \Mod^\infty_A\times^h_{\Mod^\infty} \{M((t))\}=\tilde{\Mod}_{A}^{\infty}\times_{\Mod^\infty}\{M((t))\}$
which carries a morphism $\rho:k[t^\pm]\otimes A\to \End_{k[t^\pm]}(M((t)))$ 
to the dg $k[t^\pm]$-module $M((t))$ endowed with the $k[t^\pm]\otimes A$-module structure given
by $\rho$ up to equivalence.
We may consider that $M((t))$ belongs to
$\tilde{\Mod}_{A}^{\infty}\times_{\Mod^\infty}\{M((t))\}$.
Similarly, we have
$\Map(A,\End_{k[t^\pm]}(M'((t))))\simeq \Mod^\infty_A\times^h_{\Mod^\infty} \{M'((t))\}=\tilde{\Mod}_{A}^{\infty}\times_{\Mod^\infty}\{M'((t))\}$.
The quasi-isomorphism $M((t))\to M'((t))$
induces the functorial equivalence
\[
\Mod^\infty_A\times^h_{\Mod^\infty} \{M((t))\}\to \Mod^\infty_A\times^h_{\Mod^\infty} \{M'((t))\}.
\]
To understand it, we let $\Fun([0,1],\Mod)$ be the category of the functor category endowed with the projective model structure
(see e.g. \cite[A. 2.8.2]{HTT}).
We here denote by $[0,1]$ the category arising from the linearly ordered set
$0\to 1$.
Let $\Fun([0,1],\Mod)^\infty$
be the $\infty$-category associated to $\Fun([0,1],\Mod)$ as above.
There are two maps $d_0,d_1:\Fun([0,1],\Mod)^\infty\rightrightarrows \Mod^\infty$ given by the composition with $\{0\}\to [0,1]$ and $\{1\}\to [0,1]$
respectively, and the section $\Mod^{\infty}\to \Fun([0,1],\Mod)^\infty$
induced by the constant functors. These three maps are (categorical)
equivalences
(cf. \cite[4.2.4.4]{HTT}).
Then we consider
\[
H:=\tilde{\Mod}_{A}^{\infty}\times_{\Mod^\infty, d_0}\times \Fun(\Delta^1,\Mod)^\infty\times_{d_1,\Mod^\infty}\{M((t))\}.
\]
Replacing $M$ with $M'$ we define $H'$ in a similar way.
By the canonical
categorical equivalence $\Mod^{\infty}\to \Fun([0,1],\Mod)^\infty$,
the homotopy fiber $\Fun([0,1],\Mod)^\infty\stackrel{d_1}{\to} \Mod^\infty$ over $\{M((t))\}$ is a contractible space.
(There is a split injective
$\{M((t))\}\to\Fun(\Delta^1,\Mod)^\infty\times_{d_1,\Mod^\infty}\{M((t))\}$
which makes $\{M((t))\}$ the homotopy fiber.)
The composition with $M((t))\to M'((t))$ induces
\[
\Fun(\Delta^1,\Mod)\times_{d_1,\Mod}\{M((t))\}\to \Fun(\Delta^1,\Mod)\times_{d_1,\Mod}\{M'((t))\}
\]
which carries $N\to M((t))$ to $N\to M((t))\to M'((t))$.
It gives rise to a morphism $H\to H'$.
It is induced by the morphism from the following diagram $D:I\to \textup{Set}_\Delta$ of $\infty$-categories
\[
\xymatrix{
\tilde{\Mod}_A^\infty  \ar[dr]   &  &  \Fun([0,1],\Mod)^\infty \ar[dl]^{d_0} \ar[dr]_{d_1} &   &   \{M((t))\} \ar[ld] \\
 & \Mod^\infty &   &  \Mod^\infty &  
}
\]
to another diagram $D':I\to \textup{Set}_\Delta$ in which $M((t))$ is replaced with $M'((t))$.
Here $I$ is the index category, and $\textup{Set}_\Delta$
is the category of simplicial sets.
We use the injective model structure on the functor category
$\Fun(I,\textup{Set}_\Delta)$ (see e.g. \cite[A. 2.8.2, A. 2.8.7]{HTT})
and 
take functorial fibrant replacements $D\to \overline{D}$ and $D'\to \overline{D}'$ such that the diagram
\[
\xymatrix{
D \ar[d] \ar[r] & D' \ar[d] \\
\overline{D} \ar[r] & \overline{D}'
}
\]
commutes. We let $\overline{H}$ and $\overline{H}'$
be the limits of diagram $\overline{D}$ and $\overline{D}'$
respectively. These limits are homotopy limits so that
we take $\overline{H}$ and $\overline{H}'$ as homotopy fibers
$\Mod_A^\infty\times^h_{\Mod^\infty}\{M((t))\}$ and $\Mod_A^\infty\times^h_{\Mod^\infty}\{M'((t))\}$ respectively.
The map $\overline{D} \to \overline{D}'$ induces
an equivalence $\overline{H}\to \overline{H}'$.
Let $P=(M((t))_\rho,\textup{id}_{M((t))}: M((t))\to M((t)))$
be the pair where $M((t))_\rho$ endowed with the
$k[t^\pm]\otimes A$-module structure determined by $\rho$. It belongs to $H$ (or $\Mod_A^\infty\times_{\Mod^\infty}\{M((t))\}$). 
Then $H\to H'$ sends $P$ to $Q=(M((t))_\rho, M((t))\to M'((t)))$.
The pair $Q$ is equivalent to $R=(M'((t))_{\rho'}, \textup{id}_{M'((t))}:M'((t))\to M'((t)))$ in $H'$, where $M'((t))_{\rho'}$ is the $k[t^\pm]\otimes A$-module
$M'((t))$ determined by $F'$. The image of the last pair $R$ in $\overline{H}'$ 
corresponds to
$F':A\to \End_{k[t^\pm]}(M'((t)))$ up to equivalence
(while the image of $P$ in $\overline{H}$ corresponds to $F$).
Consequently, we have
$\Map(A,\End_{k[t^\pm]}(M((t))))\simeq \overline{H} \simeq \overline{H'}\simeq \Map(A,\End_{k[t^\pm]}(M'((t))))$ which carries $F$ to $F'$ up to equivalence. Hence our claim follows.
\QED

Let $\Alg_{\mathcal{C}alc}$ (resp. $\Alg_{\mathcal{L}ie^\dagger}$,
$\Alg_{\mathbf{calc}}$) be the category of $\mathcal{C}alc$-algebras
(resp.  $\mathcal{L}ie^\dagger$-algebras,
$\mathbf{calc}$-algebras)
in the category of chain complexes.
As in the case of $\mathbf{Lie}^\dagger$-algebras, by \cite{Hin}
$\Alg_{\mathcal{C}alc}$, $\Alg_{\mathcal{L}ie^\dagger}$,
and $\Alg_{\mathbf{calc}}$ admit combinatorial model category structures
where a morphism is a weak equivalence (resp. fibration) if it induces
quasi-isomorphisms (resp. termwise surjective maps) of underlying complexes.
The natural maps
$\mathcal{C}alc\to \mathbf{calc}$ and $\mathbf{Lie}^\dagger\to \mathbf{calc}$
of operads induce the pullback functors
\[
\Alg_{\mathbf{calc}} \to \Alg_{\mathcal{C}alc} \ \ \textup{and}\ \ \Alg_{\mathbf{calc}}\to \Alg_{\mathbf{Lie}^\dagger}
\]
which are right Quillen functors.

\begin{Lemma}
\label{transfer}
Let $A$ be a dg algebra.
Consider an action of $\mathcal{C}alc$ on $(C^\bullet(A),C_\bullet(A))$
that satisfies the property in Proposition~\ref{goodcalculus}.
We denote by $\mathfrak{A}$ this $\mathcal{C}alc$-algebra.
Then there exists a $\mathbf{calc}$-algebra $\mathfrak{B}:=(C,M,\cdot, [-,-],i,l,\delta)$
in the category of complexes such that its pullback along
$\mathcal{C}alc\to \mathbf{calc}$ is weak equivalence to $\mathfrak{A}$,
and its pullback $\iota^*\mathfrak{B}$ along $\iota:\mathbf{Lie}^\dagger\to \mathbf{calc}$
is weak equivalent to the $\mathbf{Lie}^\dagger$-algebra
$(C^\bullet(A),C_\bullet(A),[-,-],L,B)$.

In particular, 
there exist a cofibrant $\mathbf{Lie}^\dagger$-algebra $\mathfrak{I}$ and
a diagram of morphisms of $\mathbf{Lie}^\dagger$-algebras
\[
\iota^*\mathfrak{B}\longleftarrow \mathfrak{I} \longrightarrow \mathfrak{C}:=(C^\bullet(A),C_\bullet(A),[-,-],L,B)
\]
where both arrows are trivial fibrations.
\end{Lemma}

\Proof
We first note that by \cite[2.4.5]{Hin} $\Alg_{\mathbf{calc}} \to \Alg_{\mathcal{C}alc}$ induces an equivalence of homotopy categories.
It implies the first claim. Suppose that
$\mathfrak{B}$ is a $\mathbf{calc}$-algebra whose pullback to
$\Alg_{\mathcal{C}alc}$ is weak equivalent to $\mathfrak{A}$.
It follows that its pullback to $\Alg_{\mathcal{L}ie^\dagger}$
is weak equivalent to the pullback of $\mathfrak{C}$
to $\Alg_{\mathcal{L}ie^\dagger}$.
Notice that $\mathcal{L}ie^\dagger\to \mathbf{Lie}^{\dagger}$
is a weak equivalence and thus the pullback functor
$\Alg_{\mathbf{Lie}^\dagger}\to \Alg_{\mathcal{L}ie^\dagger}$
induces an equivalence of homotopy categories.
It follows that $\iota^*\mathfrak{B}$ is weak equivalent to
$\mathfrak{C}$.

The final statement follows from the model category structure of
$\Alg_{\mathbf{Lie}^\dagger}$.
Indeed, take a cofibrant replacement $r:\mathfrak{I}'\to \iota^*\mathfrak{B}$
that is a trivial fibration,
and choose a weak equivalence $\mathfrak{I}'\to \mathfrak{C}$.
The weak equivalence $\mathfrak{I}'\to \mathfrak{C}$ is decomposed into
$\mathfrak{I}'\stackrel{\alpha}{\to} \mathfrak{I}\stackrel{\beta}{\to} \mathfrak{C}$ where $\alpha$ is a trivial cofibration and $\beta$ is a trivial
fibration. Choose $s:\mathfrak{I}\to \iota^*\mathfrak{B}$ of $\alpha$
such that $r=s\circ \alpha$.
Then we have the desired diagram $\iota^*\mathfrak{B}\stackrel{s}{\leftarrow} \mathfrak{I}\stackrel{\beta}{\to} \mathfrak{C}$.
\QED

\subsection{}

\begin{Lemma}
\label{nicerel}
We adopt notation in Lemma~\ref{transfer}.
Suppose that $\mathfrak{B}=(C,M,\cdot, [-,-],i,l,\delta)$
is  a $\mathbf{calc}$-algebra.
Then

\begin{enumerate}
\renewcommand{\labelenumi}{(\theenumi)}

\item 
\[
[i_a,i_{[b,c]}]=0,
\]

\item
\[
[i_a,[i_b,l_c]]=0
\]

\end{enumerate}
for any $a,b,c\in C$.
\end{Lemma}

\Proof
\begin{eqnarray*}
[i_a,i_{[b,c]}] &=& i_ai_{[b,c]}-(-1)^{|i_a||i_{[b,c]}|}i_{[b,c]}i_a \\
&=& i_{a\cdot [b,c]}-(-1)^{|a|(|b|+|c|-1)}i_{[b,c]\cdot a}
\end{eqnarray*}
where we use the module structure
$i:C\otimes M\to M$ in the second equation.
Therefore it will suffices to prove that $a\cdot[b,c]-(-1)^{|a|(|b|+|c|-1)}[b,c]\cdot a=0$.  It is a direct consequence of the Koszul sign rule
and the graded commutativity.
Finally,  
$[i_a,[i_b,l_c]]=0$ follows from (1) and $[i_b,l_c]=i_{[b,c]}$.
\QED

We associate a morphism $l((t)):C[1]\to \End_{k[t^\pm]}(M((t)))$
of dg Lie algebras
to $\mathfrak{B}=(C,M,\cdot, [-,-],i,l,\delta)$ (see Section~\ref{extendlie}).
Let $i((t)):C\to \End_{k[t^\pm]}(M((t)))$ be a morphism of dg algebras
which carries $c$ to $i_c:M((t))\to M((t))$
defined by $i_c(mt^n)=i_c(m)t^n$.

\begin{Proposition}
\label{path}
We have
\[
e^{-\frac{1}{t}i((t))}\bullet 0=\sum_{n=0}^{\infty}\frac{[-\frac{1}{t}i((t)),-]^n}{(n+1)!}([-\frac{1}{t}i((t)),0]+d_{H}(\frac{1}{t}i((t))))=l((t))
\]
in the Hom complex $H:=\Hom(\overline{B_{com}(C[1])},\End_{k[t^\pm]}(M((t))))$. Here $d_H$
denotes the differential of the Hom complex.
\end{Proposition}

\Proof
For ease of notation, we put $I:=i((t))$ and $L:=l((t))=L$ in this proof.
Also, if $c\in C$, we write $I_c$ and $L_c$ for $i_c((t))$ and $l_c((t))$
of degree $-1$ and $0$ respectively.
We first calculate $[-I,0]+d_{H}(\frac{1}{t}I)=d_H(\frac{1}{t}I)$.
Let $\partial$ denote the differential of $M$.
Let $d_E(-)=[\partial+t\delta,-]$ be the differential of
$\End_{k[t^\pm]}(M((t)))$.
Let $d_B$ be the differential of $B_{com}(C[1])$.
Then we have
\begin{eqnarray*}
d_H(\frac{1}{t}I) &=& d_E(\frac{1}{t}I)-(-1)^{0\cdot |d_B|}\frac{1}{t}I\circ d_B \\
&=& \frac{1}{t}[\partial,I]+[t\delta,\frac{1}{t}I]-\frac{1}{t}I\circ d_B \\
&=& L+\frac{1}{t}([\partial,I]-I\circ d_B).
\end{eqnarray*}
Note that we here used $L=[\delta,I]$ and $\frac{1}{t}I$ has degree $0$ in $H$.
Notice also that $I$ killes all higher part $\oplus_{n>1}\Sym^n(C[2])$,
and remember that $d_B$ is the sum $d_1+d_2$:
$d_1$ comes from that of $C[2]$, $d_2$ is generated by $Q_2:\Sym^2(C[2])\stackrel{(s^{-1})^{\otimes 2}}{\simeq}C[1]\wedge C[1]\stackrel{-[-,-]}{\to} C[1] \stackrel{s}{\simeq} C[2]$. Here $s:C[1]\to C[1][1]$ is the obvious ``identity'' 
map. (See \cite{LM}, \cite{Sch} for the detailed account on sign issue.)
We deduce that $I\circ d_B=I_{d_C(-)}-I_{[-,-]}$.
Note that $i:C\to \End(M)$ is a dg map.
It follows that $[\partial,I]=I_{d_C(-)}$.
Hence
\[
d_H(\frac{1}{t}I)=L+\frac{1}{t}([\partial,I]-I_{d_C(-)}+I_{[-,-]})=L+\frac{1}{t}I_{[-,-]}.
\]
Therefore, $e^{-\frac{1}{t}I}\bullet 0$ equals to
\begin{eqnarray*}
L+\frac{1}{t}I_{[-,-]}-\frac{1}{2t}[I,L]_H+(\textup{terms of the forms}\ \frac{1}{t^n}[I[\cdots [I,I_{[-,-]}]\cdots]\ \textup{and}\ \frac{1}{t^n}[I[\cdots [I,L]\cdots])
%+(\textup{terms of the form} \frac{1}{t^n}[I[\cdots [I,L]\cdots]).  \\
\end{eqnarray*}
Here the bracket appearing in $I_{[-,-]}$ is that of $C[1]$,
and the bracket $[-,-]_H$ is the convolution Lie bracket of $H$.
According to Lemma~\ref{nicerel},
the terms of $n$-ary operations ($n\ge 3$) vanish.
Using the definition of convolution Lie bracket we have the formula
\[
[I,L]_{H}(a\otimes b)=[I_a,L_b]_E+(-1)^{|a||b|+|a|+|b|}[I_b,L_a]_E
\]
for $a\otimes b\in C[1]\wedge C[1]$ where $[-,-]_E$ is the bracket of $\End_{k[t^\pm]}(M((t)))$, and $|a|$, $|b|$ are degrees in $C$.
Combined with $[I_a,L_b]_E=I_{[a,b]}$ coming from the structure of calculus
we conclude that
$2I_{[-,-]}=[I,L]_H$. Hence $e^{-\frac{1}{t}I}\bullet 0=L$.
\QED

The dg Lie algebra $H=\Hom(\overline{B_{com}(C[1])},\End_{k[t^\pm]}(M((t))))$
is pro-nilpotent. Indeed,
put $F^nH=\{f\in H|\ \textup{restriction} \oplus_{i<n} \Sym^i(C[2])\to \End_{k[t^\pm]}(M((t)))\ \textup{is zero}\}$. Each $F^nH$ is a dg ideal and $[F^nH,F^mH]\subset F^{n+m}H$. We have a natural isomorphism
$H\simeq \varprojlim H/F^nH$ such that each quotient
$H/F^nH$ is nilpotent. Thus the exponential action of $H^0$ on
$\MC(H)$ makes sense. More generally, we may consider the
action $\exp(H^0[u,du])$ on $\MC(H[u,du])$.
The element $-\frac{u}{t}i((t))$ in $H^0[u,du]$ acts on $0$
in $\MC(H[u,du])$, and we obtain $e^{-\frac{u}{t}i((t))}\bullet 0$
in $\MC(H[u,du])$.
Let us consider the map
\[
\MC(\Hom(\overline{B_{com}(C[1])},E)[u,du])\to \MC(\Hom(\overline{B_{com}(C[1])},E[u,du])).
\]
induced by the canonical map $k[u,du]\otimes \Hom(B_{com}(C[1]),E)\to \Hom(B_{com}(C[1]),k[u,du]\otimes E)$ of dg Lie algebras. Here $E=\End_{k[t^\pm]}(M((t)))$.
Let $\Phi$ be the image of $e^{-\frac{u}{t}i((t))}\bullet 0$
in $\MC(\Hom(B_{com}(C[1]),k[u,du]\otimes E)$. Consider two face maps
\[
d_0,d_1:\MC(\Hom(\overline{B_{com}(C[1])},E[u,du]))\rightrightarrows \MC(\Hom(\overline{B_{com}(C[1])},E))
\]
where $d_0$ (resp. $d_1$) is determined by $u=du=0$ (resp. $u=1$ and $du=0$).
By Proposition~\ref{path} we see that $d_0(\Phi)=0$ and $d_1(\Phi)=l((t))$.
Consequently, we have

\begin{Proposition}
\label{goldenpath}
The element $\Phi$ (constructed from $-\frac{u}{t}i((t))$) gives
a morphism/homotopy from $0$ to $l((t))$.
\end{Proposition}

\subsection{}
\label{strictconst}
We continue to consider the
$\mathbf{calc}$-algebra $\mathfrak{B}=(C,M,\cdot, [-,-], i,l,\delta)$
(see Lemma~\ref{transfer}).
Let us consider the natural morphism of dg Lie algebras
\[
j:\End_{k[t]}(M[[t]]) \longrightarrow \End_{k[t^\pm]}(M((t)))
\]
induced by the base change $\otimes_{k[t]}k[t^\pm]$.
Consider a homotopy fiber of this morphism.
It is well-defined in the $\infty$-category or the model category of dg Lie algebras (or $L_\infty$-algebras).
Let $E_\ge :=\End_{k[t]}(M[[t]])$ and $E:=\End_{k[t^\pm]}(M((t)))$.
We use the product of dg lie algebras
\[
F_j:=\{(e,f(u,du))|f(0,0)=0,\ f(1,0)=j(e)\}\subset E_\ge \times E[u,du]
\]
together with the first projection to $E_\ge$ as a model of the homotopy fiber (the differential is given by $d(e,e')=(de,de')$).
The underlying complex of $F_j$ is quasi-isomorphic
to the standard mapping cocone $E_\ge\oplus E[-1]$ endowed with
the differential $d(e,e')=(de,j(e)-de')$ over $E$.
(Note that $F_j$ is a homtopy pullback exactly when the underling complex 
is a homotopy pullback at the level of complexes.)
An explicit quasi-isomorphism
$\iota:E_\ge\oplus E[-1]\to F_j$ is given by the formula
$(e,e')\mapsto (e,u\otimes e+du\otimes e')$. It has a homotopy inverse
$\pi:F_j\to E_\ge \oplus E[-1]$ defined by $(e, f(u)\otimes e'_1+g(u)du\otimes e_2')\mapsto (e,e'_2\int_0^1g(u)du)$ (cf. \cite[Section 3]{FM}).

We take the morphism/homotopy $\Phi$, constructed in Proposition~\ref{goldenpath}, from $0:C[1]\to \End_{k[t^\pm]}(M((t)))$ to $l((t)):C[1]\to \End_{k[t^\pm]}(M((t)))$. Actually, $\Phi$ is an $L_\infty$-morphism $\overline{B_{com}(C[1])}\to E[u,du]$.
It gives rise to an $L_\infty$-morphism $l[[t]]\times \Phi:B_{com}(C[1])\to F_j\subset  E_\ge \times E[u,du]$ such that the diagram
\[
\xymatrix{
  &  C[1]  \ar[ld]_{l[[t]]\times \Phi} \ar[d]^{l[[t]]} \\
F_j \ar[r]_{\textup{pr}_1} & E_{\ge}
}
\]
commutes
where arrows are implicitly considered to be $L_{\infty}$-morphisms.

\begin{Remark}
\label{contractionremark}
Observe that the morphism $l[[t]]\times \Phi:C[1]\to F_j$ at the level of
complexes without bracket
can be represented by
$\frac{1}{t}i((t)):C[1]\to (\End_{k[t^\pm]}(M((t)))/\End_{k[t]}(M[[t]]))[-1]$.
In fact, if we think of $C[1]\to E[u,du]$ as a map of complexes
by forgetting the non-linear terms,
and consider $\Hom(C[1],E[u,du])$ to be the dg Lie algebra with the
trivial bracket, the straightforward computation
of the gauge action of $-\frac{u}{t}i((t))$ on it
shows that the underlying map
$C[1]\to E[u,du]$ is given by $u\otimes l((t))+du\otimes \frac{1}{t}i((t))$.
The composite $C[1]\stackrel{l[[t]]\times \Phi}{\to} F_j\stackrel{\pi}{\to}E_\ge \oplus E[-1]$ is $l[[t]]\times \frac{1}{t}i((t))$ ($\pi$ is the quasi-isomorphism).
Since $E_\ge\to E$ is injective, there is a natural quasi-isomorphism
$E_\ge \oplus E[-1]\to (E/E_{\ge})[-1]$ given by the second projection.
Namely, $(E/E_{\ge})[-1]$ is another model of the mapping cocone.
It follows that $l[[t]]\times\Phi$ at the level
of complexes gives rise to
$C[1]\to (E/E_{\ge})[-1]$ which carries $P$ to $\frac{1}{t}i_P((t))$.
\end{Remark}

\subsection{}
\label{trivialityend}
We now apply the construction in Section~\ref{strictconst}
to the $\mathbf{Lie}^\dagger$-algebra
\[
\mathfrak{C}=(C^\bullet(A),C_\bullet(A),[-,-]_G,L,B)
\]
considered in Lemma~\ref{transfer}.
Consider the morphisms of dg Lie algebras
\[
L[[t]]:C^\bullet(A)[1]\to \End_{k[t]}(C_\bullet(A)[[t]]),\ \ \textup{and}\ \ \ L((t)):C^\bullet(A)[1]\to \End_{k[t^\pm]}(C_\bullet(A)((t)))
\]
(see Section~\ref{cylicify}).

\begin{Lemma}
\label{nice}
There is a homotopy $\Psi$ defined in $\MC(C^\bullet(A)[1],\End_{k[t^\pm]}(C_\bullet(A)((t)))[u,du])$,
from $0$ to $L((t))$ that comes from $\Phi$.
\end{Lemma}

\Proof
We apply 
Lemma~\ref{mapfunct} and Lemma~\ref{transfer}
(or one can simply use Lemma~\ref{alternative}).
Namely, by Lemma~\ref{transfer}
there is a diagram $\iota^*\mathfrak{B}\leftarrow \mathfrak{I} \rightarrow \mathfrak{C}$ of trivial fibrations of $\mathbf{Lie}^\dagger$-algebras.
Then by Lemma~\ref{mapfunct}, using a zig-zag of homotopy equivalences
(see the proof of Lemma~\ref{mapfunct})
we can transfer $\Phi$ in Proposition~\ref{goldenpath}
to a 1-simplex $\Psi$ of the Kan complex
\[
\MC(C^\bullet(A)[1],\Omega_\bullet\otimes \End_{k[t^\pm]}(C_\bullet(A)((t)))).
\]
\QED

{\it Proof of Proposition~\ref{triviality}.}
Proposition~\ref{triviality} follows from Lemma~\ref{nice}.
\QED

\subsection{}
\label{const}
For simplicity, we put $\EE_\ge=\End_{k[t]}(C_\bullet(A)[[t]])$,
$\EE=\End_{k[t^\pm]}(C_\bullet(A)((t)))$ and $\EE[u,du]=k[u,du]\otimes \EE$.
Let $\iota:\EE_\ge\to \EE$ be the natural injective morphism of dg Lie algebras
induced  by the base change $\otimes_{k[t]}k[t^\pm]$.
Let $\Gr:=F_{\iota}\subset \EE_{\ge}\times \EE[u,du]$
be a (model of) homotopy fiber of $\iota$ defined in a similar way as
$F_j$ in Section~\ref{strictconst}:
\[
\Gr=\{(p,q(u,du))|\  d_0(q(u,du))=0, d_1(q(u,du))=\iota(p),\ \textup{i.e.,}\ q(0,0)=0,\ q(1,0)=\iota(p),\}.
\]
Note that $\Gr$ depends on $C_\bullet(A)[[t]]$, and 
we here
abuse notation by omitting the subscript that indicates $C_\bullet(A)[[t]]$.
As in Section~\ref{strictconst}, we define an $L_\infty$-morphism
\[
L[[t]]\times \Psi:B_{com}(C^\bullet(A) [1])\to \Gr\subset \EE_{\ge}\times \EE[u,du]
\]
where $\Psi$ is the $L_\infty$-morphism in Lemma~\ref{nice}.

For ease of notation we write
\[
\mathcal{P}:=L[[t]]\times \Psi:C^\bullet (A)[1]\longrightarrow \Gr
\]
for this $L_\infty$-morphism.
In summary, we have the following commutative diagram:
\[
\xymatrix{
C^\bullet(A)[1]  \ar[ddr]_{\Psi} \ar[dr]^(0.8){\mathcal{P}} \ar[drr]^{L[[t]]} &    &  \\
 & \Gr \ar[r]_{\textup{pr}_1} \ar[d]^{\textup{pr}_2} &  \EE_\ge \ar[d]^{\iota} \\
 & \EE[u,du] \ar[r]_{d_1} &  \EE.
}
\]
The $L_\infty$-morphism $\mathcal{P}:C^\bullet (A)[1]\longrightarrow \Gr$
plays the main role in this paper.
In the next Section, it turns out
that this $L_\infty$-morphism ``amounts to'' a period mapping for
infinitesimal (curved) deformations of $A$ via a moduli-theoretic interpretation. 
Note that there is another map $d_0:\EE[u,du]\to \EE$
which is a part of data of the path object $\EE\hookrightarrow \EE[u,du]\stackrel{(d_0,d_1)}{\to} \EE\times \EE$, and 
the composite $C^\bullet(A)[1]\to \EE[u,du]\stackrel{d_0}{\to} \EE$
is the zero morphism.

\subsection{}
\label{peri}
The dg Lie algebra $\Gr$ is a homotopy fiber of
$\iota:\EE_\ge\to \EE$, and $\EE_\ge$ and $\EE$ ``represent''
the deformations of negative and periodic cyclic complexes respectively
(cf. Section~\ref{premoduli}).
The purpose of Section~\ref{peri} is 
to explain a modular interpretation of $\Def_{\mathbb{F}}$
associated to the homotopy fiber $\mathbb{F}$.
It turns out that $\Def_{\mathbb{F}}$ is a generalization
of formal Sato Grassmannian to the level of complexes.

\begin{Definition}
\label{PTD}
Let $R$ be an artin local $k$-algebra with residue field $R/m_R\simeq k$,
that is, $R$ belongs to $\Art$.
Let $Z$ be a dg $k[t]$-module.
%Suppose that $Z$ is free over $k[t]$ if one forgets the differential.
Let $\tilde{Z}$
be a deformation of $Z$ to $R$, that is,
a dg $R[t]$-module $\tilde{Z}$ such that
the underlying graded $R[t]$-module of $\tilde{Z}$ is
$Z\otimes_{k[t]}R[t]$, and its reduction $\tilde{Z}\otimes_{R[t]}R/m_R[t]$ is
the dg $k[t]$-module $Z$ (namely,
the differential on $\tilde{Z}$ induces that of $Z$ via the canonical
identification $\tilde{Z}\otimes_{R[t]}R/m_R[t] \simeq Z$ of graded complexes, see Section~\ref{premoduli2}).

A periodically trivialized deformation of $Z$ to $R$
is a pair
\[
(\tilde{Z}, \phi: Z\otimes_{k[t]} R[t^\pm]^{tr} \stackrel{\sim}{\to} \tilde{Z}\otimes_{R[t]}R[t^\pm])
\]
where $\tilde{Z}$ is a deformation of $Z$,
and $Z\otimes_{k[t]} R[t^\pm]^{tr}$ denotes the trivial deformation of the dg $k[t^\pm]$-module $Z\otimes_{k[t]}k[t^\pm]$,
and
$\phi$ is an isomorphism of deformations of dg $k[t^\pm]$-modules to $R$.

Suppose that we are given two periodically trivialized deformations $(Z_1,\phi_1)$ and $(Z_2,\phi_2)$.
An isomorphism $(Z_1,\phi_1)\to (Z_2,\phi_2)$
is an isomorphism $h:Z_1\to Z_2$ of deformations
such that there exists some $a$ in $(\End_{k[t^\pm]}(Z\otimes_{k[t]}k[t^\pm])\otimes m_R)^{-1}$
such that the diagram
\[
\xymatrix{
Z\otimes_{k[t]} R[t^\pm]^{tr} \ar[r]^{e^{da}} \ar[d]_{\phi_1} & Z\otimes_{k[t]} R[t^\pm]^{tr}   \ar[d]^{\phi_2} \\
Z_1\otimes_{R[t]}R[t^\pm]    \ar[r]_{h\otimes_{R[t]}R[t^\pm]} & Z_2\otimes_{R[t]}R[t^\pm]
}
\]
commutes
where $d$ is the differential of $\End_{k[t^\pm]}(Z\otimes_{k[t]}k[t^\pm])\otimes m_R$.

\end{Definition}

\begin{Remark}
There is a more conceptual description
of the commutativity of the square in Definition~\ref{PTD}
(however, we will not need it in this paper).
Let $V$ be a deformation of the dg $k[t^\pm]$-module
$Z\otimes_{k[t]}k[t^\pm]$ to $R\in \Art$.
We denote by $Z\otimes_{k[t]} R[t^\pm]^{tr}$ the trivial deformation
as above.
Let $\psi_0:Z\otimes_{k[t]} R[t^\pm]^{tr}\to V$
and $\psi_1:Z\otimes_{k[t]} R[t^\pm]^{tr}\to V$
be isomorphisms of deformations.
A homotopy from $\psi_0$ to $\psi_1$ is defined to be
a $k[u,du]\otimes R[t^\pm]$-linear isomorphism  $h:k[u,du]\otimes_k(Z\otimes_{k[t]} R[t^\pm]^{tr})\to k[u,du]\otimes_kV$ such that (i) its reduction
$k[u,du]\otimes Z\otimes_{k[t]}k[t^\pm]\to k[u,du]\otimes Z\otimes_{k[t]}k[t^\pm]$ is the identity, (ii) the reduction by $u=du=0$ is $\psi_0$, and (iii)
the reduction by $u=1,\ du=0$, is $\psi_1$.
Put $\psi_0=e^{f}$ and $\psi_1=e^{g}$ where $f,g:Z\otimes_{k[t]} k[t^\pm]\to Z\otimes_{k[t]} k[t^\pm]\otimes_k m_R$ (use the logarithms).
Also, the isomorphism $h$ is of the form
$e^{n(u,du)}$ where $n(u,du):k[u,du]\otimes_k Z\otimes_{k[t]} k[t^\pm]\to k[u,du]\otimes_k Z\otimes_{k[t]} k[t^\pm]\otimes_k m_R$.
Then there is a homotopy from $\psi_0$ to $\psi_1$
if and only if there is an element
$a$ of degree $-1$ in $\End_{k[t^\pm]}(Z\otimes_{k[t]}k[t^\pm])\otimes m_R$
such that $\psi_1=\psi_0\circ e^{da}$.
To observe the ``only if'' direction, consider the composite
$k[u,du]\otimes e^{-f}\circ h:k[u,du]\otimes_k(Z\otimes_{k[t]} R[t^\pm]^{tr})\to k[u,du]\otimes_k(Z\otimes_{k[t]} R[t^\pm]^{tr})$.
It gives rise to the equation $e^{(-k[u,du]\otimes f) \cdot n(u,du)}\bullet 0=0$
in $\MC(\End_{k[u,du]((t))}(k[u,du]\otimes Z\otimes_{k[t]}k[t^\pm])\otimes m_R)$.
To simplify the notation, let $l^0(u)+du\otimes l^{-1}(u)=l(u,du)=(-k[u,du]\otimes f) \cdot n(u,du)$
(the multiplication ``$\cdot$'' is given by the Baker-Campbell-Hausdorff product). Then $e^{l(u,du)}\bullet 0=0$ corresponds to $dl(u,du)=0$
where $d$ is the differential.
 It is equivalent to simultaneous equations
\[
dl^0(u)=0\ \ \textup{and}\ \ l^0(u)'=d l^{-1}(u)
\]
where $l^0(u)'$ denotes the formal derivative of the polynomial.
Note also that $l(0,0)=(-f)\cdot f=0$.
It follows that there is an element
$a$ of degree $-1$ in $\End_{k[t^\pm]}(Z\otimes_{k[t]}k[t^\pm])\otimes m_R$
such that $l(1,0)=da$. Thus, $e^{-f}\circ e^{g}=e^{da}$.
To see the ``if'' direction, suppose that
we have $\psi_1=\psi_0\circ e^{da}$.
If we put $n(u,du)=e^{f}\circ e^{u\otimes da+du\otimes a}$, then
$n(u,du)$ gives a homotopy. Consequently, we see that
the commutativity of the square in Definition~\ref{PTD} amounts to
the homotopy between
$h\otimes_{R[t]}R[t^\pm]\circ \phi_1$ and $\phi_2$.

Finally, let us consider homotopies from the viewpoint of
the space of morphisms of deformations.
Let $\Hom_{\Omega_n\otimes k[t^\pm]}(\Omega_n\otimes Z\otimes_{k[t]} k[t^\pm]),\Omega_n\otimes Z\otimes_{k[t]} k[t^\pm])$ be the hom set of morphisms of
dg $\Omega_n\otimes k[t^\pm]$-modules.
These sets form a simplicial set $\Hom_{\Omega_\bullet\otimes k[t^\pm]}(\Omega_\bullet\otimes Z\otimes_{k[t]} k[t^\pm],\Omega_\bullet\otimes Z\otimes_{k[t]} k[t^\pm])$, which is a Kan complex of the mapping space
because $Z\otimes_{k[t]} k[t^\pm]$ is cofibrant with respect to the projective model structure.
Similarly, we have a Kan complex
$\Hom_{\Omega_\bullet\otimes R[t^\pm]}(\Omega_\bullet\otimes Z\otimes_{k[t]} R[t^\pm]^{tr},\Omega_\bullet\otimes V)$, that is a model of the mapping space.
The reduction by $R[t^\pm]\to k[t^\pm]$ induces
\[
\Hom_{\Omega_\bullet\otimes R[t^\pm]}(\Omega_\bullet\otimes Z\otimes_{k[t]} R[t^\pm]^{tr},\Omega_\bullet\otimes V)\to \Hom_{\Omega_\bullet\otimes k[t^\pm]}(\Omega_\bullet\otimes Z\otimes_{k[t]} k[t^\pm],\Omega_\bullet\otimes Z\otimes_{k[t]} k[t^\pm]).
\]
It is a Kan fibration since $V\to Z\otimes_{k[t]} k[t^\pm]$ is a surjective morphism of dg $R[t^\pm]$-modules (i.e., a fibration),
and the (projective) model category of dg $R[t^\pm]$-modules
is simplicial in the sense of \cite[1.4.2]{Hin2}.
The fiber $F$ of this Kan fibration over the constant simplicial
set determined by the identity map of $Z\otimes_{k[t]} k[t^\pm]$
is a homotopy fiber which is a Kan complex.
This Kan complex should be regarded as the $\infty$-groupoid
of morphisms of deformations from $Z\otimes_{k[t]} R[t^\pm]^{tr}$ to $V$.
An edge in this homotopy fiber $F$ is a homotopy
$h:k[u,du]\otimes_k(Z\otimes_{k[t]} R[t^\pm]^{tr})\to k[u,du]\otimes_kV$
defined above.
\end{Remark}

We apply Definition~\ref{PTD} to the dg $k[t]$-module $C_\bullet(A)[[t]]$
of negative cyclic complex.
We shall denote by $(C_\bullet(A)\otimes R)((t))^{tr}$
the trivial deformation $((C_\bullet(A)\otimes R)((t)),(\partial_{\textup{Hoch}}+tB)\otimes R)$ of the dg $k[t^\pm]$-module $C_\bullet(A)((t))$ (see Section~\ref{premoduli2}).
Let
\[
(Q=((C_\bullet(A)\otimes R)[[t]],\ \tilde{\partial}),\phi:(C_\bullet(A)\otimes R)((t))^{tr} \stackrel{\sim}{\to} Q\otimes_{R[t]}R[t^\pm] )
\]
be a pair
where $Q$ is a deformation of the dg $k[t]$-module
$C_\bullet(A)[[t]]$ to $R$, and $\phi$ is an isomorphism of deformations
of the dg $k[t^\pm]$-module $C_\bullet(A)((t))$.
Let $\dGr(R)$ be the set of isomorphism classes of
periodically trivialized deformations of $C_\bullet(A)[[t]]$ to $R$.
The assignment $R\mapsto \dGr(R)$
gives rise to a functor 
\[
\dGr:\Art\to \textup{Sets}.
\]
Again, $\dGr$ depends on $C_\bullet(A)[[t]]$
though we omit the subscript indicating $C_\bullet(A)[[t]]$.

The following is the moduli-theoretic presentation
of $\Def_{\mathbb{F}}$ in explicit terms:

\begin{Proposition}
\label{FMrevisited}
There is a natural equivalence $\Def_{\Gr}\to \dGr$ of functors.
This equivalence carries an element $c\in \Def_{\Gr}(R)$ represented by a Maurer-Cartan element $(\alpha,\beta)$
in $\MC(\Gr\otimes m_R)\subset \MC(\EE_\ge\otimes m_R)\times \MC(\EE[u,du]\otimes m_R)$
to the class of a periodically trivialized deformation of the form
\[
(Q_\alpha, (C_\bullet(A)\otimes R)((t))^{tr} \stackrel{\sim}{\longrightarrow} Q_\alpha\otimes_{R[t]}R[t^\pm]).
\]
Here $Q_\alpha$ denotes the
element in $\mathsf{D}_{C_\bullet(A)[[t]]}(R)$ which corresponds to $\alpha$
(see Proposition~\ref{compdeform2}).
\end{Proposition}

{\it Proof of Proposition~\ref{FMrevisited}.}
By \cite{FM},
there exist an $L_\infty$-structure on the mapping
cocone $\CC:=\mathbb{E}_\ge\oplus \mathbb{E}[-1]$ of $\iota:\mathbb{E}_\ge\to \mathbb{E}$ and an $L_\infty$-quasi-isomorphism
$\CC\to \mathbb{F}$.
(In this paper, an $L_\infty$-structure on a graded vector space $V$
is defined to be a unital graded cocommutative coalgebra $B_{com}V$
endowed with a square-zero coderivation $b:\overline{B_{com}V}\to \overline{B_{com}V}$, see Section~\ref{mappingsp} for the notation.
The coderivation $b:\overline{B_{com}V}\to \overline{B_{com}V}$
amounts to $\{b_n:\Sym^n(V[1])\to V[1]\}_{n\ge 1}$ satisfying
certain identities.)
Moreover, according to
\cite[Theorem 2]{FM} the set of Maurer-Cartan elements is given by
\[
\MC(\CC\otimes m_R)=\{(\alpha,a)\in \MC(\mathbb{E}_\ge\otimes m_R)\times \mathbb{E}^0\otimes m_R |\ e^{-a}\bullet 0=\iota(\alpha)\}.
\]
A Maurer-Cartan element $(\alpha,a)$ is gauge equivalent to another element
$(\beta,b)$
if and only if there exist $m\in \mathbb{E}^0\otimes m_R$
and $q \in \mathbb{E}^{-1}\otimes m_R$ such that
$e^m\bullet \alpha=\beta$ and $b={dq}\cdot a \cdot(-\iota(m))$ where in the last formula
we use Baker-Campbell-Hausdorff product.
Notice that there is a natural equivalence between
$\Def_{\CC}$ and $\dGr$.
In fact, there is a natural bijective map
$\Def_{\CC}(R)\to \dGr(R)$ which carries the class represented by
$(\alpha,a)$ to the class represented by
\[
(Q_\alpha, e^{-a}:(C_\bullet(A)\otimes R)((t))^{tr} \stackrel{\sim}{\longrightarrow} Q_\alpha \otimes_{R[t]}R[t^\pm]).
\]
By the invariance of $\Def$ with respect to $L_\infty$-quasi-isomorphisms,
the $L_\infty$-quasi-isomorphism $\CC\to \mathbb{F}$
induces an equivalence $\Def_{\CC}\to \Def_{\mathbb{F}}$ of functors.
In addition, the composition $\Def_{\CC}\to \Def_{\mathbb{F}}\to \Def_{\mathbb{E}_{\ge}}$
sends the class represented by $(\alpha,a)$ to the class of $\alpha$
(cf. \cite[Remark 5.3]{FM}).
Choose an $L_\infty$-quasi-isomorphism
$\mathbb{F}\to \mathbb{C}$ which is a homotopy inverse of
$\CC\to \mathbb{F}$.
Then $\mathbb{F}\to \mathbb{C}$ induces the inverse 
$\Def_{\mathbb{F}}\to \Def_{\mathbb{C}}$ of $\Def_{\CC}\to \Def_{\mathbb{F}}$.
We have the equivalence $\Def_{\mathbb{F}}\to \Def_{\mathbb{C}}\simeq \dGr$,
as desired.
\QED

We should think of the functor $\Def_{\Gr}\simeq \dGr$ as a formal neighborhood of a point on a
generalized Sato Grassmannian.
To understand it, we begin by reviewing the Sato Grassmannian.
Let $V$ be a finite dimensional $k$-vector space.
Let $R$ be a commutative $k$-algebra.
Consider a pair
$(W,f:(V\otimes R)((u))\simeq W\otimes_{R[[u]]}R((u)))$
such that $W$ is a finitely generated projective $R[[u]]$-module,
and $f$ is an isomorphism of $R((u))$-modules.
Here $R[[u]]$ is
the ring of formal power series, and
$R((u))$ is the ring of formal Laurent series.
An isomorphism $(W,f)\to (W',f')$ of such pairs
is defined to be an isomorphism $W\to W$ of $R[[u]]$-modules
that commutes with $f$ and $f'$.
Let $\mathsf{SGr}(R)$ be the set of isomorphism classes of pairs.
If $\CAlg_k$ denotes the category of commutative $k$-algebras,
the Sato Grassmannian as a functor
is the functor
$\CAlg_k \to \Set$ which assigns to each $R$ the set $\mathsf{SGr}(R)$.
Moreover, it is well-known that
this functor $\mathsf{SGr}$ is represented by th colimit
of the sequence of closed immersions of schemes $\varinjlim_{i\in\NN }X_i$.
Fix a pair $(W,f:V((u))\simeq W\otimes_{k[[u]]}k((u)))$ in $\mathsf{SGr}(k)$.
Let $R$ be in $\Art$.
A deformation of the pair $(W,f)$ is a pair $(\tilde{W}, \phi:(V\otimes R)((u))\simeq \tilde{W}\otimes_{R[[u]]}R((u)))$
such that $\tilde{W}$ is a deformation of $W$, that is,
a finitely generated projective
$R[[u]]$-module endowed with $\tilde{W}\otimes_{R[[u]]}k[[u]]\simeq W$, and $\phi:(V\otimes R)((u))\simeq \tilde{W}\otimes_{R[[u]]}R((u))$ is an isomorphism of $R((u))$-modules whose reduction to $k((u))$ is
$f$.
An isomorphism of deformations is defined in the obvious way.
Then consider the functor $\widehat{\mathsf{SGr}}_{(W,f)}:\Art \to \Set$ 
which carries $R$ to the set of isomorphism classes of deformations of
$(W,f)$. Informally, this functor $\Art \to \Set$
can be viewed as a formal neighborhood (i.e., a formal scheme)
at a point $(W,f)$ on $\mathsf{SGr}$ (the proper interpretation is left
to the interested reader).
Now we generalize it to the level of complexes.
We replace $k[[u]]$ by the dg algebra $k[t]$
where the cohomoligical degree of $t$ is $2$.
Note that $k[t]$ is $t$-adically complete in the sense that
$k[t]$ is a homotopy limit of the sequence
$\cdots \to k[t]/(t^{n+1}) \to  \cdots \to k[t]/(t^2) \to k[t]/(t)$.
Consider the $k[t]$-module $C_\bullet(A)[[t]]$, and the natural
isomorphism $C_\bullet(A)((t)) \stackrel{\sim}{\longrightarrow} C_\bullet(A)[[t]]\otimes_{k[t]}k[t^\pm]$ of $k[t^\pm]$-modules
instead of $(W,f)$.
It is natural to think of $\dGr\simeq \Def_{\mathbb{F}}$ as
a ``complicial generalization'' of $\widehat{\mathsf{SGr}}_{(W,f)}$,
and regard the dg Lie algebra $\mathbb{F}$ as the Lie-theoretic
presentation ($\mathbb{F}$ also has the ``derived structure'').
Thus, we call $\dGr\simeq \Def_{\mathbb{F}}$ the formal complicial
Sato Grassmannian or the complicial Sato Grassmannian simply.

\begin{Remark}
It is natural to ask for a global complicial Sato Grassmannian.
%Let $\CAlg_k^{cc}$ be the $\infty$-category of commutative dg algebras $A$
%over $k$
%such that $A^i=0$ for $i>0$.
It might be realized as a geometric object in derived algebraic
geometry developed by To\"en-Vezzosi and Lurie,
whose $R$-valued points informally parametrize the space
of pairs
$(W, (C_\bullet(A)\otimes_k R)((t))\simeq W\otimes_{R[t]}R[t^\pm])$ consisting
of a compact $R[t]$-module $W$ and an equivalence
of $R[t^\pm]$-modules.
\end{Remark}

\subsection{}
We now define a morphism $\mathsf{DAlg}_{A}\to \dGr$
of functors. The construction is based on the moduli-theoretic interpretation
of dg Lie algebras and their morphisms.
Let $\tilde{A}$ be a deformation of the dg algebra $A$ to $R$.
It gives rise to the negative cyclic complex
$C_\bullet(\tilde{A})[[t]]$ that is a deformation of 
the dg $k[t]$-module $C_\bullet(A)[[t]]$.
Namely, we associate to $\tilde{A}$ a deformation
$C_\bullet(\tilde{A})[[t]]$ of the dg $k[t^\pm]$-module $C_\bullet(A)[[t]]$
to $R$.
Let $\alpha$ be the Maurer-Cartan element in $\MC(C^\bullet(A)[1]\otimes m_R)$
that corresponds to $\tilde{A}$ (see Claim~\ref{deformclaim1}).
By Proposition~\ref{perioddeform1}, the morphism of dg Lie algebras $L[[t]]:C^\bullet(A)[1]\to \End_{k[t]}(C_\bullet(A)[[t]])$ sends $\alpha$ to the class represented by $L[[t]](\alpha)$ in $\Def_{\End_{k[t]}(C_\bullet(A)[[t]])}(R)$
that corresponds to the isomorphism class
of the deformation $C_\bullet(\tilde{A})[[t]]$ in $\mathsf{D}_{C_\bullet(A)[[t]]}(R)$
(see also Proposition~\ref{compdeform2}).
Since $L[[t]]$ has the lift $\mathcal{P}:C^\bullet(A)[1] \to \mathbb{F}$,
the element
$\Def_{\mathcal{P}}(\alpha)$ lying in $\dGr(R)$
through the isomorphism $\Def_{\mathbb{F}}(R)\simeq \dGr(R)$
in Proposition~\ref{FMrevisited}
is represented by a periodically trivialized deformation of the form
\[
(C_\bullet(\tilde{A})[[t]],\ (C_\bullet(A)\otimes R)((t))^{tr}\stackrel{\sim}{\to}C_\bullet(\tilde{A})((t))).
\]
We refer to it as (the isomorphism class of) the periodically trivialized deformation associated to $\tilde{A}$.
We obtain a morphism
\[
\mathsf{P}(R):\mathsf{DAlg}_{A}(R) \to \dGr(R),\ \ \tilde{A}\mapsto (C_\bullet(\tilde{A})[[t]],\ (C_\bullet(A)\otimes R)((t))^{tr}\stackrel{\sim}{\to}C_\bullet(\tilde{A})((t))).
\]
We refer to this morphism of functors as the period mapping for $A$.

Unwinding our construction based on results about modular interpretations
of morphisms of dg Lie algebras we can conclude:

\begin{Theorem}
\label{moduli}
Through the identifications $\mathsf{DAlg}_{A}\simeq \Def_{C^\bullet(A)[1]}$ and $\dGr\simeq \Def_{\Gr}$,
the morphism $\mathsf{P}:\mathsf{DAlg}_{A}\to \dGr$ can be identified with
$\Def_{\mathcal{P}}:\Def_{C^\bullet(A)[1]}\to \Def_{\Gr}$.
\end{Theorem}

We shall call both morphisms $\mathsf{P}:\mathsf{DAlg}_{A}\to \dGr$ and
$\Def_{\mathcal{P}}:\Def_{C^\bullet(A)[1]}\to \Def_{\Gr}$
the period mapping for $A$.

\begin{Remark}
In Introduction, for simplicity, we do not distinguish
$\mathsf{P}:\mathsf{DAlg}_{A}\to \dGr$ from 
$\Def_{\mathcal{P}}:\Def_{C^\bullet(A)[1]}\to \Def_{\Gr}$.
\end{Remark}

\section{Hodge-to-de Rham spectral sequence and complicial Sato Grassmannians}

In the last section, we constructed a period mapping from
the dg Lie algebra representing deformations of an algebra
to the dg Lie algebra of the complicial Sato Grassmanian.
In good cases, it turns
out that the complicial Sato Grassmanian admits a simple and nice structure.
Our main interest lies in the situation where
a non-commutative analogue of Hodge-to-de Rham spectral sequence
for cyclic homology theories degenerates.
Hence we start with a spectral sequence for cyclic homology theories.

\subsection{}
\label{smoothproper}
Let $A$ be a dg algebra over a field $k$ of characteristic zero.
We write $C_\bullet:=C_\bullet(A)$, $C_\bullet[[t]]:=(C_\bullet(A)[[t]],\partial_{\textup{Hoch}}+tB)$ and $C_\bullet((t)):=(C_\bullet(A)((t)),\partial_{\textup{Hoch}}+tB)$ for the Hochschild chain complex, the negative cyclic complex and
the periodic cyclic complex respectively (see Section~\ref{cyclic}).
(Note that the degree of $t$ is two.)
We let $C_\bullet[t^{-1}]:=(C_\bullet(A)[t^{-1}],\partial_{\textup{Hoch}}+tB)$
be the cyclic complex of $A$.
There is an exact sequence of complexes
\[
0\to C_\bullet[[t]]\to C_\bullet((t)) \to C_{\bullet}[t^{-1}]\cdot t^{-1}\to 0,
\]
whose associated long exact sequence
\[
\cdots \to HN_n(A)\to HP_n(A)\to HC_{n-2}(A)\to \cdots
\]
relates the negative and periodic
cyclic homology with the cyclic homology $HC_*(A)=H_*(C_{\bullet}[t^{-1}])$.
We define a decreasing filtration of $C_\bullet((t))$
by the formula
\[
F^i C_\bullet((t))_n=\prod_{r\ge i}  C_{n+2r}(A) \cdot t^r\subset C_\bullet((t))_n=\prod_{r\in \ZZ}  C_{n+2r}(A) \cdot t^r
\]
where $C_l(A)$ is the homologically $l$-th term of $C_\bullet(A)$.
This filtration forms a family of subcomplexes of $C_\bullet((t))$.
Note that $F^iC_\bullet((t))/F^{i+1}C_\bullet((t))$ is isomorphic to
$C_\bullet\cdot t^i$.
This filtration gives rise to a spectral sequence which we denote by
\[
HH_*(A)((t)) \Rightarrow HP_*(A).
\]
Each term in the $E_1$-stage is of the form $HH_j(A)\cdot t^i$.
The restriction $F^iC_\bullet[[t]]=C_\bullet[[t]]\cap F^i C_\bullet((t))$
also yields a spectral sequence
$HH_*(A)[[t]] \Rightarrow HN_*(A)$.
Since $C_\bullet[t^{-1}]\simeq C_\bullet((t))/t\cdot C_\bullet [[t]]$,
the filtration $F^i C_\bullet((t))/t\cdot C_\bullet [[t]]$ gives rise to
a spectral sequence $HH_*(A)[t^{-1}]\Rightarrow HC_{*}(A)=H_*(C_\bullet[t^{-1}])$.
We call these spectral sequences the Hodge-to-de Rham spectral sequences.
We refer the reader to \cite{Kal}
for the relation to the classical Hodge-to-de Rham
spectral sequence (see also Remark~\ref{classHodge}).
We now recall some properties of dg algebras.

\begin{Definition}
\label{propersmooth}
Let $A$ be a dg algebra over $k$
\begin{enumerate}
\renewcommand{\labelenumi}{(\theenumi)}

\item $A$ is proper if the underlying complex has finite dimensional cohomology,
and $H^n(A)=0$ for $|n|>>0$.

\item $A$ is smooth if $A$ is a compact object in the triangulated category
of dg $A\otimes A^{op}$-modules.

\end{enumerate}

\end{Definition}

\begin{Remark}
Put another way, the condition (2) in Definition~\ref{propersmooth}
is equivalent to the condition that $A$ belongs to
the smallest triangulated subcategory
which includes free $A\otimes A^{op}$-modules of finite rank
and is closed under direct summands.

\end{Remark}

\begin{Remark}
If $A$ is smooth and proper, it is easy to see that $HH_n(A)$ is finite dimensional for $n\in \ZZ$, and 
$HH_n(A)=0$ for $|n|>>0$.
(Keep in mind that $HH_n(A)$ is
$\textup{Tor}^{A\otimes A^{op}}_n(A,A)$.)
\end{Remark}

\subsection{}
\label{Kaledin}
In the rest of this section,
we assume that the dg algebra $A$ is smooth and proper.
In the smooth proper case, a conjecture of Kontsevich and Soibelman
predicts that
the Hodge-to-de Rham spectral sequences degenerate at $E_1$-stage
\cite{KS}.
On the basis of his ingenious generalization
of the Deligne-Illusie approach by positive characteristic methods
to a noncommutative setting,
Kaledin proved
this degeneration conjecture in \cite{Kal0} under some technical condition,
and in \cite{Kal2} without any assumption:

\begin{Theorem}[\cite{Kal0}, \cite{Kal}, \cite{Kal2}]
(Suppose that the dg algebra $A$ is smooth and proper.)
Then the Hodge-to-de Rham spectral sequences
\[
HH_*((t))\Rightarrow HP_{*}(A), \ \  \ HH_*[t^{-1}]\Rightarrow HC_{*}(A)
\]
degenerate at $E_1$-stage.
\end{Theorem}

The filtration on $HP_n(A)$ is defined
by $F^iHP_n(A)=\textup{Image}(H_n(F^iC_\bullet((t)))\to H_n(C_\bullet((t)))$.
Filtrations $F^iHN_n(A)$ and $F^iHC_n(A)$
are defined in a similar way.
By the degeneration, identifying $\textup{Gr}_F^iHP_n(A)$ with
$HH_{n+2i}(A)\cdot t^i$ we obtain a non-canonical isomorphism of vector spaces
\[
\oplus_{i}HH_{n+2i}(A)\cdot t^i\simeq HP_n(A).
\]
We here forget the degree of $t^i$, and hope that no confusion is
likely to arise.
Likewise, we have a non-canonical
isomorphism
$\oplus_{i\le 0}HH_{n+2i}(A)\cdot t^i\simeq HC_n(A).$
By the long exact sequence 
\[
\cdots \to HN_n(A)\to HP_n(A)\to HC_{n-2}(A)\to \cdots
\]
and reason of dimension, we also have
$\oplus_{i\ge 0}HH_{n+2i}(A)\cdot t^i\simeq HN_n(A)$
where we identify $HH_{n+2i}(A)\cdot t^i$ with 
$\textup{Gr}_F^iHN_n(A)$ (i.e., for dimension resasons, the spectral sequence
for $HN_*(A)$ degenerates).
The filtration $\{F^iHP_*(A)\}_{i\in \ZZ}$
may be considered as a noncoomutative analogue of
Hodge filtration on the de Rham cohomology
of a smooth projective variety (Remark~\ref{classHodge}). The $0$-th part $F^0HP_*(A)$
can be identified with the image of $HN_*(A)\hookrightarrow HP_*(A)$.

\begin{Remark}
\label{classHodge}
To illustrate
the analogy with the classical situation,
let us recall the comparison results.
For simplicity, $X$ is a smooth projective variety over the complex
number field, though the resuts hold more generally.
By Hochschild-Kostant-Rosenberg theorem, there is an
isomorpism
\[
HN_n(X)\simeq \oplus_{i\in \ZZ}\mathbb{H}^{2i-n}(X,\Omega_{X}^{\ge i}),\ \ \textup{and}\ \ HP_n(X)\simeq \oplus_{i\in \ZZ}H_{dR}^{2i-n}(X)=\oplus_{i\in \ZZ}\mathbb{H}^{2i-n}(X,\Omega_{X}^\bullet).
\]
where $\Omega_X^\bullet$ is the algebraic de Rham complex, $\mathbb{H}$ indicates
the hypercohomology, and 
$NH_n(X)$ and $HP_n(X)$ are the negative cyclic homology
and the periodic cyclic homology of $X$ respectively.
The classical Hodge theory implies that
$\oplus_{i\in \ZZ}\mathbb{H}^{2i-n}(X,\Omega_{X}^{\ge i})$
is equal to $\oplus_{i\in \ZZ}F^iH_{dR}^{2i-n}(X)$
where $F^i$ is the Hodge filtration.

\end{Remark}

\subsection{}

Let us recall the smoothness of functors.

\begin{Definition}
Let $\mathsf{X},\mathsf{Y}:\Art\to \Set$ be functors.
Let $F:\mathsf{X}\to \mathsf{Y}$ be a morphism.
We say that $F$ is formally smooth
if for any surjective morphism $R'\to R$ in
$\Art$, the following natural map is surjective:
\[
\mathsf{X}(R')\to \mathsf{X}(R)\times_{\mathsf{Y}(R)}\mathsf{Y}(R')
\]
induced by the morphisms $\mathsf{X}(R')\to \mathsf{X}(R)$
and $\mathsf{X}(R')\to \mathsf{Y}(R')$.

Let $*$ is the functor $\Art\to \Set$ such that
  $*(R)$ is the set consisting of
 a single element for any $R$.
 A functor $\mathsf{X}$ is said to be formally smooth
if the natural morphism $\mathsf{X}\to *$ is formally smooth.

\end{Definition}

\begin{Theorem}
\label{smoothGr}
The functor $\mathsf{Gr}$ is formally smooth.
\end{Theorem}

Let $HH_*(A)((t))$ be the dg $k[t^\pm]$-module with zero differential
whose term of cohomological degree $-n$ is the vector
space $\oplus_{i\in \ZZ}HH_{n+2i}(A)\cdot t^i$
(we abuse notation by forgetting the degree of $t^i$).
Let $HH_*(A)[[t]]$ be the dg $k[t]$-module with zero differential
whose term of
cohomological degree $-n$ is $\oplus_{i\ge 0}HH_{n+2i}(A)\cdot t^i$.
In other words, $HH_*(A)((t))=HH_*(A)\otimes_kk[t^\pm]$
and $HH_*(A)[[t]]=HH_*(A)\otimes_kk[t]$
since $HH_n(A)=0$ for $|n|>>0$.

\vspace{2mm}

We need some technical Lemmata.

\begin{Lemma}
There is an injective quasi-isomorphism of dg $k[t^\pm]$-modules
\[
HH_*(A)((t)) \to C_\bullet((t)).
\]
Similarly, 
there is an injective quasi-isomorphism of dg $k[t]$-modules
\[
HH_*(A)[[t]] \to C_\bullet[[t]].
\]
Moreover, one can choose a quasi-isomorphism in such a way that
$HH_*(A)((t)) \to C_\bullet((t))$ is obtained 
from $HH_*(A)[[t]] \to  C_\bullet[[t]]$ by
tensoring with $k[t^\pm]$.
\end{Lemma}

\Proof
We construct a quasi-isomorphism $HH_*(A)((t)) \to C_\bullet((t))$.
Suppose that $HH_m(A)=0$ for $m>N$, and $HH_N(A)\neq 0$ (if $HH_*(A)=0$,
the assertion is obvious).
We will construct an injective linear map $\oplus_{r\in \ZZ} HH_{2r}(A)\to \prod_{r}C_{2r}(A)\cdot t^r=C_\bullet((t))_0$.
Let $h$ is the maximal integer such that $2h\le N$. Choose a section
\[
q_h:HH_{2h}(A)\cdot t^h \simeq H_0(F^hC_{\bullet}((t)))\to \prod_{r\ge h}C_{2r}(A)\cdot t^r=F^hC_{\bullet}((t))_0
\]
where (the existence of) the first isomorphism follows
from the degeneration and
the vanishing of the higher term of Hochschild homology.
Next consider the long exact sequence
\[
\to HH_{2h}(A)\cdot t^h \simeq H_0(F^hC_\bullet((t)))\hookrightarrow H_0(F^{h-1}C_\bullet((t)))\twoheadrightarrow H_0(F^{h-1}/F^{h})\simeq HH_{2h-2}(A)\cdot t^{h-1}\to
\]
arising from $0 \to F^hC_\bullet((t))\to F^{h-1}C_\bullet((t))\to F^{h-1}/F^{h}\to 0$
(in fact, it is a short exact sequence by degeneration).
Choose sections $HH_{2h-2}(A)\cdot t^{h-1} \to H_0(F^{h-1}C_\bullet((t)))$
and $H_0(F^{h-1}C_\bullet((t)))\to F^{h-1}C_\bullet((t))_0$
that extends $H_0(F^hC_{\bullet}((t)))\to F^hC_{\bullet}((t))_0$.
Let $q_{h-1}:HH_{2h-2}(A)\cdot t^{h-1} \to F^{h-1}C_\bullet((t))$ be the composite.
Then we have 
\[
q_h\oplus q_{h-1}:HH_{2h}(A)\cdot t^h\oplus HH_{2h-2}(A)\cdot t^{h-1} \to F^{h-1}C_\bullet((t))_0.
\]
By the construction it is injective.
We repeat this procedure to obtain an injective linear map
\[
Q_0:=\oplus_{h\ge i}q_i: \oplus_{h\ge i} HH_{2i}(A)\cdot t^i \to C_\bullet((t))_0.
\]
%Note that this procedure terminates since $HH_j(A)=0$ for $j<<0$.
Moreover, the image of $Q_0$ is contained in the 
kernel of the differential of $C_\bullet((t))$, and the image of each
$q_{i}$
is contained in $F^iC_\bullet((t))$.
For any $z\in \ZZ$, we define
\[
Q_{2z}=Q_0 \cdot t^{-z}: \oplus_{h\ge i} HH_{2i}(A)\cdot t^{i-z} \to C_\bullet((t))_{2z}=t^{-z}\cdot C_\bullet((t))_{0}.
\]
Consider odd degrees. As in the case of degree $0$,
we construct an injective linear map
\[
P_1=\oplus_{N\ge 2i+1}p_i: \oplus_{N\ge 2i+1} HH_{2i+1}(A)\cdot t^i \to C_\bullet((t))_{1}
\]
such that
the image of $P_1$ is contained in the 
kernel of the differential of $C_\bullet((t))$, and the image of each
$p_i:HH_{2i+1}(A)t^i\to C_\bullet((t))_{1}$
is contained in $F^iC_\bullet((t))$.
For an arbitrary odd degree $2z+1$ we put
\[
P_{2z+1}=P_1\cdot t^{-z}: \oplus_{N\ge 2i+1} HH_{2i+1}(A)\cdot t^{i-z} \to C_\bullet((t))_{2z+1}=t^{-z}\cdot C_\bullet((t))_{1}.
\]
Since $\oplus_{i}HH_{n+2i}(A)\cdot t^i\simeq HP_n(A)$, these linear maps $\{Q_{2z}, P_{2z+1}\}_{z\in \ZZ}$ define an injective quasi-isomorphism
$HH_*(A)((t)) \to C_\bullet((t))$. Note also that
it is a dg $k[t^\pm]$-module map.
The case of $C_\bullet[[t]]$ is similar.
We define $Q_{2z}':\oplus_{h\ge i,\atop 0\le i-z} HH_{2i}(A)\cdot t^{i-z} \to C_\bullet[[t]]_{2z}$ and $P_{2z+1}': \oplus_{N\ge 2i+1,\atop 0\le i-z} HH_{2i+1}(A)\cdot t^{i-z} \to C_\bullet[[t]]_{2z+1}$ to be the restrictions of $Q_{2z}$ and
$P_{2z+1}$.
We then obtain an injective quasi-isomorphism
$HH_*(A)[[t]] \to C_\bullet[[t]]$ which has the desired compatibility.
\QED

\begin{Lemma}
\label{injection}
There exist a dg Lie algebra $\mathbf{E}$ and 
quasi-isomorphisms of dg Lie algebras
\[
\End_{k[t^\pm]}(C_\bullet((t)))  \leftarrow  \mathbf{E} \to  \End_{k[t^\pm]}(HH_*(A)((t))).
\]
If we replace $((t))$ by $[[t]]$, then the same assertion holds.
Namely, there exist a dg Lie algebra $\mathbf{E}_{\ge}$ and 
quasi-isomorphisms of dg Lie algebras
\[
\End_{k[t]}(C_\bullet[[t]])  \leftarrow  \mathbf{E}_\ge \to  \End_{k[t]}(HH_*(A)[[t]]).
\]
Moreover, there is a morphism of dg Lie algebras $\mathbf{E}_\ge\to \mathbf{E}$ such that the diagram
\[
\xymatrix{
\End_{k[t]}(C_\bullet[[t]]) \ar[d]  & \mathbf{E}_\ge \ar[d] \ar[r] \ar[l]   &  \End_{k[t]}(HH_*(A)[[t]]) \ar[d] \\
\End_{k[t^\pm]}(C_\bullet((t)))   &  \mathbf{E}  \ar[r] \ar[l]  &  \End_{k(((t))}(HH_*(A)((t)))
}
\]
commutes
where the right and left vertical arrows are induced  by tensoring with
$k[t^\pm]$.
\end{Lemma}

\Proof
Let us define $\mathbf{E}$ to be the pullback in the diagram of complexes
\[
\xymatrix{
\mathbf{E} \ar[r] \ar[d] &  \End_{k[t^\pm]}(HH_*(A)((t))) \ar[d] \\
\End_{k[t^\pm]}(C_\bullet((t))) \ar[r] & \Hom_{k[t^\pm]}(HH_*(A)((t)),C_\bullet((t)))
}
\]
where the lower horizontal map and the right vertical map are
induced by the composition with the quasi-isomorphism
$i:HH_*(A)((t)) \hookrightarrow C_\bullet((t))$
in Lemma~\ref{injection}.
Note that $\mathbf{E}$ is a dg Lie subalgebra
of $\End_{k[t^\pm]}(C_\bullet((t)))$
which consists of linear maps preserving $HH_*(A)((t))$.
The upper horizontal arrow and the left vertical
arrow are morphisms of dg Lie algebras.
Since $HH_*(A)((t))$ is cofibrant, the right vertical arrow is a
quasi-isomorphism.
(we here employ the complicial model category with $h$-model structure
on the category of dg $k[t^\pm]$-modules (see \cite[Theorem 3.5]{BMR}),
where a morphism is a $h$-weak equivalence if the map of dg $k[t^\pm]$-modules
is a homotopy equivalence. Every object is $h$-cofibrant and $h$-fibrant.)
Observe that the injective map $i:HH_*(A)((t)) \hookrightarrow C_\bullet((t))$ is
a trivial cofibration.
It follows that the lower horizontal arrow is a trivial fibration,
and thus the left vertical arrow and the upper horizontal arrow
are quasi-isomorphisms.
To see that $i$ is a trivial cofibration,
by \cite[Proposition 3.7]{BMR} it is enough
to show that $HH_*(A)((t)) \hookrightarrow C_\bullet((t))$
is isomorphic to $HH_*(A)((t)) \to HH_*(A)((t))\oplus C_\bullet((t))/HH_*(A)((t))$, and $C_\bullet((t))/HH_*(A)((t))$ is contractible.
For this purpose,
take decompositions $C_0=d C_1 \oplus H_0\oplus N_0$ and
$C_1=dC_2\oplus H_1\oplus N_1$ where $C_l$ (resp. $H_l$)
denotes the homologically $l$-th degree
of $C_\bullet((t))$ (resp. $HH_*(A)((t))$), and $d$ is the differential.
We here choose subspaces $N_0$ and $N_1$.
By $C_{2z}=t^{-z}\cdot C_0$ and $C_{2z+1}=t^{-z}\cdot C_1$ for $z\in \ZZ$,
we put $C_{2z}=t^{-z}\cdot (d C_1 \oplus H_0\oplus N_0)$ and
$C_{2z+1}=t^{-z}\cdot (dC_2\oplus H_1\oplus N_1)$.
Then $t^{-z}\cdot H_0=H_{2z}\to t^{-z}\cdot (H_0\oplus (dC_1\oplus N_0))$
and 
$t^{-z}\cdot H_1=H_{2z+1}\hookrightarrow t^{-z}\cdot (H_1\oplus (dC_2\oplus N_1))$
determines an injective homotopy equivalence with a splitting
$C_\bullet((t))\to HH_*(A)((t))$.
The case of $C_\bullet[[t]]$ is a variant of the above proof
($\mathbf{E}_\ge$ is the dg Lie subalgebra of $\End_{k[t]}(C_\bullet[[t]])$ which consists of linear maps preserving $HH_*(A)[[t]]$).
\QED

Theorem~\ref{smoothGr} is a consequence of the following Proposition:

\begin{Proposition}
\label{quasi-abelian}
The dg Lie algebra $\mathbb{F}$
is quasi-isomorphic
to an abelian dg Lie algebra. Here by abelian we mean the vanishing of the bracket.
\end{Proposition}

\Proof
By Lemma~\ref{injection} we
may replace $\iota:\End_{k[t]}(C_{\bullet}[[t]])\to \End_{k[t^\pm]}(C_\bullet((t)))$ by the natural injective morphism
$\iota':\End_{k[t]}(HH_*(A)[[t]])\to \End_{k[t^\pm]}(HH_*(A)((t)))$.
Then according to \cite[Proposition 3.4]{IM}
a homotopy fiber of $\iota'$ is quasi-isomorphic (as dg Lie algebras or $L_\infty$-algebras)
to an abelian dg Lie algebra.
\QED

{\it Proof of Theorem~\ref{smoothGr}.}
By Proposition~\ref{quasi-abelian} we may replace $\mathbb{F}$ by an abelian dg Lie algebra.
Thus we will assume that $\mathbb{F}$ is abelian.
Then for $R\in \Art$,
$\MC(\mathbb{F}\otimes m_R)=Z^1(\mathbb{F})\otimes m_R$.
Here $Z^1(-)$ is the space of closed elements of degree one.
Therefore, for any surjective homomorphism $R'\to R$ in $\Art$,
we see that $\MC(\mathbb{F}\otimes m_{R'})\to \MC(\mathbb{F}\otimes m_R)$
is surjective.
\QED

\subsection{}
\label{concludeSP}

We conclude this Section by some observations
which reveals a simple structure of $\mathbb{F}$ and $\mathsf{Gr}$.
According to Lemma~\ref{injection},
a homotopy fiber of the natural injective morphism
\[
\End_{k[t]}(HH_*(A)[[t]])\to \End_{k[t^\pm]}(HH_*(A)((t)))
\]
is equivalent (quasi-isomorphic) to the homotopy fiber
$\mathbb{F}$ of 
$\End_{k[t]}(C_\bullet[[t]])\to \End_{k[t^\pm]}(C_\bullet((t)))$.
The underlying complex of
a homotopy fiber of
$\End_{k[t]}(HH_*(A)[[t]])\to \End_{k[t^\pm]}(HH_*(A)((t)))$
is quasi-isomorphic to the mapping cocone
\[
\End_{k[t]}(HH_*(A)[[t]])\oplus \End_{k[t^\pm]}(HH_*(A)((t))[-1]
\]
with differential $(a,b)\mapsto (0,a)$. There is a natural
projection to the graded vector space with the zero differential:
\[
\End_{k[t^\pm]}(HH_*(A)((t)))/\End_{k[t]}(HH_*(A)[[t]])[-1],
\]
which is a quasi-isomorphism.
In addition, by Proposition~\ref{quasi-abelian}
the homotopy fiber is equivalent to
an abelian dg Lie algebra.
Consequently, we can choose the homotopy fiber to be the graded vector space
$G:=\End_{k[t^\pm]}(HH_*(A)((t)))/\End_{k[t]}(HH_*(A)[[t]])[-1]$
endowed with the zero differential and the zero bracket.
We take an $L_\infty$-morphism
\[
C^\bullet(A)[1]\to G
\]
which is equivalent to the period map $\mathcal{P}:C^\bullet(A)[1] \to \mathbb{F}$
(via an $L_\infty$-quasi-isomorphism $\mathbb{F}\to G$).
By Remark~\ref{contractionremark} the induced map
\[
HH^*(A)[1]\to \End_{k[t^\pm]}(HH_*(A)((t)))/\End_{k[t]}(HH_*(A)[[t]])[-1]
\]
carries $P$ to $H_*(\frac{1}{t}I_P((t)))\ \textup{modulo}\ \End_{k[t]}(HH_*(A)[[t]])$.
Therefore, we can summarize this observation:

\begin{Proposition}
\label{niceperiod}
There is an $L_\infty$-morphism
\[
C^\bullet(A)[1]\to G=\End_{k[t^\pm]}(HH_*(A)((t)))/\End_{k[t]}(HH_*(A)[[t]])[-1]
\]
which is equivalent to
the period map $\mathcal{P}:C^\bullet(A)[1]\to \mathbb{F}$
constructed in Section~\ref{const}.
The induced morphism
\[
HH^*(A)[1]\to G=\End_{k[t^\pm]}(HH_*(A)((t)))/\End_{k[t]}(HH_*(A)[[t]])[-1]
\]
carries $P$ in $HH^*(A)[1]$ to $H^*(\frac{1}{t}I_P((t)))$.
\end{Proposition}

\begin{Remark}
\label{easyrem}
The graded vector space
$\End_{k[t^\pm]}(HH_*(A)((t)))/\End_{k[t]}(HH_*(A)[[t]])$
is isomorphic to
\[
\oplus_{i\in \ZZ,j\in \ZZ, r<0}\Hom_k(HH_i(A), HH_j(A)\cdot t^r).
\]
where $\Hom_k(-,-)$ indicates the space of $k$-linear maps,
and the (cohomological) degree of elements in $\Hom_k(HH_i(A), HH_j(A)\cdot t^r)$
is $2r-j+i$.
Notice that an element of the image of 
$HH^*(A)[1]\to \End_{k[t^\pm]}(HH_*(A)((t)))/\End_{k[t]}(HH_*(A)[[t]])[-1]$
does not preserve the filtration on $HH_*(A)((t))\simeq HP_*(A)$.
However, $H^*(\frac{1}{t}I_P((t)))$
carries $F^{i}HP_*(A)$ to $F^{i-1}HP_*(A)$.
It can be thought of as
the Griffiths' transversality in the noncommutative situation.
\end{Remark}

\section{Unobstructedness of deformations of algebras}
\label{BTTsection}

We apply our period mapping to study unobstructedness of deformations,
and quasi-abelian properties of Hochschild cochain complexes.
We prove a noncommutative generalization of Bogomolov-Tian-Todolov theorem
(cf. Corollary~\ref{BTT}).
The argument is based on the idea of a purely
algebraic proof of Bogomolov-Tian-Todorov theorem by
Iacono and Manetti \cite{IM}.
We continue to assume that $A$ is smooth and proper
(so that by the theorem of Kaledin, Hodge-to-de Rham degenerate).

\begin{Theorem}
\label{generalBTT}
Suppose that the following condition: the linear map
\[
HH^s(A)\to \oplus_{i\in \ZZ}\Hom_{k}(HH_i(A),HH_{i-s}(A))
\]
given by $P\mapsto I_P$
is injective for any integer $s$.
We here abuse notation by writing $I_P$ for $H^*(I_P)$.
Then $C^\bullet(A)[1]$ is quasi-abelian,
 namely,
 it is quasi-isomorphic to an abelian dg Lie algebra.
 In particular, the functor $\Def_{C^\bullet(A)[1]}\simeq \mathsf{DAlg}_{A}$
 is formally smooth.
\end{Theorem}

\Proof
Note first that
the linear map
\[
H^*(\mathcal{P}):HH^*(A)[1]\to G=\End_{k[t^\pm]}(HH_*(A)((t)))/\End_{k[t]}(HH_*(A)[[t]])[-1]
\]
induced by the period mapping (see Proposition~\ref{niceperiod})
can be naturally identified with the linear map
$HH^*(A)\to \oplus_{i\in \ZZ,j\in \ZZ, r<0}\Hom_k(HH_i(A), HH_j(A)\cdot t^r)$
which sends $P$ to $\frac{1}{t}I_P$
(see Remark~\ref{easyrem}).
Thus, our condition amounts to the injectivity of the first linear graded
map
$H^*(\mathcal{P})$.
Clearly, this graded map admits a left inverse (graded) map
$G\to HH^*(A)[1]$.
If we equip $HH^*(A)[1]$ with the zero bracket,
then this inverse is a morphism of dg Lie algebras.
It gives rise to
\[
C^\bullet(A)[1]\to \End_{k[t^\pm]}(HH_*(A)((t)))/\End_{k[t]}(HH_*(A)[[t]])[-1]\to HH^*(A)[1]
\]
where the left morphism is the period mapping in Proposition~\ref{niceperiod}
which is an $L_\infty$-morphism.
The composite is an $L_\infty$-quasi-isomorphism.
Thus $C^\bullet(A)[1]$ is quasi-isomorphic to
$HH^*(A)[1]$ with the zero bracket.
The final assertion follows from the same argument with the proof of Theorem~\ref{smoothGr}.
\QED

Let us recall Calabi-Yau condition
on $A$.
We say that a smooth dg algebra $A$ is Calabi-Yau of dimension $d$
if
there is an isomorphism $f:A\to \textup{RHom}_{A^e}(A,A\otimes A)[-d]=:A^{!}[-d]$
in the triangulated category of dg $A$-bimodules (i.e. left dg $A^e:=A\otimes A^{op}$-modules), such that $f=f^{!}[d]$ (see \cite{Gin}).
Here $A$ has the $A$-bimodule structure given by the left and right multiplications,
$A\otimes A$ is endowed with the outer $A$-bimodule action of $A\otimes A$, i.e. $a\cdot (a_1\otimes a_2)\cdot  a':=aa_1\otimes a_2a'$,
and $\textup{RHom}_{A^e}(A,A\otimes A)$ is a derived Hom complex.
We use the projective model structure of left/right dg $A\otimes A^{op}$-modules. 
The $A$-bimodule structure on $\textup{RHom}_{A^e}(A,A\otimes A)$
is given by the inner bi-module action, i.e.
$a\cdot (a_1\otimes a_2)\cdot  a':=a_1a'\otimes aa_2$.

\begin{Corollary}
\label{BTT}
Let $A$ be a Calabi-Yau algebra of dimension of $d$.
Then $C^\bullet(A)[1]$ is quasi-abelian.
\end{Corollary}

\Proof
It is enough to prove that the Calabi-Yau condition
implies the hypothesis in Theorem~\ref{generalBTT}.
Assume that $A$ is Calabi-Yau of dimension $d$.
Then by Van den Bergh duality \cite[Theorem 1]{VdB}
there is an element $\pi \in HH_d(A)$
such that
\[
HH^s(A)\to HH_{d-s}(A),
\]
which carries $P$ to $H_*(I_P(\pi))$, is an isomorphism for $s\in \ZZ$.
Indeed, this isomorphism is given by
\begin{eqnarray*}
HH^s(A) & = & H^s(\textup{RHom}_{A^e}(A,A)) \\
        & \simeq & H^s(\textup{RHom}_{A^e}(A,A\otimes A)\otimes^L_{A^e} A) \\
      & \simeq & H^s(A[-d]\otimes^L_{A^e} A) \\
      & \simeq & H^{s-d}(A\otimes^L_{A^e} A) \\
     & = & HH_{d-s}(A),
\end{eqnarray*}
where the the third identification is induced by a fixed morphism
$f:A\to A^{!}[-d]$, and $\otimes^L$ is the derived tensor product.
Let $\pi$ be the image of $1_A\in HH^0(A)$ under $HH^0(A)\simeq HH_d(A)$.
(It is natural to think that $\pi$ is the fundamental class of $A$.)
As observed in  \cite[Theorem 3.4.3 (i)]{Gin} and its proof,
the above isomorphism carries a class $P$ in $HH^s(A)$ to $H_*(I_P(\pi))$
in $HH_{d-s}(A)$.
Thus, the Calabi-Yau condition implies that the hypothesis in Theorem~\ref{generalBTT}.
\QED

\begin{Example}
Let $X$ be a smooth projective Calabi-Yau variety over $k$.
Here Calabi-Yau means that the canonical bundle is trivial.
There is a dg algebra $A$
such that the dg category (or $\infty$-category) of dg $A$-modules are equivalent to
that of (unbounded) quasi-coherent complexes on $X$ (see Remark~\ref{comparisonscheme}).
Then $A$ is an example of a Calabi-Yau algebra.
In this case, it seems that Corollary~\ref{BTT}
may be proved by using Kontsevich formality type result
for the Hochschild cochain
complex $C^\bullet(X)[1]$ of $X$
and an analytic method by $\partial\bar{\partial}$-lemma (over the complex number).
However, our proof is purely algebraic.
Moreover, it is also a  ``non-commutative proof'' 
of the noncommutative problem in the sense that
the dg Lie algebra
$C^\bullet(A)[1]$ depends only on the derived
Morita equivalence class of $A$,
and
the proof does not rely on the commutative world.
\end{Example}

\begin{Example}
There are interesting constructions of Calabi-Yau algebras
due to the work of Kuznetsov (see \cite{KuzICM} and references therein).
We consider the famous example that comes from a smooth cubic fourfold $X\subset \mathbb{P}^5$. Let $\textup{Perf}(X)$ be the triangulated
category (or the enhancement by a dg category) of perfect comlexes on
$X$, that is, the derived category $D^b(X)$.
Let
\[
\mathcal{A}_X=\{F\in \textup{Perf}(X)|\ \textup{RHom}_{\textup{Perf}(X)}(\mathcal{O}_X(i),F)\simeq 0\ \ \textup{for}\ i=0,1,2\}
\]
be the orthogonal subcategory to the set of line
bundles $\OO_X$, $\OO_X(1)$, $\OO_X(2)$.
Then $\mathcal{A}_X$ is a smooth and proper 2-dimensional
Calabi-Yau category, and thus $\mathcal{A}_X$ is equivalent
to the triangulated category (or the dg category) of perfect dg modules
over a (smooth and proper) 2-dimensional Calabi-Yau algebra $A$
(see Remark~\ref{moritaremark} for perfect modules).
%We here mean by a perfect dg $A$-module that
%it belongs to the smallest triangulated
%subcategory of the category of dg $A$-modules which includes
%$A$ and is closed under direct summands.
The Hochschild homology coincides with that of K3 surfaces,
i.e., $\dim HH_0(A)=22$, $\dim HH_{-2}(A)=\dim HH_2(A)=1$ and the other terms
are zero. Moreover, there is a smooth cubic fourfold $X$ which
does not admit a K3 surface $S$ such that $\mathcal{A}_X\simeq \textup{Perf}(S)$. These categories/algebras are called ``noncommutative K3 surfaces''.
The fascinating other examples can be found in \cite{KuzICM}, \cite{Kuz1}.
We refer the reader also to \cite{Orl} for the various
constructions and facts on smooth and proper algebras.

\end{Example}

\section{Infinitesimal Torelli theorem}

The purpose of this Section is to prove the following
infinitesimal Torelli theorem:

\begin{Theorem}
Let $\mathsf{P}:\mathsf{DAlg}_A \to \mathsf{Gr}$
be the period mapping constructed in Section~\ref{peri}.
Suppose that $A$ is Calabi-Yau of dimension $d$.
Then $\mathsf{P}$ is a monomorphism.
Namely, for each $R\in \Art$ the induced map
\[
\mathsf{P}(R):\mathsf{DAlg}_A(R) \to \mathsf{Gr}(R)
\]
sending the isomorphism class of a deformation
$\tilde{A}$ of $A$ to $R$ to the isomorphism class of
the associated periodically trivialized deformation
\[
(C_\bullet(\tilde{A})[[t]],(C_\bullet(A)\otimes R)((t))^{tr}\simeq C_\bullet(\tilde{A})[[t]]\otimes_{R[t]}R[t^\pm])
\]
is injective. 
\end{Theorem}

\Proof
We first consider the case when $R=k[\epsilon]/(\epsilon^2)$.
According to Theorem~\ref{moduli}, we can interpret
$\mathsf{P}(R):\mathsf{DAlg}_A(R) \to \mathsf{Gr}(R)$
as $\Def_\mathcal{P}(R):\Def_{C^\bullet(A)[1]}(R)\to \Def_{\mathbb{F}}(R)$.
In addition, 
we may and will replace $\mathcal{P}$ by
\[
C^\bullet(A)[1]\to \End_{k[t^\pm]}(HH_*(A)((t)))/\End_{k[t]}(HH_*(A)[[t]])[-1]
\]
in Proposition~\ref{niceperiod}. We denote by $\mathcal{P}$ this
$L_\infty$-morphism.
Since the bracket on $C^\bullet(A)[1]\otimes m_R$ is zero,
there is the natural isomorphism $\Def_{C^\bullet(A)[1]}(R)\simeq HH^2(A)$.
By Proposition \ref{niceperiod}, $\Def_{\mathcal{P}}$
sends an element  $P$ in $HH^2(A)$
to $\frac{1}{t}I_P((t))$ in $\End_{k[t^\pm]}^0(HH_*(A)((t)))/\End_{k[t]}^0(HH_*(A)[[t]])$.
As in the proof of Theorem \ref{generalBTT} and Corollary \ref{BTT},
the Calabi-Yau condition implies that
\[
\Def_{\mathcal{P}}(R)=HH^2(A)\to \End_{k[t^\pm]}^0(HH_*(A)((t)))/\End_{k[t]}^0(HH_*(A)[[t]])
\]
is injective.
Hence $\mathsf{P}(R)$ is injective when $R=k[\epsilon]/(\epsilon^2)$.
Next we prove the general case by induction on the length of the maximal
ideal $m_R$ of $R$. To this end, let
\[
0 \to (\epsilon)/(\epsilon^2) \to R \to R'\to 0
\]
be an exact sequence where $(\epsilon)/(\epsilon^2)$ is
a nonunital square-zero $1$-dimensional $k$-algebras
which is the kernel of the surjective homomorphism $R\to R'$
of Artin local $k$-algebra.
Suppose that $\Def_{\mathcal{P}}(R')$ is injective.
It is enough to prove that $\mathsf{P}(R)\simeq \Def_{\mathcal{P}}(R)$
is injective.
To simplify the notation, put $L:=C^\bullet(A)[1]$ and 
$E:=\End_{k[t^\pm]}(HH_*(A)((t)))/\End_{k[t]}(HH_*(A)[[t]])[-1]$.
Let $\alpha,\beta$ be two elements in $\Def_{L}(R)$ such that
$\Def_{\mathcal{P}}(R)(\alpha)=\Def_\mathcal{P}(R)(\beta)$ in $\Def_{E}(R)$.
Note that $(\Def_{\mathcal{P}}(R)(\alpha),\Def_\mathcal{P}(R)(\beta))$
belongs to $\Def_{E}(R)\times_{\Def_{E}(R')}\Def_{E}(R)$.
It follows from the injectivity of $\Def_{\mathcal{P}}(R')$
that the images of $\alpha$ and $\beta$ in $\Def_L(R')$
coincide.
Namely, $(\alpha,\beta)$ lies in
$\Def_{L}(R)\times_{\Def_{L}(R')}\Def_{L}(R)$.
Note that $\Def_L$ satisfies the Schlessinger's condition ``$(H_1)$'' 
(see \cite[2.21]{Sch}, \cite[I.3.31]{I}).
In particular,
the natural map
$\Def_{L}(R\times_{R'}R)\to \Def_{L}(R)\times_{\Def_{L}(R')}\Def_{L}(R)$
is surjective. Here $R\times_{R'}R$ is not the tensor product but the fiber
product of artin local $k$-algebras.
We write $(r,\bar{r}+a\epsilon)$
for an element in $R\times_{k} k[\epsilon]/(\epsilon^2)$
where $\bar{r}$ is the image of $r$ in $k\simeq R/m_R$.
There is an isomorphism of artin local $k$-algebras
$R\times_{k} k[\epsilon]/(\epsilon^2)\simeq R\times_{R'}R$
which carries $(r,\bar{r}+a\epsilon)$ to $(r,r+a\epsilon)$.
We identify $\Def_L(R\times_{R'}R)$
with
$\Def_{L}(R\times_k k[\epsilon]/(\epsilon^2))\simeq \Def_{L}(R)\times \Def_L(k[\epsilon]/(\epsilon^2))$.
We choose an element $(\alpha, q)$ in
$\Def_{L}(R)\times \Def_L(k[\epsilon]/(\epsilon^2))$
which is a lift of $(\alpha,\beta)$.
It will suffice to prove that $q$ is zero.
Since $\Def_{\mathcal{P}}(k[\epsilon]/(\epsilon^2))$ is injective,
we are reduced to showing that $\Def_{\mathcal{P}}(k[\epsilon]/(\epsilon^2))(q)$
is zero in $\Def_{E}(k[\epsilon]/(\epsilon^2))$.
For this, notice that
$E$ has the zero bracket and the zero differential,
so that the natural map
\[
\Def_E(R)\times \Def_E(k[\epsilon]/(\epsilon^2))\simeq \Def_{E}(R\times_{R'}R)\to \Def_{E}(R)\times_{\Def_{E}(R')}\Def_{E}(R)
\]
is an isomorphism.
The composition $\Def_E(R)\times \Def_E(k[\epsilon]/(\epsilon^2))\to \Def_{E}(R)$ with the first projection (resp. the second projection)
$\Def_{E}(R)\times_{\Def_{E}(R')}\Def_{E}(R)\to \Def_E(R)$
is induced by the first projection (resp. $R\times_kk[\epsilon]/(\epsilon^2)\to R$ defined by $(r,\bar{r}+a\epsilon)\mapsto r+a\epsilon$).
Thus,
$(\Def_{\mathcal{P}}(R)(\alpha),\Def_\mathcal{P}(R)(\beta))$
corresponds to $(\Def_{\mathcal{P}}(R)(\alpha),0)$.
Taking account of the functoriality of $\Def_{\mathcal{P}}$,
we conclude that $\Def_{\mathcal{P}}(k[\epsilon]/(\epsilon^2))(q)$ is zero.
\QED

\end{document}